\documentclass[a4paper,11pt]{article}

\usepackage{color}
\usepackage[latin1]{inputenc}
\usepackage[T1]{fontenc}
\usepackage[normalem]{ulem}
\usepackage[english]{babel}
\usepackage{verbatim}
\usepackage{graphicx}
\usepackage{enumerate,enumitem}
\usepackage{amsmath,amssymb,amsfonts,amsthm,mathrsfs}
\usepackage{rotating,relsize}
\usepackage{float}
\usepackage{array}
\usepackage{MnSymbol}
\usepackage[all,cmtip]{xy}
\usepackage[hidelinks]{hyperref}
\xyoption{all}
\setlist[enumerate]{leftmargin=.7cm,label=\roman*)}
\usepackage{mathbbol}
\usepackage{authblk}

\newtheorem{theorem}{Theorem}[section]
\newtheorem{theoremA}{Theorem}

\newtheorem{corA}[theoremA]{Corollary}
\newtheorem*{theorem*}{Theorem}
\newtheorem*{cor*}{Corollary}

\newtheorem{notation}[theorem]{Notation}
\newtheorem{lemma}[theorem]{Lemma}
\newtheorem{prop}[theorem]{Proposition}
\newtheorem{cor}[theorem]{Corollary}

\theoremstyle{definition}
\newtheorem*{remark*}{Remark}
\newtheorem*{rem*}{Remark}
\newtheorem{defn}[theorem]{Definition}
\newtheorem{rem}[theorem]{Remark}
\newtheorem{example}[theorem]{Example}

\DeclareMathOperator{\id}{id}
\DeclareMathOperator{\im}{Im}

\DeclareMathOperator*{\colim}{colim}

\DeclareMathOperator{\THH}{THH}
\DeclareMathOperator{\THR}{THR}
\DeclareMathOperator{\TRR}{TRR}
\DeclareMathOperator{\TR}{TR}

\DeclareMathOperator{\res}{res}

\DeclareMathOperator{\K}{K}

\DeclareMathOperator{\Z}{\mathbb{Z}}
\DeclareMathOperator{\F}{\mathbb{F}}
\DeclareMathOperator{\proj}{proj}
\DeclareMathOperator{\tran}{tran}

\DeclareMathOperator{\cross}{cr}
\DeclareMathOperator{\Com}{Ring}
\DeclareMathOperator{\Ker}{Ker}

\DeclareMathOperator{\crr}{\mathbb{cr}}
\DeclareMathOperator{\Poly}{Ring^{poly}}
\DeclareMathOperator{\OO}{\mathcal{O}}
\DeclareMathOperator{\C}{\mathcal{C}}
\DeclareMathOperator{\Pro}{Pro}

\DeclareMathOperator{\D}{\mathcal{D}}
\DeclareMathOperator{\FreeT}{Free^{\mathcal{T}}}
\DeclareMathOperator{\A}{\mathcal{A}}
\DeclareMathOperator{\Sets}{Sets}
\DeclareMathOperator{\Mon}{Mon}

\title{Witt Vectors, Polynomial Maps, and Real Topological Hochschild Homology}
\author[1]{Emanuele Dotto}
\author[2]{Irakli Patchkoria}
\author[3]{Kristian Jonsson Moi}
\affil[1]{TeX.SX}
\affil[2]{Both on a bus}
\date{}

\begin{document}

\begin{center}\LARGE{Witt Vectors, Polynomial Maps, and Real Topological Hochschild Homology}
\end{center}

\begin{center}\Large{Emanuele Dotto, Kristian Moi, Irakli Patchkoria}
\end{center}

\vspace{.05cm}

\abstract{
We show that various flavours of Witt vectors are functorial with respect to multiplicative polynomial laws of finite degree. We then deduce that the $p$-typical Witt vectors are functorial in multiplicative polynomial  \textit{maps} of degree at most $p-1$. This extra functoriality allows us to extend the $p$-typical Witt vectors functor from commutative rings to $\Z/2$-Tambara functors, for odd primes $p$. We use these Witt vectors for Tambara functors to describe the components of the dihedral fixed-points of  the real topological Hochschild homology spectrum at odd primes.
}

\vspace{.05cm}

\tableofcontents

\newpage

\section*{Introduction}
\phantomsection\addcontentsline{toc}{section}{Introduction}

The various rings of Witt vectors have a prominent role in number theory and algebraic topology. They are defined as endofunctors on the category of commutative rings, and they provide a functorial way of passing from characteristic $p$ to characteristic zero. The prototypical example is the ring of $p$-typical Witt vectors of the field $\mathbb{F}_p$, which is isomorphic to the ring of $p$-adic integers $\Z_p$. These Witt vectors functors exhibit some fundamental extra structures, such as $\lambda$-operations, $\delta$-ring structures, and Frobenius lifts, which determine several of their universal properties (see e.g. \cite{AT}, \cite{Joyal}). In topology, Witt vectors appear in calculations related to topological cyclic homology \cite{Wittvect}, cyclic $K$-theory \cite{Almkvist1}, and in chromatic homotopy theory. Here they also exhibit extra structure, as they relate to the free Tambara functors of the cyclic groups \cite{BrunTamb}. In this paper we will provide novel additional structure on the Witt vectors, related to polynomial laws and polynomial maps.


We recall from \cite{Roby1} that a \emph{multiplicative polynomial law} $f$ from a commutative ring $A$ to a commutative ring $B$ is a collection of multiplicative maps
\[
f_R\colon A\otimes_{\Z}R\longrightarrow B\otimes_{\Z}R
\]
for every commutative ring $R$, which is natural with respect to ring homomorphisms in $R$. Every multiplicative polynomial law $f$ of finite degree $n$ has an underlying multiplicative map
\[f_{\Z}\colon A\longrightarrow B\]
which is $n$-polynomial, in the sense that its $(n+1)$-st cross-effect, or deviation, vanishes.
%
The main goal of this paper is to show that various Witt vectors functors extend from the category of commutative rings and ring homomorphisms to the category of commutative rings and polynomial laws of finite degree, or to polynomial maps.

In \S\ref{secGamma} we introduce an axiomatic framework of ``PD-functors'', to study this extended functoriality in polynomial laws. A \emph{PD-functor} is an endofunctor of the category of commutative rings
\[
F\colon \Com\longrightarrow \Com
\]
which commutes with certain limits and colimits. Examples of these functors include the Witt vectors $W_S$ for any truncation set $S\subset\mathbb{N}$, so in particular the big and $p$-typical Witt vectors, as well as their truncated versions. They also include the rational Witt vectors, the subring of the big Witt vectors of those power series with constant term one which are rational functions, which by a theorem of Almkvist \cite{Almkvist1, Almkvist2} is isomorphic to the cyclic $K$-theory ring. The following is the main result of \S\ref{secGamma}.

\begin{theoremA}
Any PD-functor $F\colon \Com\to \Com$ extends canonically to an endofunctor on the category $\Poly$ of commutative rings and multiplicative polynomial laws. For any of the Witt vectors functors $W$ listed above, this is the unique extension such that for any multiplicative polynomial law $f \colon A \to B$ the diagram
\[\xymatrix@C=70pt@R=15pt{W(A) \ar[d]_w \ar[r]^{W(f)} & W(B) \ar[d]^w \\ \prod A \ar[r]^{\prod f} & \prod B  }\]
commutes in $\Poly$, where $w$ is the ghost map of $W$ and $\prod f$ is the product polynomial law.
\end{theoremA}

The theorem for a general PD-functor is proved in \S\ref{secGammathm}. We reduce the construction to torsion-free rings by means of a resolution argument. We then use that the $n$-homogeneous polynomial laws out of a torsion-free ring $A$ are classified by the ``\emph{universal polynomial law}'' $\gamma_A=(-)^{\otimes n}\colon A\to (A^{\otimes n})^{\Sigma_n}$, which we by definition send to the map 
\[
F(\gamma_A)\colon F(A)\xrightarrow{\gamma_{F(A)}} (F(A)^{\otimes n})^{\Sigma_n}\longrightarrow F(A^{\otimes n})^{\Sigma_n}\stackrel{\cong}{\longleftarrow}F((A^{\otimes n})^{\Sigma_n}),
\] 
where the last map is an isomorphism by the axioms of a PD-functor. In fact we show that our extension of $F$ to $\Poly$ is the unique one that sends $\gamma_A$ to this map. When $A$ has torsion, the universal polynomial law has value in the divided powers $\Gamma_nA$, which motivates the name PD-functor, where PD stands for  ``puissances divis\'ees'' (that is ``divided powers'').
 In \S\ref{ghostlaw} we describe the ghost components of a polynomial law for the Witt vectors functors. 
\begin{rem*}
By a theorem of Almkvist  \cite{Almkvist1, Almkvist2}, our result shows that the cyclic $K$-group $\K^{cy}_0(A)$, defined as $K$-group of the exact category of endomorphisms of finitely generated projective $A$-modules modulo the zero endomorphisms, is functorial in multiplicative polynomial laws. It is well known that $\K^{cy}_0$ and $\K_0$, as functors from additive categories, are functorial in polynomial \emph{functors} (see \cite{BGMN} for a highly structured statement).  It is however not clear how multiplicative polynomial laws of commutative rings relate to polynomial functors on the respective module categories. We also remark that $\K_0$ is not a PD-functor (Example \ref{K0noGamma}), and therefore that our theorem does not provide this extra functoriality for $\K_0$.
\end{rem*}

In \S\ref{secpolymap} we turn our attention to $n$-\emph{polynomial maps}. These are the multiplicative maps $f\colon A\to B$ which satisfy the additive condition
\[
(\cross_{n+1}f)(a_1,\dots,a_{n+1}):=\sum_{\substack{U\subset \{1,\cdots, n+1\} }} (-1)^{n+1-|U|}f(\sum_{l\in U}a_l)=0.
\]
As remarked above polynomial laws forget to polynomial maps, but this correspondence is neither surjective nor injective in general. However, it is bijective when the target ring is $p$-local and the degrees are at most $p-1$. By combining this observation with the theorem above we prove the following, in \S\ref{secWpoly}. For any integer or infinity $1\leq m\leq \infty$, let $W_m(A;p)$ denote the ring of $p$-typical $m$-truncated Witt vectors.
\begin{corA}
The functor $W_m(-;p)$ extends to the partial category of multiplicative polynomial maps of degree at most $p-1$. That is, a multiplicative $n$-polynomial map $f\colon A\to B$ induces a multiplicative $n$-polynomial map 
\[W_m(f)\colon W_m(A;p)\longrightarrow W_m(B;p)\]
for every $n<p$, with the property that if $f\colon A\to B$ and $g\colon B\to C$ are multiplicative and $n$ and $k$-polynomial, respectively, and $nk<p$, then
$
W_m(g)\circ W_m(f)=W_m(g\circ f)$.
This extension is unique with the property that the diagram
\[\xymatrix@C=70pt@R=15pt{W_m(A;p) \ar[d]^w \ar[r]^-{W_m(f)} & W_m(B;p) \ar[d]^w
 \\
  \prod_{j=0}^{m-1} A \ar[r]^-{\prod_{j=0}^{m-1} f} & \prod_{j=0}^{m-1} B }\]
commutes. 

\end{corA}
Much like the universal polynomials for the sum and multiplication of $W_m(A;p)$ it does not seem to be possible to give an explicit description of the Witt components of the map $W_m(f)$, but there is an inductive procedure for finding them. For odd $p$ the first two components of $W_m(f)$ are
\[
W_m(f)(a_0,a_1,\dots)=\big(f(a_0),\sum_{i=1}^{p-1}(-1)^i({p\choose i}/p)f(a_{0}^p+ia_1),\dots\big)
\]
(see Example \ref{formula}).
When $f$ is a ring homomorphism, one can verify using standard binomial identities that the second component is equal to $f(a_1)$, so that this construction indeed extends the usual functoriality of $W_m(-;p)$ in ring homomorphisms.
In \S\ref{secWpoly} we also discuss how the hypotheses of this corollary are necessary. Most notably, there is no further extension of this functoriality on multiplicative polynomial maps of degree $p$. For example, the map
\[
N\colon \Z\longrightarrow \Z[x]/(x^2-px)
\]
that sends $a$ to $N(a)=a+\frac{a^p-a}{p}x$ is of degree $p$ and does not induce a map on $W_2(-;p)$ with the ghost components as in the corollary (see Example \ref{cex}).

Our motivation for considering polynomial maps is rooted in topology. A large supply of polynomial maps is provided by the multiplicative transfers, or norms, of \emph{Tambara functors}. 
A Tambara functor is a 
Mackey functor with a multiplication (a Green functor) and multiplicative transfers subject to certain axioms. Tambara functors naturally occur in topology, in particular in equivariant stable homotopy theory, as the components of genuine $G$-equivariant commutative ring spectra. For example the map $N$ above is the norm of the Burnside Tambara functor for the group $G=C_p$ which corresponds to the inclusion $e\to C_p$ of index $p$. This Tambara functor is the components of the initial $G$-equivariant commutative ring spectrum, namely the sphere spectrum.
For the group $G=\Z/2$, a Tambara functor $T$ consists of two  commutative rings $A$ and $B$, and maps \vspace{-.2cm}
\[
\xymatrix@C=70pt{A \ar@<1ex>[r]^-{\tran}\ar@<-1ex>[r]_-{N}\ar@(ul,dl)[]_{\tau}&B\ar[l]|-{\res},
}
\]
subject to the axioms of \cite{Tambara}. In particular the involution $\tau$ and $\res$ are ring homomorphisms, $N$ is multiplicative $2$-polynomial, and $\tran$ is additive and determined by the Tambara reciprocity relation $\tran(a)=N(a+1)-N(a)-1$.
Therefore for every odd prime $p$, we can define a diagram
\[
W_m(T;p):=\big(\xymatrix@C=70pt{\ar@(ul,dl)[]_{W_m(\tau)}&\hspace{-3.5cm} W_m(A;p) 
&W_m(B;p) \ar[ll(.71)]|-{W_m(\res)} \ar@{<-}@<-1ex>[ll(.71)]_-{\tran}\ar@{<-}@<1ex>[ll(.71)]^-{W_m(N)}\  \big),
}
\]
where $W_m(\tau)$ and $W_m(\res)$ are induced by the usual functoriality in ring homomorphisms, $W_m(N)$ is induced by the functoriality of the Corollary, and $\tran(x):=W_m(N)(x+1)-W_m(N)(x)-1$. The following is proved in \S\ref{secTambara}.

\begin{theoremA} Let $T$ be a $\Z/2$-Tambara functor, $p$ an odd prime and $1\leq m \leq\infty$ an integer or infinity. The diagram $W_m(T;p)$ is a $\Z/2$-Tambara functor. It is the unique Tambara functor functorial in $T$ with underlying rings $W_m(A;p)$ and $W_m(B;p)$ such that the ghost maps define a natural morphism of Tambara functors
\[
\xymatrix@C=70pt@R=18pt{W_{m}(A;p)\ar[d]_w \ar@<1ex>[r]^-{\tran}\ar@<-1ex>[r]_-{N}&W_{m}(B;p)\ar[d]^w\ar[l]|-{\res}
\\
\prod_{m}A \ar@<1ex>[r]^-{\prod\tran}\ar@<-1ex>[r]_-{\prod N}& \prod_{m}B\ar[l]|-{\prod\res}\rlap{\ .}
}
\]
\end{theoremA}
The reader should not confuse our construction with those arising in 
the theory of Witt vectors for Green functors of \cite{GreenWitt}. 

In \S\ref{secTHR} we use this theorem to describe the components of the dihedral fixed points of the \emph{real topological Hochschild homology} $\THR(E)$ of a connective commutative $\Z/2$-equivariant ring spectrum $E$. Let $D_{p^m}$ be the dihedral group of order $2p^m$. Then $\THR(E)$ is a commutative $D_{p^m}$-equivariant ring spectrum, for all $m\geq 0$, defined as the (derived) dihedral bar construction
\[
\THR(E):=B^{di}E=|[k]\longmapsto E^{\wedge k+1}|
\]
with the usual cyclic structure of $\THH(E)$, and the involution of $E^{\wedge k+1}$ defined as the indexed smash product over the $\Z/2$-set $\{0,1,\dots,k\}$ with the involution which reverses the order of $\{1,\dots,k\}$ (see \cite{THRmodels} and \cite{Amalie}). 
For every $m\geq 0$ one can define a $\Z/2$-spectrum
\[
\TRR^{m+1}(E;p):=\THR(E)^{C_{p^m}},
\]
which is a $\Z/2$-equivariant refinement of the $\TR$-spectra of \cite{BHM}.  In Corollary 5.2 of  \cite{THRmodels} we identify the $\Z/2$-Tambara functor of components $\underline{\pi}_0\THR(E)$, in the case $m=0$. In \cite[Theorem 3.3]{Wittvect}  Hesselholt and Madsen show that the components of the underlying ring spectrum 
\[\pi_0\TR^{m+1}(E;p)=\pi_0\THH(E)^{C_{p^m}}\cong W_{m+1}(\pi_0E;p)=W_{m+1}(\pi_0\THH(E);p)\]
are naturally isomorphic to the ring of $p$-typical $(m+1)$-truncated Witt vectors of $\pi_0E$.
 Here we establish a real version of this statement for odd primes.
We call a $\Z/2$-Tambara functor \emph{cohomological} if $N\res=(-)^{2}$ (in particular the underlying Mackey functor is cohomological, that is $\tran \circ \res=2$). For example the Tambara functor defined by the fixed-points of a commutative ring with involution is cohomological.

\begin{theoremA}
Let $E$ be a connective $\Z/2$-equivariant flat commutative orthogonal ring spectrum, with $\underline{\pi}_0E$ cohomological. Then for every odd prime $p$ and $m\geq 0$, there is an isomorphism of $\Z/2$-Tambara functors \vspace{-.3cm}
\[
\underline{\pi}_0\TRR^{m+1}(E;p)=\big(
\xymatrix@C=30pt{\pi_0\THH(E)^{C_{p^m}}
&\pi_0\THR(E)^{D_{p^m}} \ar[l]|-{\res} \ar@{<-}@<-1ex>[l]_-{\tran}\ar@{<-}@<1ex>[l]^-{N}}
\big) \cong W_{m+1}(\underline{\pi}_0\THR(E);p)
\]
which is natural in $E$. Here $W_{m+1}(-;p)$ is the Tambara functor of the previous theorem. In particular there is a natural ring isomorphism $\pi_0\THR(E)^{D_{p^m}}\cong W_{m+1}(\pi_0\THR(E)^{\Z/2};p)$.
\end{theoremA}

At the prime $p=2$, or if $\underline{\pi}_0E$ is not cohomological, the ring $\pi_0\THR(E)^{D_{p^m}}$ is in general not the Witt vectors of a ring.
Moreover on the algebraic side there is no reason for the norm of $\underline{\pi}_0\THR(E)$ to induce a map on $2$-typical Witt vectors, since the condition $n<p$ of Corollary B is violated. For example when $E=\mathbb{S}$ is the sphere spectrum, whose components are not cohomological, $\THR(\mathbb{S})=\mathbb{S}$ and therefore
\[
\pi_0\TRR^{m+1}(\mathbb{S};p)^{\Z/2}=\pi_0\mathbb{S}^{D_{p^m}}\cong \mathbb{A}(D_{p^m})
\]
is the Burnside ring of the dihedral group, which is not the $p$-typical $(m+1)$-truncated Witt vectors of $\pi_0 \THR(\mathbb{S})^{\Z/2}\cong \mathbb{A}(\Z/2)$ (see Example \ref{different-witt}).
The  dihedral fixed-points can still be described by a variant of the Witt vectors construction. For odd $p$, there are ``twisted ghost maps'' $\tilde{w}_j\colon \prod_{i=0}^m\pi_0\THR(A)^{\Z/2}\to \pi_0\THR(A)^{\Z/2}$, defined by the formula
\[
\tilde{w}_j(x_0,\dots,x_m):=\sum_{i=0}^j(1+\frac{(p^i-1)}{2}\tran^{\Z/2}_e(1))x_i(N_{e}^{\Z/2}\res_{e}^{\Z/2}(x_i))^{\frac{p^{j-i}-1}{2}}.
\]
When $\underline{\pi}_0E$, and therefore $\underline{\pi}_0\THR(E)$, is cohomological $\tran(1)=2$ and this is the usual ghost map $w_j$ of the ring $\pi_0\THR(E)^{\Z/2}$. In Theorem \ref{THRTamb} we describe the dihedral fixed-points in the following terms.
\begin{theoremA}
Let $E$ be a connective $\Z/2$-equivariant flat commutative orthogonal ring spectrum, and $p$ an odd prime.  There is a unique ring structure $\tilde{W}_{m+1}(\pi_0\THR(E)^{\Z/2};p)$ on the set $\prod_{i=0}^{m}\pi_0\THR(E)^{\Z/2}$ such that the maps $\tilde{w}_j$ are natural ring homomorphisms, and a natural ring isomorphism
\[\pi_0\THR(E)^{D_{p^m}}\cong  \tilde{W}_{m+1}(\pi_0\THR(E)^{\Z/2};p)\]
for every $1\leq m\leq \infty$.
\end{theoremA}

There is a similar description of $\pi_0\THR(E)^{D_{2^m}}$ for the prime $2$. The ring structure is again determined by a twisted version of the ghost maps $\tilde{w}_j$, which additionally take into account the action of the non-trivial Weyl group of $\Z/2$ in $D_{2^m}$. However, as a set, $\pi_0\THR(E)^{D_{2^m}}$ is a quotient of the product $\prod_{i=0}^m\pi_0\THR(E)^{\Z/2}$. This quotient accounts for the fact that the transfer maps $\pi_0\THR(E)^{D_{2^i}}\to \pi_0\THR(E)^{D_{2^{i+1}}}$ are not injective for the prime $2$. This situation is analogous to the Witt vectors for non-commutative rings of \cite{HesselholtncW} and \cite{HesselholtncWcorr} where the Verschiebung is generally not injective. The description of this quotient requires a choice of free resolution of $\underline{\pi}_0E$ as a $\Z/2$-Tambara functor, which is unsatisfying if one is interested only in the case where $E$ is a discrete ring with involution. This situation will be analysed in a forthcoming paper with different methods.

We note that the topological applications are independent on the rest of the paper, and the readers who are only interested in the algebraic results can safely ignore \S \ref{secTHR}.

\subsection*{Acknowledgements}

We  thank Peter Scholze for generously sharing his ideas, particularly for making us aware that in this context it is more natural to consider polynomial laws than polynomial maps. We thank the anonymous referees for their helpful comments, especially for suggesting to formulate some of the constructions in terms of Kleisli categories.
We would also like to thank Benjamin B\"ohme, Christopher Davis and Christian Wimmer for useful conversations. 

Dotto and Patchkoria were supported by the German Research Foundation Schwerpunktprogramm 1786 and the Hausdorff Centre for Mathematics at the University of Bonn. Moi was supported by the K\&A Wallenberg  Foundation. 

Patchkoria would like to thank the Isaac Newton Institute for Mathematical Sciences for support and hospitality during the programme ``Homotopy harnessing higher structures'' when work on this paper was undertaken. This work was supported by EPSRC grant number EP/R014604/1.

\section*{Notation and Conventions}

All the rings of this paper will be unital and commutative, and ring homomorphisms will not necessarily be unital. We will denote by $\Com$ the category of unital commutative rings and not necessarily unital ring homomorphisms (or non-unital ring homomorphisms). The category of commutative rings and unital ring homomorphisms will be denoted by $\Com_{1}$. The tensor product $\otimes$ defines a symmetric monoidal structure on $\Com$, which is not the coproduct. It is given by the usual tensor product of unital rings and the tensor product of non-unital maps.
%

\section{PD-functors and their functoriality in polynomial laws}\label{secGamma}

In this section we present an axiomatic framework that includes the various flavours of Witt vector functors. The notion of \emph{PD-functor}, named after the French ``puissances divis{\'e}es'', is an axiomatisation of the properties which allow extra functoriality in polynomial laws. Since divided powers govern polynomial laws, the axioms can be thought of as sufficient conditions for compatibility with divided power structures.

\subsection{Review of polynomial laws and divided powers}\label{recallpolylaw}

We begin by recalling some basic results about polynomial laws, mostly from \cite{Roby1}.

\begin{defn} Let $A$ and $B$ be abelian groups.
A \emph{polynomial law} from $A$ to $B$ is a collection of maps of sets $f_R\colon A\otimes_{\Z} R\to B\otimes_{\Z}R$ for every commutative ring $R$, which is natural with respect to unital ring homomorphisms $R\to R'$ (i.e., a natural transformation of set valued functors). Such a collection is called \emph{$n$-homogeneous} if
\[
f_R(x\cdot r)=f_R(x)\cdot r^{n}
\]
for every $x\in A\otimes_{\Z}R$, $r\in R$, and any commutative ring $R$. 
If $A$ and $B$ are rings, we call a polynomial law \emph{multiplicative} if each map $f_R$ is multiplicative (but not necessarily unital).
\end{defn}

All the results of the present section apply equally well to non-unital rings, which is the generality employed in \cite{Roby2}. We will however need the rings to have a unit for the constructions of \S\ref{secGammathm} (see in particular (\ref{F(gamma)})), so we will consider unital rings already from now.

Commutative rings and multiplicative polynomial laws form a category under the composition of natural transformations.

\begin{example} \label{first example}
\
\begin{enumerate}
\item A homomorphism of abelian groups $f\colon A\to B$ defines a canonical $1$-homogeneous polynomial law 
\[f_R:=f\otimes_{\Z} R\colon A\otimes_{\Z} R\longrightarrow B\otimes_{\Z}R,\]
which is multiplicative if $f$ is a ring homomorphism. This defines an embedding of the category of abelian groups and group homomorphisms into the category of abelian groups and polynomial laws. 
\item
Let $A$ be a ring, and $n$ a non-negative integer. The $n$-th power maps
\[
(-)^{n}\colon A\otimes_{\Z}R\longrightarrow A\otimes_{\Z}R
\]
of the rings $A\otimes_{\Z}R$ define a multiplicative $n$-homogeneous polynomial law.
\end{enumerate}
\end{example}

\begin{rem} \label{conventionpoly} To denote a polynomial law from $A$ to $B$, we will use the usual arrow $A \to B$. In view of Example \ref{first example} (i) this should not cause any confusions since homomorphisms uniquely correspond to $1$-homogenous polynomial laws \cite[Section I.11]{Roby1}. The same applies in the multiplicative context. In particular any commutative diagram of rings and homomorphisms uniquely determines a commutative diagram of polynomial laws. 
This convention will be in place in all of Section \ref{secGamma}. 
\end{rem}

We recall from \cite{Roby1} that there is a universal $n$-homogeneous polynomial law
\[
\gamma_n\colon A\longrightarrow \Gamma_n(A),
\]
where  $\Gamma_n(A)$ is the $n$-th graded piece of the divided power algebra of the abelian group $A$. This is universal in the following sense. Let $H_n(A,B)$ be the set of $n$-homogeneous polynomial laws from $A$ to $B$, and $M_n(A,B)$ the subset of the multiplicative laws.

\begin{prop}[\cite{Roby1, Roby2}]\label{classification}
For every pair $A$ and $B$ of abelian groups, there is a natural bijection
\[
\hom_{Ab}(\Gamma_{n}(A),B)\cong H_n(A,B)
\]
which sends a group homomorphism $\varphi\colon \Gamma_{n}(A)\to B$ to $\varphi\circ\gamma_n$.  If $A$ is a commutative ring, then $\Gamma_n(A)$ possesses a natural ring structure \cite{Roby2}, and restricting to ring homomorphisms gives a bijection. 
\[
\hom_{\Com}(\Gamma_{n}(A),B)\cong M_n(A,B)
\]
for every commutative ring $B$.
\end{prop}

For the convenience of the reader we recall the construction of divided powers (see e.g., \cite[III.1]{Roby1}). Let $A$ be an abelian group. The divided power algebra of $A$, denoted by $\Gamma(A)$, is a commutative ring generated by the symbols $\gamma_n(a)$ for all $n \geq 0$ and $a \in A$ subject to the following relations:
\begin{enumerate}
\item $\gamma_0(a)=1$, for all $a\in A$;
\item $\gamma_n(ka)=k^n\gamma_n(a)$, for all $a \in A$ and $k \in \Z$;
\item $\gamma_n(a)\gamma_m(a)=\frac{(n+m)!}{n!m!}\gamma_{n+m}(a)$, for all $a \in A$;
\item $\gamma_n(a+b)=\Sigma_{k=0}^n\gamma_k(a)\gamma_{n-k}(a)$, for all $a \in A$.
\end{enumerate}
Then $\Gamma_n(A)$ is the subgroup of $\Gamma(A)$ generated by elements of the form $\gamma_{i_1}(a_1)\gamma_{i_2}(a_2)\dots\gamma_{i_l}(a_l)$, where $i_1 + i_2 +\dots+ i_l=n$, and  $\Gamma(A)$ decomposes as a direct sum
\[\bigoplus_{n\geq0} \Gamma_n(A).\]
The universal $n$-homogeneous polynomial law $\gamma_n \colon A \to \Gamma_n(A)$ defined in \cite[IV.1]{Roby1}, is given by 
\[a \otimes r \mapsto \gamma_n(a) \otimes  r^n\]
for an elementary tensor $a \otimes r \in A \otimes R$. We note that for the sums of elementary tensors, the formulas in general are more complicated and involve lower divided powers. Thus homogeneous polynomial laws are generally not determined by their value at $R=\Z$. See for example \cite[Example 1.2]{Chenev}.

We recall from \cite{Roby1} that any polynomial law $f$ has a ``Taylor decomposition'': it decomposes uniquely into a locally finite sum of homogeneous polynomial laws. If there are finitely many non-zero homogeneous pieces, we say that the polynomial law is of degree $n$, where $n$ is the largest degree of its homogeneous summands.

\begin{notation}
Let $\Poly$ denote the category whose objects are commutative rings and whose morphisms are the multiplicative polynomial laws of finite degree, and where the composition is the composition of natural transformations. The category $\Com$ is a subcategory of $\Poly$ of the multiplicative polynomial laws of degree $1$ which preserve zero.
\end{notation}

Let us denote the Taylor decomposition of a polynomial law $f\colon A\to B$ of degree $n$ by
\[f=f_0+f_1+\cdots+ f_n,  \]
where $f_i$ is a homogeneous polynomial law of degree $i$ (the $f_i$ are however not unital even if $f$ is, in particular $f_0$ is a constant in $B$).
Each $f_i$ in the Taylor decomposition of $f$ corresponds to a unique additive map
\[\varphi_i \colon \Gamma_i A \longrightarrow B,\] 
such that $\varphi_i \circ \gamma_i = f_i$. When $f$ is multiplicative, the maps $\varphi_i$ have the following properties (see \cite[Page 71]{Zieplies}): 
\begin{itemize}
\item The map $\varphi_i$ is multiplicative, i.e., a not necessarily unital homomorphism of rings;
\item The maps $\varphi_i$ and $\varphi_j$ are orthogonal for $i\neq j$, i.e., for any $x \in \Gamma_i A$ and $y \in \Gamma_j A$, one has
\[ \varphi_i(x) \cdot \varphi_j(y)=0. \]
\end{itemize}
Altogether  the multiplicative polynomial law $f$ admits a unique factorisation in $\Poly$
\[\xymatrix@C=70pt@R=13pt{A \ar[r]^{f}  \ar[d]_-{\prod_i \gamma_i} & B, \\ \prod_{i=0}^{n} \Gamma_i A \ar[ur]_{\varphi=\oplus_i \varphi_i}    & } \]
and the map $\oplus_i \varphi_i$ is a morphism in the subcategory $\Com$ of ring homomorphisms (see Remark \ref{conventionpoly}). Note that if $f$ is unital, then so is $\oplus_i \varphi_i$ (but not the individual $\varphi_i$).

We now give a more conceptual description of the category $\Poly$, by expressing it as a Kleisli-type category of something that resembles a comonad on $\Com$. Let us look more closely at the structure of the sequence of endofunctors $\{\prod_{i=0}^{n} \Gamma_i\}_{n \geq 0}$ of $\Com$. We have natural ring homomorphisms
\[\pi_{n+1} \colon \prod_{i=0}^{n+1} \Gamma_i A \longrightarrow \prod_{i=0}^{n} \Gamma_i A\]
for every $n\geq 0$, defined by the projections. Additionally, for every $n,m\geq 0$, there is a natural ring homomorphism
\[\Delta_{n,m} \colon \prod_{k=0}^{nm} \Gamma_k A \longrightarrow \prod_{i=0}^{n} \Gamma_i (\prod_{j=0}^{m} \Gamma_j A)\]
which classifies the composite polynomial law
\[\xymatrix{ A \ar[r]^-{\prod_j \gamma_j} & \prod_{j=0}^{m} \Gamma_j A \ar[r]^-{\prod_i \gamma_i} &   \prod_{i=0}^{n} \Gamma_i (\prod_{j=0}^{m} \Gamma_j A) } \]
of degree $nm$ under the isomorphism of Proposition \ref{classification}. Finally, we also have the projection morphism $\varepsilon \colon \prod_{i=0}^{1} \Gamma_i (A)=\Z \times A \to A$. These morphisms are compatible in a way which resembles the axioms of a comonad, and govern the composition of polynomial laws. This point of view will be very useful in the subsequent parts of the paper. To formalise this structure we introduce the following:

\begin{defn} \label{mufiltered} Let $\C$ be a category. A $\mu$-filtered comonad $T_\bullet$ on $\C$ is a sequence of endofunctors $T_n\colon\C\to \C$, for $n\in\mathbb{N}$, and natural transformations
\[\pi_{n+1}\colon T_{n+1}\to T_n,\quad \Delta_{n,m} \colon T_{nm} \to T_n \circ T_m,\quad \varepsilon \colon T_1 \to \id\]
for all $n,m \geq 0$, satisfying the following axioms:
\begin{enumerate}
\item (Coassociativity) For any $n, m ,k \geq 0$, the diagram 
\[\xymatrix@C=85pt{T_{nmk} \ar[d]_-{\Delta_{nm, k } } \ar[r]^-{\Delta_{n, mk}} & T_n \circ T_{mk} \ar[d]^-{T_n \Delta_{m,k}} \\  T_{nm} \circ T_k \ar[r]^-{\Delta_{n,m} T_k } &  T_n \circ T_m \circ T_k}\]
commutes. 

\item (Unitality) For any $n \geq 0$, the diagram
\[\xymatrix@C=50pt{T_n \circ T_1 \ar[dr]_{T_n \varepsilon} & T_n \ar[l]_-{\Delta_{n,1}}  \ar@{=}[d] \ar[r]^-{\Delta_{1,n}} & T_1 \circ T_n  \ar[dl]^{\varepsilon T_n} \\ & T_n & }\]
commutes. 

\item (Filtration) For any $n' \geq n$ and $m' \geq m$, the diagrams
\[\xymatrix@C=35pt{T_{n'm} \ar[d]_{\pi_{n'm, nm}} \ar[r]^-{\Delta_{n',m}} & T_{n'} \circ T_m \ar[d]^{\pi_{n',n} T_m} & &   T_{nm'} \ar[d]_-{\pi_{nm', nm}} \ar[r]^-{\Delta_{n,m'}} & T_{n} \circ T_{m'} \ar[d]^{T_n \pi_{m',m}} \\ T_{nm} \ar[r]^-{\Delta_{n,m}} & T_{n} \circ T_m  & &   T_{nm}  \ar[r]^-{\Delta_{n,m}} & T_{n} \circ T_{m}}\]
commute. Here $\pi_{i,j}=\pi_{j+1} \circ \dots \circ \pi_{i} \colon T_i \to T_j$ for $i >j$.
\end{enumerate}
\end{defn}
Here $\mu$-filtered stands for multiplicatively filtered, as $T_\bullet$ is indexed over the multiplicative monoid of natural numbers.
When the maps $\pi_{n+1}$ are identities this definition recovers the usual notion of comonad. We recall that given an ordinary comonad $T$ on a category $\C$, one can define the \emph{Kleisli category} $\C_T$ with the same objects as $\C$ and the morphisms
\[\C_T(X,Y):=\C(TX, Y).\]
The compositions and identities are defined using the structure maps of $T$, and the associativity and unitality of $\C_T$ follow from those for $T$ \cite{Kleis}. The following generalises the Kleisli category to $\mu$-filtered comonads.
 
\begin{prop} \label{Kleisli} Let $T_{\bullet}$ be a $\mu$-filtered comonad on a category $\C$. Then there is a category $C_{T_{\bullet}}$, which we call the Kleisli category of $T_{\bullet}$, with the same objects as $\C$ and with morphism sets
\[\C_{T_{\bullet}}(X, Y)=\colim_n\C(T_nX, Y),\]
where the colimit is taken in the category of sets along the maps $\pi_n^*$. The composition is defined on representatives by
\[
(T_mY\xrightarrow{g} Z)\circ (T_nX\xrightarrow{f} Y)=(T_{mn}X\xrightarrow{\Delta_{m,n}}T_{m}(T_n X)\xrightarrow{T_mf}T_mY\xrightarrow{g} Z),
\]
and the identity of $X$ is $\varepsilon\colon T_1X\to X$.
\end{prop}

\begin{proof} This follows immediately from the axioms of a $\mu$-filtered comonad. In particular, the Filtration Axiom of Definition \ref{mufiltered} implies that the compositions do not depend on the choice of colimit  representatives. 
\end{proof}

Our main example of a $\mu$-filtered comonad is the sequence of functors $T_n=\prod_{i=0}^{n} \Gamma_i$ together with the canonical projections, the map $\varepsilon$, and the maps
\[\Delta_{n,m} \colon \prod_{k=0}^{nm} \Gamma_k \longrightarrow \prod_{i=0}^{n} \Gamma_i (\prod_{j=0}^{m} \Gamma_j )\]
defined above. Checking that the axioms of Definition \ref{mufiltered} are satisfied is straightforward and just uses that $\Delta_{n,m}$ classifies the composition of the universal polynomial laws of degrees $n$ and $m$. We let $\Com_{\Gamma}$ denote the Kleisli category associated to the $\mu$-filtered comonad $\prod_{i=0}^{\bullet} \Gamma_i (-)$.

\begin{prop} \label{alternativepolycat} The categories $\Com_{\Gamma}$ and $\Poly$ are equivalent.
\end{prop}

\begin{proof} We define a functor $E: \Poly \to \Com_{\Gamma}$. It is the identity on the objects. A multiplicative polynomial law $f \colon A \to B$ of finite degree $n$ uniquely factors in $\Poly$ as
\[\xymatrix@C=70pt@R=13pt{A \ar[r]^{f}  \ar[d]_-{\prod_i \gamma_i} & B. \\ \prod_{i=0}^{n} \Gamma_i A \ar[ur]_{\varphi}    & } \]
where $\varphi$ is a ring homomorphism. We let $E(f)$ be the element in $\colim _m \Com(\prod_{i=0}^{m} \Gamma_i A, B)$ determined by $\varphi$. By the surjectivity of the projections this determines a bijection between finite degree polynomial laws from $A$ to $B$ and $\Com_{\Gamma}(A,B)$. Clearly $E$ preserves the identities. It remains to show that $E$ preserves compositions. Given two composable arrows in $\Poly$,
\[A \stackrel{f}{\longrightarrow} B \stackrel{g}{\longrightarrow}  C, \]
let us consider the corresponding ring homomorphisms $\varphi \colon  \prod_{j=0}^{n} \Gamma_j A \to B$ and $\psi \colon  \prod_{i=0}^{m} \Gamma_i B \to C$. It follows immediately from the definition of $\Delta_{m,n}$ that the composite $E(g)\circ E(f)$ in the Kleisli category
\[ \xymatrix{\prod_{k=0}^{mn} \Gamma_k A \ar[r]^-{\Delta_{m,n}} &  \prod_{i=0}^{m} \Gamma_i(\prod_{j=0}^{n} \Gamma_j A) \ar[rr]^-{\prod_{i=0}^{m} \Gamma_i(\varphi)} & & \prod_{i=0}^{m} \Gamma_i B \ar[r]^-{\psi} & C }\]
is the unique ring homomorphism corresponding to the composite polynomial law $g \circ f$, that is $E(g\circ f)$.
\end{proof}

\begin{rem} The Kleisli category $\C_{T_{\bullet}}$ of a $\mu$-filtered comonad $T_{\bullet}$ can be related to the Kleisli category of an honest comonad. Indeed $T_\bullet$ defines a comonad $T$ on the category $\Pro(\C)$ of pro-objects in $\C$, whose underlying functor assigns to a pro-object $\{X_{\lambda}\}_{\lambda \in \Lambda}$  the pro-object  $\{T_n(X_{\lambda})\}_{(n, \lambda) \in \mathbb{N} \times \Lambda}$. 
Then $\C_{T_{\bullet}}$ is the full subcategory of the Kleisli category $\Pro(\C)_T$ spanned by the constant pro-objects. 
\end{rem}


%

\subsection{Definition and examples of PD-functors}

In this section we introduce the notion of \emph{PD-functor} which defines a class of well-behaved endofunctors on $\Com_{1}$. We will show in \S \ref{secGammathm} that these extend in a canonical way to functors on polynomial laws, and the definition given here is tailored to this purpose. 

We recall that $\Com$ is the category of unital rings and not necessarily unital ring homomorphisms, and $\Com_1$ is the subcategory of unital ring homomorphisms.




\begin{defn}A functor $F \colon \Com_1 \to \Com_1$ is a PD functor if it preserves the following universal constructions:
\begin{enumerate}
	\item finite products,
	\item reflexive coequalisers,
	\item fixed points of finite group actions.
\end{enumerate}
\end{defn}

\begin{example}\label{exproduct}
Let $S$ be a set, and $(-)^{\times S}\colon \Com\to \Com$ be the functor that takes a ring $A$ to the product ring $A^{\times S}=\prod_{S}A$, and a ring homomorphism to the product map $f^{\times S}$. Then $(-)^{\times S}$ is a PD-functor. Conditions i) and iii) are satisfied, since products commute with limits. Given a reflexive coequaliser of commutative rings, it is also a reflexive coequaliser of underlying abelian groups and sets. If the set $S$ is finite, then it is clear that $(-)^{\times S}$ commutes with reflexive coequalisers, since finite products in sets do. Infinite products in sets do not commute with reflexive coequalisers in general. However, infinite products of abelian groups do commute with reflexive coequalisers and reflexive coequalisers in abelian groups are reflexive coequalisers of underlying sets.

\end{example}

\begin{example}\label{tensor}
Let $R$ be a torsion-free commutative ring. The functor $(-)\otimes_{\Z}R\colon \Com_{1}\to \Com_{1}$ that sends a commutative ring $A$ to $A\otimes_{\Z}R$ is a PD-functor. The fact that $R$ has no torsion is used to show that  $(-)\otimes_{\Z}R$ commutes with invariants of finite groups.
\end{example}

\begin{example} \label{Witt vectors} For any truncation set $S$, the $S$-truncated Witt vectors functor $W_S \colon \Com_{1} \to \Com_{1}$ is a PD-functor. We recall that as a set
\[
W_S(A)=\prod_{s\in S}A,
\]
and that the functoriality is given by taking the product map on underlying sets. As in Example \ref{exproduct} this implies that $W_S$ satisfies Conditions i) and iii). A reflexive coequaliser of commutative rings is also a reflexive coequaliser of underlying abelian groups and sets. To show Condition ii) one can argue exactly as in Example \ref{exproduct} using that $W_S$ and $(-)^{\times S}$ have the same underlying set.

By choosing the appropriate $S$, we see that big Witt vectors, $n$-truncated big Witt vectors, $p$-typical Witt vectors and $n$-truncated $p$-typical Witt vectors are all PD-functors. See for example \cite{LarsWitt, Rabinoff} for more background on Witt vectors. 

\end{example}

\begin{example} \label{rationalwitt} For a commutative ring $A$, the ring of big Witt vectors $\mathbb{W}(A)$ is isomorphic to a ring whose underlying additive group is the subgroup of the units of $A[[x]]$ consisting of the power series with constant term equal to $1$. The ring of rational Witt vectors $W_{\mathbb{Q}}(A)$ is the subring of $\mathbb{W}(A)$ corresponding to those power series which are rational functions. See \cite{Almkvist1, Almkvist2, Haze} for more details. 

It is clear that $W_{\mathbb{Q}}$ preserves finite products. We verify Conditions ii) and iii) directly.
For a reflexive coequaliser of unital ring maps
\[\xymatrix{A \ar@<0.5ex>[r]^{f}  \ar@<-0.5ex>[r]_{g} & B \ar[r]^{h}  & C } \]
with a common section $s \colon B \to A$, we would like to show that
\[\xymatrix{W_{\mathbb{Q}}(A) \ar@<0.5ex>[rr]^-{W_{\mathbb{Q}}(f)}  \ar@<-0.5ex>[rr]_-{W_{\mathbb{Q}}(g)} & & W_{\mathbb{Q}}(B) \ar[rr]^-{W_{\mathbb{Q}}(h)} & & W_{\mathbb{Q}}(C) } \]
is a coequaliser. One only needs to check that $\Ker W_{\mathbb{Q}}(h) \subset \im (W_{\mathbb{Q}}(f) -W_{\mathbb{Q}}(g) )$. Suppose that we have 
\[\frac{1+b_1x + \cdots + b_n x^n}{1+b'_1x + \cdots + b'_n x^n} \in \Ker W_{\mathbb{Q}}(h).\] 
Then $h(b_i)=h(b_i')$ for all $i=1, \cdots, n$. This implies that there are $a_i \in A$, $i=1, \dots, n$ such that $f(a_i)=b_i$ and $g(a_i)=b'_i$ for every $i$. Consider the rational function
\[\frac{1+a_1x + \cdots + a_n x^n}{1+s(b_1)x + \cdots + s(b_n) x^n} \in W_{\mathbb{Q}}(A).\]
Then the map $W_{\mathbb{Q}}(f) -W_{\mathbb{Q}}(g)$ sends this rational function to
\begin{align*} &\frac{1+f(a_1)x + \cdots + f(a_n) x^n}{1+f(s(b_1))x + \cdots + f(s(b_n)) x^n} \cdot  \frac{1+g(s(b_1))x + \cdots + g(s(b_n)) x^n}{1+g(a_1)x + \cdots + g(a_n) x^n}
\\&= \frac{1+b_1x + \cdots + b_n x^n}{1+b_1x + \cdots + b_n x^n} \cdot  \frac{1+b_1x + \cdots + b_n x^n}{1+b'_1x + \cdots + b'_n x^n}
=\frac{1+b_1x + \cdots + b_n x^n}{1+b'_1x + \cdots + b'_n x^n}, \end{align*}
where we have used that $fs=gs=1$. This shows Condition ii).

Let $G$ be a finite group and $A$ a commutative ring with $G$-action. To show Condition iii) we want to show that the canonical map
\[W_{\mathbb{Q}}(A^G) \to W_{\mathbb{Q}}(A)^G\]
is an isomorphism. Since $W_{\mathbb{Q}}$ preserves injections this map is injective, so we must show surjectivity. Given a polynomial $\alpha \in A[x]$ and an element $g \in G$, we denote by $g \cdot \alpha$ the polynomial obtained by acting on the coefficients of $\alpha$ by $g$. Now let
\[\frac{\alpha}{\beta} \in  W_{\mathbb{Q}}(A)^G \]
be an invariant rational function.
Then by definition for any element $g \in G$, we have
\[g \cdot \frac{\alpha}{\beta}= \frac{g \cdot \alpha}{g \cdot \beta}= \frac{\alpha}{\beta}. \]
Hence $(g \cdot \alpha) \beta= \alpha (g \cdot \beta)$. Note that $\alpha$ and $\beta$ do not have to be $G$-invariant. However, the latter equation allows us to replace them with $G$-invariant polynomials. Indeed, we observe that
\[\frac{\alpha}{\beta}= \frac{\alpha (\prod_{g \in G, g \neq 1} g \cdot \beta ) } {\prod_{g \in G} g \cdot\beta}.\]
It is clear that the denominator of the right hand fraction is $G$-invariant. We check that the numerator is invariant as well:
\[\begin{aligned} h\cdot (\alpha (\prod_{g \in G, g \neq 1} g \cdot \beta ))= (h \cdot \alpha) (h h^{-1} \cdot \beta  )(\prod_{g \in G, g \neq 1, h^{-1}} hg \cdot \beta )= \\  (h \cdot \alpha) \beta (\prod_{g \in G, g \neq 1, h^{-1}} hg \cdot \beta ) =  \alpha (h \cdot \beta) (\prod_{g \in G, g \neq 1, h^{-1}} hg \cdot \beta )=   \alpha (\prod_{g \in G, g \neq 1} g \cdot \beta ).\end{aligned} \]
The third equality uses the relation  $(h \cdot \alpha) \beta= \alpha (h \cdot \beta)$. Now the quotient 
\[\frac{\alpha (\prod_{g \in G, g \neq 1} g \cdot \beta ) } {\prod_{g \in G} g \cdot\beta}\]
is an element of $W_{\mathbb{Q}}(A^G)$ that maps under the canonical map to $\frac{\alpha}{\beta}$ in $W_{\mathbb{Q}}(A)^G$. 
\end{example}

\begin{example}\label{K0noGamma}
By Almkvist's theorem \cite{Almkvist1, Almkvist2}, the functor $W_{\mathbb{Q}}$ is isomorphic to the cyclic $K$-theory functor $\K^{cy}_0$, which is therefore a PD-functor. We remark however that the Grothendieck group $\K_0$ is not a PD-functor. Indeed, since any surjective ring homomorphism is the projection onto a reflexive coequaliser, any PD functor must preserve surjections, whereas $\K_0$ does not.
\end{example}

We end the section by remarking that any product-preserving endofunctor of $\Com_1$ extends canonically to an endofunctor of $\Com$. For this reason we will often implicitly consider a PD-functors as being defined on the larger category $\Com$.

\begin{lemma}\label{lemma:extension} Let $F \colon \Com_{1} \to \Com_{1}$ be a functor that preserves finite products. Then there is a functor $E \colon \Com \to \Com$ and an isomorphism $\phi \colon E|_{\Com_1} \to F$. The pair $(E,\phi)$ is unique up to unique isomorphism over $F$.
\end{lemma}

\begin{proof} We first prove existence. For a possibly non-unital ring $A$ we write $A^+$ for its unitalisation and $p_A \colon A^+ \to \Z$ for the canonical projection. We set $E(A) = \ker(F(A^+) \to F(\Z))$ and note that this defines an endofunctor of the category of non-unital rings and non-unital ring homomorphisms. If $A$ is unital there is a natural unital map $A^+ \to A$ and the induced map $A^+ \to A\times \Z$ is an isomorphism of unital rings. Since $F$ preserves products, taking kernels of the projection to $F(\Z)$ gives an isomorphism $\phi_A \colon E(A) \to F(A)$, so we see that $E$ restricts to an endofunctor $F$ on $\Com$. The maps $\phi_A$ are natural in unital ring homomorphisms so we get a natural isomorphism $\phi \colon E|_{\Com_1} \to F$ as desired.

We now prove uniqueness. Let $(E,\phi)$ and $(E',\phi')$ be pairs as in the statement of the lemma and write $\psi = (\phi')^{-1} \circ \phi$. We will show that $\psi$ is natural with respect to all maps in $\Com$. Any such map $f\colon A \to B$ factors as 
\[A \stackrel{e\cdot f}{\longrightarrow} e B \xrightarrow{\id \times 0} e B \times (1-e)B \cong B,\]
where $e = f(1)$ is idempotent. The left and right hand maps are unital, so it suffices to check $\psi$ is natural with respect to maps of the form $\id \times 0\colon R \to R \times S$. This follows since $F$, and therefore also $E$ and $E'$, preserve finite products. The uniqueness of $\psi$ is clear.
\end{proof}

By abuse of notation we will simply write $F$ for a choice of extension of a product preserving functor $F$ on $\Com_1$ to non-unital ring homomorphisms. A pair of maps $f \colon A \to C$ and $g\colon B \to C$ in $\Com$ is called \emph{orthogonal} if $f(a)\cdot g(b)=0\in C$ for all $a\in A$ and $b\in B$.

\begin{lemma}\label{lemma:ext-prop}Let $F \colon \Com_1 \to \Com_1$ be a functor that preserves finite products. Then the extension to the category $\Com$ has the following properties:
\begin{enumerate}
	\item\label{lemma:ext-prop1} For any pair of maps $f \colon A \to A'$ and $g \colon B \to B'$ in $\Com$, the diagram
\[\xymatrix{F(A) \otimes F(B) \ar[d]^c \ar[rr]^-{F(f) \otimes F(g)} && F(A') \otimes F(B') \ar[d]^{c} \\   F(A \otimes B) \ar[rr]^{F(f \otimes g)} && F(A' \otimes B') } \]
commutes, where $c$ is the canonical coproduct map for unital rings.
	\item\label{lemma:ext-prop2} For any pair of orthogonal maps $f \colon A \to C$ and $g \colon B \to C$ the maps $F(f)$ and $F(g)$ are orthogonal.
\item\label{lemma:ext-prop3} For any pair of orthogonal maps $f\colon A\to C$ and $g\colon B\to C$ in $\Com$ there is a commutative diagram in $\Com$
\[
\xymatrix@C=70pt@R=17pt{
F(A\times B)\ar[r]^-{F(f+g)}\ar[d]_-{\cong}&F(C)
\\
F(A)\times F(B).\ar[ur]_-{\ F(f)+F(g)}
}
\]

\end{enumerate}
\end{lemma}

\begin{proof}Part i): The statement is true for unital maps. By factoring $f\otimes g$ as $(f\otimes \id)\circ (\id\otimes g)$ we may assume that $f$ is the identity. Moreover since every map in $\Com$ factors as a unital map followed by an inclusion $\id \times 0 \colon B \to B \times C$ (see the proof of \ref{lemma:extension}), we may assume that $g$ is of the form $\id \times 0$. We must show that the left hand square in the following diagram commutes
\[\xymatrix{F(A) \otimes F(B) \ar[d]^c \ar[rr]^-{\id \otimes F(\id\times 0)} && F(A) \otimes F(B \times C) \ar[d]^{c} \ar[rrr]_-{\cong}^-{(\id \otimes F(p_1),\id \otimes F(p_2))} &&&  (F(A)\otimes F(B)) \times (F(A)\otimes F(C)) \ar[d]^{c\times c}
 \\ 
  F(A \otimes B) \ar[rr]^-{F(\id \otimes (\id\times 0))}  && F(A \otimes (B\times C)) \ar[rrr]_-{\cong}^-{(F(\id \otimes p_1) , F(\id \otimes p_2))} &&& F(A \otimes B) \times F(A\otimes C).} \]
The right hand square commutes and the right hand horizontal maps are isomorphisms, so it suffices to show that the outer rectangle commutes. The composite through the upper right hand corner is the map $(c \circ \id\otimes \id) \times c \circ (\id \otimes 0)$ and the other composite equals $(F(\id \otimes \id) \circ c) \times (F(\id \otimes 0) \circ c)$. The second coordinates of both maps are equal to the $0$ map and the first coordinates are equal because $\otimes$ is the coproduct on unital rings.


Part ii): Orthogonality of $f$ and $g$ is equivalent to the composite map
\[A\otimes B \xrightarrow{f\otimes g} C \otimes C \xrightarrow{\mu_C}  C \]
factoring through $0$, where $\mu_C$ is the multiplication of $C$. Since $F(0)=0$, Part i) implies that 
\[F(A)\otimes F(B) \xrightarrow{F(f)\otimes F(g)} F(C) \otimes F(C) \xrightarrow{\mu_{F(C)}} F(C) \]
factors through $0$ as well.

Part iii): The orthogonal maps $f$ and $g$ give a splitting of $C$ as $C \cong C_f \times C_g \times C'$ with $f$ and $g$ factoring through $C_f$ and $C_g$, respectively. Here $C_f=f(1)C$ and $C_g=g(1)C$. The map $f+g$ is the composite
\[A \times B \xrightarrow{f^u \times g^u} C_f \times C_g \hookrightarrow C_f \times C_g \times C' \stackrel{\cong}{\longrightarrow} C,\]
where $f^u$ is the corestriction of $f$, which is unital, and similarly for $g^u$.
Now consider the diagram
\[
\xymatrix@C=70pt@R=17pt{
F(A\times B)\ar[r]^-{F(f^u \times g^u)}\ar[d]_{\cong}& F(C_f \times C_g) \ar[d]^-{\cong} \ar[r] &F(C)
\\
F(A)\times F(B)  \ar[r]^{F(f^u) \times F(g^u)} & F(C_f) \times F(C_g) \ar[ur] &
}
\]
where the vertical maps are the canonical isomorphisms. This diagram commutes in $\Com$ and the composite of the upper row is $F(f + g)$ while the composite from the lower left hand corner to the right hand corner is the map $F(f) + F(g)$.
\end{proof}

\subsection{Functoriality of PD-functors in multiplicative polynomial laws}\label{secGammathm}

We now explain how to extend a PD-functor $F\colon \Com \to \Com$ to the category $\Poly$ of commutative rings and polynomial laws of finite degree.

If $A$ is a torsion-free ring, the natural transformation $\xi_A\colon \Gamma_n(A)\to (A^{\otimes n})^{\Sigma_n}$ from \cite[Section III.6]{Roby1} is an isomorphism, and we can consider the polynomial law defined as the composite
\begin{equation}\label{F(gamma)}
\xymatrix{
 F(A)\xrightarrow{(-)^{\otimes n}}(F(A)^{\otimes n})^{\Sigma_n}\longrightarrow F(A^{\otimes n})^{\Sigma_n}\stackrel{\cong}{\longleftarrow}F((A^{\otimes n})^{\Sigma_n})&F(\Gamma_n(A))\ar[l]^-{\cong}_-{F(\xi_A)}.
}
\end{equation}

\begin{theorem} \label{polylawfunct}
Any PD-functor $F\colon \Com\to \Com$ extends canonically to an endofunctor of $\Poly$ which preserves the degree. This is the unique extension which sends the universal $n$-homogeneous polynomial law $\gamma_n\colon A\to \Gamma_n(A)$ to the map (\ref{F(gamma)}) for any torsion-free $A$.
\end{theorem}

The proof of this theorem will occupy the rest of the section. The key ingredient of the proof is the construction of a natural unital ring homomorphism
\[c_n \colon \Gamma_n F(A) \longrightarrow F(\Gamma_n A),\]
for every $n\geq 0$.
For $n=0$, the functor $\Gamma_0$ is constant with value $\mathbb{Z}$, and we define $c_0$ to be the unit map 
\[c_0\colon \mathbb{Z} \longrightarrow F(\mathbb{Z})\]
of the ring $F(\Z)$.
Now suppose that $n\geq 1$. We start by defining $c_n$ for torsion-free rings, as the composite
\[\xymatrix{c_n\colon \Gamma_n F(A) \ar[r]^-{\xi_{F(A)}} & (F(A)^{\otimes n})^{\Sigma_n} \ar[r]^-{c^{\Sigma_n}} & F(A^{\otimes n})^{\Sigma_n}   & F((A^{\otimes n})^{\Sigma_n})\ar[l]_-{\cong}   & F(\Gamma_n A).\ar[l]^-{\cong}_-{F(\xi_A)} }\]
Here  $c \colon F(A)^{\otimes n} \to F(A^{\otimes n})$ is the canonical map that commutes $F$ and the coproduct, and it is $\Sigma_n$-equivariant. 
The third map is an isomorphism since by assumption PD-functors commute with invariants of finite groups. 

In order to construct  $c_n$ for rings which are not necessarily torsion-free, we use resolutions by polynomial rings. We recall the formal setup of the monadicity theorem, see e.g. \cite[Section VI.7]{MacL}. Let $\C$ and $\D$ be categories, such that $\C$ has reflexive coequalisers, and
\[U \colon \C \longrightarrow \D\]
a functor which has a left adjoint $L$, reflects isomorphisms, and preserves reflexive coequalisers. The Barr-Beck theorem (crude version) then states that $U$ is \emph{monadic}, i.e. that the category $\C$ is equivalent to the category $\D^\mathcal{T}$ of algebras over the monad $\mathcal{T}=U \circ L$, and the functor $U$ corresponds to the forgetful functor $\D^\mathcal{T} \to \D$. 

Let us now identify the full subcategory $\FreeT$ of free $\mathcal{T}$-algebras in $\D^\mathcal{T}$ with a full subcategory of $\C$. Then the inclusion $\iota \colon \FreeT \subset \C$ enjoys the following universal property: For any category $\A$ which admits reflexive coequalisers and any functor 
\[G \colon \FreeT \to \A\] 
which sends reflexive coequalisers in $\C$ to reflexive coequalisers in $\A$, there exists a unique (up to natural isomorphism) extension $\widehat{G} \colon \C \to \A$ which preserves reflexive coequalisers and a natural isomorphism 
\[\alpha \colon G \overset{\cong} \to \widehat{G} \circ \iota.\]
The value of $\widehat{G}$ on an object $C$ is given by the reflexive coequaliser $\xymatrix@C=10pt{G(TTC) \ar@<0.5ex>[r]  \ar@<-0.5ex>[r] & G(TC) \ar[r] & \widehat{G}(C)}$.
We see in particular that $\widehat{G}$ is the left Kan extension of $G$ along $\iota$, that is given any functor $H \colon \C \to \A$ and natural transformation $\beta \colon G \to H \circ \iota$, there exists a unique natural transformation $\gamma \colon  \widehat{G} \to H$ such that $\gamma \iota \circ \alpha=\beta$. In particular, if a functor $F \colon \C \to \A$ preserves reflexive coequalisers, then it is left Kan extended from its restriction to $\FreeT$. More generally, $F$ is in fact left Kan extended from any full subcategory of $\C$ containing $\FreeT$. 


We now apply this discussion to the forgetful functor $U \colon \Com_1 \to \Sets$, whose left adjoint $\Z[-] \colon \Sets \to \Com$ is the polynomial algebra functor,
to define the unital ring homomorphism $c_n \colon \Gamma_n F(A) \to F\Gamma_n (A)$ for a possibly non-torsion free ring. 

\begin{prop} \label{well-defined} For any PD-functor $F$, the natural transformation $c_n\colon \Gamma_n F \to F \Gamma_n$ on the category of torsion-free rings defined above extends uniquely to the category $\Com_1$ of all rings and unital ring homomorphisms. The value at a ring $A$ is the map induced on coequalisers
\[
\xymatrix@R=20pt{
\Gamma_n F(\Z[\Z[A]])\ar[d]^-{c_n} \ar@<0.5ex>[r]  \ar@<-0.5ex>[r] & \Gamma_n F(\Z[A]) \ar[d]^{c_n}\ar[r] & \Gamma_n F(A)\ar@{-->}[d]
\\
F\Gamma_n(\Z[\Z[A]]) \ar@<0.5ex>[r]  \ar@<-0.5ex>[r] & F\Gamma_n(\Z[A]) \ar[r] & F\Gamma_n(A).
}
\]
\end{prop}

\begin{proof} The functor $\Gamma_n(-)$ preserves reflexive coequalisers by \cite[Section IV.10]{Roby1} and the functor $F$ does so by Axiom ii), hence so do the composites $\Gamma_n F$ and $F\Gamma_n$. Thus both functors are left Kan extended from the subcategory of torsion-free rings and unital ring homomorphisms, and the natural transformation $c_n$ uniquely extends. 
\end{proof}

We will need to use that this natural transformation $c_n$ is in fact natural also with respect to the non-unital ring homomorphisms. 

\begin{prop} \label{naturaloncom} For any $n \geq 0$, the map $c_n \colon \Gamma_n F\to F\Gamma_n$ is natural with respect to morphisms in $\Com$. In other
words, for any not necessarily unital ring homomorphism $f \colon A \to B$, the diagram
\[\xymatrix@C=60pt@R=17pt{ \Gamma_n F(A) \ar[r]^{c_n} \ar[d]_-{\Gamma_n F(f)} & F\Gamma_n (A) \ar[d]^{F\Gamma_n (f)} \\  \Gamma_n F(B) \ar[r]_{c_n}  & F\Gamma_n (B)} \]
commutes.
\end{prop}

\begin{proof} 
We obtain this extra naturality again by applying the monadicity theorem. Since we are considering non-unital ring homomorphisms the forgetful functor $\Com \to \Sets$ does not admit a left adjoint, and therefore we need to find a substitute forgetful functor. Let $\Mon$ be the category of commutative monoids and not necessarily unital monoid homomorphisms. Then the functor $U \colon \Com \to \Mon$ which sends a commutative ring to its underlying multiplicative monoid has a left adjoint $\Z(-)$ which is the usual monoid ring functor. The functor $U$ preserves reflexive coequalisers (they are all computed in sets), and by the monadicity theorem $\Com$ is equivalent to the category of $U\Z(-)$-algebras.
The free objects in this case are monoid rings, which are in particular torsion-free.


Now $\Gamma_n$ and $F$ also preserve reflexive coequalisers in $\Com$ (the first preserves them in abelian groups by \cite[Section IV.10]{Roby1}, and the second by Lemma \ref{lemma:extension}), and therefore $F\Gamma_n$ and $\Gamma_nF$ are left Kan extended from the full subcategory of $\Com$ on torsion-free rings. Thus $c_n$ uniquely extends to a natural transformation on $\Com$.
\end{proof}

\begin{rem} The extension of $c_n$ in the proof of Proposition \ref{naturaloncom} coincides on $\Com_1$ with the extension obtained in Proposition \ref{well-defined}, by the uniqueness of the latter. In fact Proposition \ref{well-defined} could be considered redundant, but we prefer to define $c_n$ using the more common resolutions by polynomial rings.
\end{rem}

We are now ready to define a multiplicative polynomial law $F(f)$ for any multiplicative polynomial law $f \colon A \to B$ of finite degree. We recall from \S\ref{recallpolylaw} that the Taylor decomposition
\[f=f_0+f_1+\cdots+ f_n, \]
where $f_i$ is a homogeneous polynomial law of degree $i$, provides a factorisation in $\Poly$
\[\xymatrix@C=70pt@R=17pt{A \ar[r]^{f}  \ar[d]_-{\prod_i \gamma_i} & B, \\ \prod_{i=0}^{n} \Gamma_i A \ar[ur]_{\varphi=\oplus_i \varphi_i}    & } \]
where $\varphi_i \circ \gamma_i = f_i$.
We define the multiplicative polynomial law $F(f)$ to be the composite
\[\xymatrix{F(A) \ar[r]^-{\prod_i \gamma_i}  & \prod_{i=0}^{n} \Gamma_i F(A) \ar[r]^{\prod_i c_i} &  \prod_{i=0}^{n} F(\Gamma_i A) \ar[r]^-{\cong} & F(\prod_{i=0}^{n} \Gamma_i A) \ar[r]^-{F(\varphi)} & F(B), }  \]
where the third morphism is the inverse of the canonical morphism, which is an isomorphism since PD-functors commute with finite products. This definition recovers the original $F(f)$ when $f$ is a morphism in $\Com$, and it clearly preserves the degree. In order to prove Theorem \ref{polylawfunct} we must therefore show that this construction respects the composition of polynomial laws.

\begin{theorem} \label{functorialityinroby} Let $F$ be a PD-functor and $f \colon A \to B$ and $g \colon B \to C$ two multiplicative polynomial laws of finite degree. Then $F(g \circ f)= F(g) \circ F(f)$.
\end{theorem}

In order to prove Theorem \ref{functorialityinroby} it is convenient to recall that by Proposition \ref{alternativepolycat} the category $\Poly$ of rings and polynomial laws is equivalent to the Kleisli category $\Com_\Gamma$ of the $\mu$-filtered comonad $\prod_{i=0}^\bullet\Gamma_i$. We will then make use of the following general construction regarding Kleisli categories.  

\begin{lemma}\label{Fkleisli}
Let $T_\bullet$ be a $\mu$-filtered comonad  on a category $\C$, let $F\colon \C\to \C$ be a functor and $\eta_n\colon T_nF\to FT_n$ be natural transformation such that the diagrams
\begin{enumerate}
\item (Comultiplicativity)\hspace{.2cm}
$\xymatrix@C=50pt{
T_{nm}F\ar[d]_{\Delta_{n,m}}\ar[rr]^-{\eta_{nm}}&&FT_{nm}\ar[d]^{F(\Delta_{n,m})}
\\
T_nT_mF\ar[r]_{T_n(\eta_m)}&T_nFT_m\ar[r]_-{\eta_nT_m}&FT_nT_m
}$
\item (Unitality) \hspace{3.5cm}
$
\xymatrix@C=50pt{
T_1F\ar[r]^-{\varepsilon}\ar[d]_{\eta_1}&F
\\
FT_1\ar[ur]_-{F(\varepsilon)}
}
$
\item (Filtration) \hspace{2.8cm}
$
\xymatrix@C=50pt{
T_{n+1}F\ar[d]_{\pi_{n+1}}\ar[r]^-{\eta_{n+1}}&FT_{n+1}\ar[d]^{F(\pi_{n+1})}
\\
T_nF\ar[r]_-{\eta_n}&FT_n
}
$
\end{enumerate}
commute for every $n,m\geq 0$. Then $F$ can be extended to a functor $\C_{T_\bullet}\to \C_{T_\bullet}$  by sending an object $X$ to $F(X)$, and a morphism in $\C_{T_\bullet}$ represented by $\varphi\colon T_nX\to Y$ in $\C$ to
\[
T_nF(X)\stackrel{\eta_n}{\longrightarrow}FT_n(X)\xrightarrow{F(\varphi)} F(Y).
\] 
\end{lemma}

\begin{proof}
The proof is immediate from the definitions. The third axiom guarantees that this map does not depend on the colimit representative, and the first and the second that it preserves the compositions and the identities, respectively.  
\end{proof}

\begin{proof}[Proof of \ref{functorialityinroby}]
Our definition of $F$ on $\Poly$ arises precisely as in Lemma \ref{Fkleisli} from the maps
\[
\eta_n\colon \prod_{i=0}^n\Gamma_iF\xrightarrow{\prod_ic_i}\prod_{i=0}^nF\Gamma_i\cong F(\prod_{i=0}^n\Gamma_i)
\]
which are natural on $\Com$ by Proposition \ref{naturaloncom}. It is immediate to verify that these natural transformations $\eta_n$ are compatible with the projections $\pi_{n+1}$ and with the counit $\varepsilon$. It therefore suffices to show that $\eta_n$ is compatible with the comultiplication $\Delta_{n,m}$ as in condition $i)$ of Lemma \ref{Fkleisli}.

Let us recall from \cite[Theorem III.4]{Roby1} that for every $n\geq 0$ and rings $A$ and $B$, there is a natural group isomorphism
\[
\Gamma_n(A\times B)\stackrel{\cong}{\longrightarrow}\prod_{i+j=n}\Gamma_i(A)\otimes\Gamma_j(B)
\]
whose $(i,j)$-component classifies the $n$-homogeneous polynomial law $\gamma_i\hat{\otimes}\gamma_j\colon A\times B\to \Gamma_i(A)\otimes\Gamma_j(B)$. Here $\hat{\otimes}$ is the tensor product of multiplicative polynomial laws $\alpha\colon A\to C$ and $\beta\colon B\to D$, defined as the composite
\[
A\times B\xrightarrow{\alpha\times \beta} C\times D\stackrel{\pi}{\longrightarrow} C\otimes D,
\]
where the map $\pi$ is the canonical map from the product to the tensor product, which is a polynomial law by \cite[Section I.7]{Roby1}. Since $\gamma_i\hat{\otimes}\gamma_j$ is multiplicative the isomorphism above is a natural ring isomorphism.
Under this isomorphism the comultiplication of the $\mu$-filtered comonad is the ring homomorphism
\[
\Delta_{n,m}\colon \prod_{k=0}^{nm}\Gamma_k\longrightarrow \prod_{i=0}^{n}\Gamma_i(\prod_{j=0}^m\Gamma_j)\cong \prod_{0\leq j_0+\dots+j_m\leq n}(\Gamma_{j_0}\Gamma_0)\otimes\dots\otimes (\Gamma_{j_m}\Gamma_m)
\]
whose $(j_0,\dots,j_m)$-component is the composite of ring homomorphisms
\[
\prod_{k=0}^{nm}\Gamma_kA\longrightarrow\Gamma_{0j_0+\dots+mj_m}A\xrightarrow{\gamma_{0j_0,\dots,mj_m}}(\Gamma_{0j_0}A)\otimes\dots\otimes (\Gamma_{mj_m} A)\xrightarrow{\otimes_l\chi_{l,j_l}}(\Gamma_{j_0}\Gamma_0A)\otimes\dots\otimes (\Gamma_{j_m}\Gamma_mA)
\]
where the first is the projection, the second classifies the tensor product $\gamma_{0j_0}\hat{\otimes}\dots\hat{\otimes}\gamma_{mj_m}$, and the third is the tensor product of the ring homomorphisms $\chi_{l,k}\colon \Gamma_{lk}\to\Gamma_k\Gamma_l$ which classify $\gamma_k\circ\gamma_l$.
By using that $F$ commutes with finite products, the comultiplicativity of the $\eta_n$
 is thus equivalent to the commutativity of the outer rectangle in the diagram
\begin{equation}\label{Kleislidiag}\xymatrix{
\Gamma_{0j_0+1j_1+\dots+mj_m}F\ar[d]_{\gamma_{0j_0,\dots,mj_m}}\ar[rrr]^-{c_{0j_0+1j_1+\dots+mj_m}}&&&F\Gamma_{0j_0+1j_1+\dots+mj_m}\ar[d]^{F(\gamma_{0j_0,\dots,mj_m})}
\\
\bigotimes_{l=0}^m\Gamma_{lj_l}F\ar[d]_{\otimes_l\chi_{l,j_l}}\ar[rr]^-{\otimes_lc_{lj_l}}
&&
 \bigotimes_{l=0}^mF\Gamma_{lj_l}\ar[r]\ar[d]_-{\otimes_lF(\chi_{l,j_l})}
&
F(\bigotimes_{l=0}^m\Gamma_{lj_l})\ar[d]^-{F(\otimes_l\chi_{l,j_l})}
\\
\bigotimes_{l=0}^m\Gamma_{j_l}\Gamma_lF\ar[r]_-{\otimes_l\Gamma_{j_l}(c_{l})}
&
\bigotimes_{l=0}^m\Gamma_{j_l}F\Gamma_l
\ar[r]_-{\otimes_lc_{j_l}}
&
\bigotimes_{l=0}^mF\Gamma_{j_l}\Gamma_l\ar[r]
&
F(\bigotimes_{l=0}^m\Gamma_{j_l}\Gamma_l)\rlap{\ ,}
}
\end{equation}
for every $0\leq j_0+\dots+j_m\leq n$. The horizontal maps in the bottom right square are the canonical maps that commute $F$ and the tensor products, and that square commutes by naturality. We therefore need to show that the remaining two rectangles commute, and since all the functors involved commute with reflexive coequalisers it is sufficient to prove it on torsion-free rings.

We start with the bottom left rectangle of diagram (\ref{Kleislidiag}), where we need to show that for every pair of integers $j,k\geq 0$ and torsion-free $A$, the top rectangle in the diagram
\[
\xymatrix@C=50pt{
\Gamma_{kj}F(A)\ar[d]_{\chi_{k,j}}\ar[rr]^-{c_{kj}}&&F\Gamma_{kj}(A)\ar[d]^{F(\chi_{k,j})}
\\
\Gamma_{j}\Gamma_kF(A)\ar[r]_-{\Gamma_j(c_k)}\ar[d]_{\xi}&\Gamma_{j}F\Gamma_k(A)\ar[r]_-{c_j}\ar[d]^\xi&F\Gamma_{j}\Gamma_k(A)\ar@{>->}[d]^{F(\xi)}
\\
(\Gamma_kF(A))^{\otimes j}\ar[r]_-{c^{\otimes j}_k}&(F\Gamma_k(A))^{\otimes j}\ar[r]_-{c}&F(\Gamma_k(A)^{\otimes j})
}
\]
commutes. The vertical maps from the second to the third row are induced by the natural ring homomorphisms $ \Gamma_j   \to ((-)^{\otimes j})^{\Sigma_{j}}$ of \cite[Proposition III.1]{Roby1}, followed by the inclusion of $\Sigma_j$ invariants in  the $j$-fold tensor product  (see also \cite{Lak}). The right-hand one is injective, since it is the composite
\[
F\Gamma_{j}\Gamma_k(A)\xrightarrow{\cong}F((\Gamma_k(A)^{\otimes j})^{\Sigma_j})\cong F(\Gamma_k(A)^{\otimes j})^{\Sigma_j}\subset F(\Gamma_k(A)^{\otimes j})
\]
where the first map is an isomorphism since $A$, and therefore $\Gamma_k(A)$, is torsion-free (see \cite[Section III.6]{Roby1}). The bottom left-hand square commutes by naturality of $\xi$, and the bottom right-hand square commutes by the definition of $c_j$ on torsion free rings, where $c$ is the canonical map that commutes $F$ and the tensor product. It is therefore sufficient to show that the outer rectangle commutes.
The composite of the two left vertical maps is the natural ring homomorphism
\[
\gamma_{k,\dots,k}\colon\Gamma_{kj}\longrightarrow\Gamma_k^{\otimes j}
\]
which classifies the $kj$-homogeneous polynomial law $\gamma_k\hat{\otimes}\dots\hat{\otimes}\gamma_k\colon \id\to \Gamma_k^{\otimes j}$. This is because by definition of $\chi_{k,j}$ the composite
\[
\id\xrightarrow{\gamma_{kj}}\Gamma_{kj}\xrightarrow{\chi_{k,j}}\Gamma_{j}\Gamma_k\stackrel{\xi}{\longrightarrow} \Gamma_k^{\otimes j}
\]
is $\xi\circ\gamma_j\circ\gamma_k=(\gamma_k)^{\hat{\otimes} j}$. It follows that the composite of the two right vertical maps is $F(\gamma_{k,\dots,k})$. By induction on $j$, it is sufficient to show that for every $k_1,k_2\geq 0$ the diagram
\[
\xymatrix@C=60pt{
\Gamma_{k_1+k_2} F\ar[rr]^-{c_{k_1+k_2}} \ar[d]_{\gamma_{k_1,k_2}}
&&
F\Gamma_{k_1+k_2}\ar[d]^{F(\gamma_{k_1,k_2})}
\\
(\Gamma_{k_1}F)\otimes(\Gamma_{k_2}F)\ar[r]_-{c_{k_1}\otimes c_{k_2}}
&
(F\Gamma_{k_1})\otimes(F\Gamma_{k_2})\ar[r]_-c
&
F(\Gamma_{k_1}\otimes\Gamma_{k_2})
}
\]
commutes. We moreover observe that by induction on $m$ also the commutativity of the top rectangle of diagram (\ref{Kleislidiag}) can be reduced to the commutativity of the latter, and this will conclude the proof.

In order to show that this last diagram commutes, we need understand the map $\gamma_{k_1,k_2}$. By definition of $\gamma_{k_1,k_2}$ the diagram of ring homomorphisms
\[
\xymatrix@C=60pt@R=17pt{
 \Gamma_{k_1+k_2} A \ar[d]_-{\gamma_{k_1,k_2}}\ar[rr]^-{\Gamma_{k_1+k_2}(\Delta) }&&\Gamma_{k_1+k_2}(A\times A)\ar[d]_-{\lambda^{-1}}^-{\cong}
\\
(\Gamma_{k_1}A )\otimes (\Gamma_{k_2}A) && \displaystyle\prod_{i+j=k_1+k_2}(\Gamma_i A)\otimes (\Gamma_j A)\ar[ll]^-{\proj_{k_1, k_2}} 
}
\]
commutes for every ring $A$, where $\Delta\colon A\to A\times A$ is the diagonal map, and the lower map projects onto the summand $(i,j)=(k_1,k_2)$. The map $\lambda$ is explicitly described in \cite[Theorem III.4]{Roby1}, and it sends an elementary tensor $a\otimes b$ in the $(i,j)$-summand to
\[((\id,0)_\ast(a))\star ((0,\id)_\ast(b)),\]
where $(\id,0)_\ast\colon \Gamma_i A \to \Gamma_i (A\times A)$ is induced by the inclusion $(\id,0)\colon A\to A\times A$ in the first summand, and similarly $(0,\id)_\ast$ is induced by the inclusion in the second summand. The map $\star$ is the graded multiplication of the divided power algebra $\Gamma(A \times A)=\oplus_{n}\Gamma_n (A \times A)$ as defined in \cite[Sections III.3-5]{Roby1}.

Now we can replace $\gamma_{k_1,k_2}$ by the composite $\proj_{k_1, k_2} \circ \lambda^{-1} \circ \Gamma_{k_1+k_2}(\Delta)$ into the diagram above. Using the definition of the maps $c_{k_1+k_2}$, $c_{k_1}$ and $c_{k_2}$, and Property \ref{lemma:ext-prop1} of Lemma \ref{lemma:ext-prop}, one can reduce the commutativity of this diagram  to the commutativity of 
\[
\xymatrix@C=30pt@R=15pt{
\Gamma_{k_1+k_2} (F(A)\times F(A)) \ar[dd]_{\lambda^{-1}}^{\cong}
&
\Gamma_{k_1+k_2} F(A\times A)\ar[l]_-{\cong}\ar[r]^-{c_{k_1+k_2}}
&
F\Gamma_{k_1+k_2}(A\times A)\ar[d]^{F(\lambda^{-1})}_{\cong}
\\
&&\displaystyle F(\prod_{i+j=k_1+k_2}(\Gamma_{i}A)\otimes(\Gamma_{j}A))\ar[d]^-{\cong}
\\
\displaystyle\prod_{i+j=k_1+k_2}(\Gamma_{i}FA)\otimes(\Gamma_{j}FA)\ar[r]_-{\prod c_{i}\otimes c_{j}}
&
\displaystyle\prod_{i+j=k_1+k_2}(F\Gamma_{i}A)\otimes(F\Gamma_{j}A)\ar[r]_-c
&
\displaystyle\prod_{i+j=k_1+k_2}F((\Gamma_{i}A)\otimes(\Gamma_{j}A))\rlap{\ .}
}
\]
Since $\lambda^{-1}$ is an isomorphism we can equivalently verify that the rectangle obtained by replacing  $\lambda^{-1}$ with $\lambda$ commutes, and this can be verified one component at the time. By inputting the definition of $\lambda$ we need to verify that the outer  rectangle of
\[\hspace{-1.3cm}
\xymatrix@C=15pt@R=18pt{\Gamma_{k_1} F(A) \otimes \Gamma_{k_2} F(A) \ar[d]_-{\xi_{F(A)}\otimes \xi_{F(A)}} \ar[rr]_-{(\id,0)_\ast\otimes (0,\id)_\ast} 
&&
 \Gamma_{k_1} (F(A)^{\times 2}) \otimes \Gamma_{k_2} (F(A)^{\times 2}) \ar[r]^-{\star} \ar[d]_-{\xi_{F(A)^{\times 2}}\otimes \xi_{F(A)^{\times 2}}} 
 & \Gamma_{k_1+k_2}(F(A)^{\times 2}) \ar[d]^{\xi_{F(A)^{\times 2}}}
\\
(F(A)^{\otimes k_1})^{\Sigma_{k_1}}\otimes (F(A)^{\otimes k_2})^{\Sigma_{k_2}}\ar[dd]_-{c^{\Sigma_{k_1}}\otimes c^{\Sigma_{k_2}}}\ar[rr]_-{(\id,0)_\ast\otimes (0,\id)_\ast}
&&((F(A)^{\times 2})^{\otimes k_1})^{\Sigma_{k_1}}\otimes ((F(A)^{\times 2})^{\otimes k_2})^{\Sigma_{k_2}}\ar[r]^-{\star}
& ((F(A)^{\times 2})^{\otimes k_1+k_2})^{\Sigma_{k_1+k_2}}
\\
&&(F(A^{\times 2})^{\otimes k_1})^{\Sigma_{k_1}}\otimes (F(A^{\times 2})^{\otimes k_2})^{\Sigma_{k_2}}\ar[r]^-{\star}\ar[u]_-{\cong}\ar[d]^-{c^{\Sigma_{k_1}}\otimes c^{\Sigma_{k_2}}}
&(F(A^{\times 2})^{\otimes k_1+k_2})^{\Sigma_{k_1+k_2}}\ar[d]^{c^{\Sigma_{k_1+k_2}}}\ar[u]_-{\cong}
\\
F(A^{\otimes k_1})^{\Sigma_{k_1}}\otimes F(A^{\otimes k_2})^{\Sigma_{k_2}}\ar[rr]_-{F(\id,0)_\ast\otimes F(0,\id)_\ast}
&&
F((A^{\times 2})^{\otimes k_1})^{\Sigma_{k_1}}\otimes F((A^{\times 2})^{\otimes k_2})^{\Sigma_{k_2}}\ar@{-->}[r]^-s
&F((A^{\times 2})^{\otimes k_1+k_2})^{\Sigma_{k_1+k_2}}
\\
F((A^{\otimes k_1})^{\Sigma_{k_1}})\otimes F((A^{\otimes k_2})^{\Sigma_{k_2}})\ar[u]^-{\cong}\ar[d]_-{c}&&
\\
F((A^{\otimes k_1})^{\Sigma_{k_1}}\otimes (A^{\otimes k_2})^{\Sigma_{k_2}})\ar[rrr]_-{F((\id,0)_\ast\star (0,\id)_\ast)}&&&F(((A^{\times 2})^{\otimes k_1+k_2})^{\Sigma_{k_1+k_2}})\ar[uu]^-{\cong}
\\
F(\Gamma_{k_1}A\otimes \Gamma_{k_2}A)\ar[u]^-{F(\xi_{A}\otimes\xi_{A})}_-{\cong}\ar[rrr]_-{F((\id,0)_\ast\star (0,\id)_\ast)}&&&F(\Gamma_{ k_1+k_2}(A\times A))\ar[u]_-{F(\xi_{A^{\times 2}})}^-{\cong}
}
\]
commutes  for every torsion-free ring $A$. Here the map $\star$ denotes the graded multiplication of the divided power algebra, as well as the shuffle product of the symmetric tensor algebra, with respect to which $\xi\colon \bigoplus_{k\geq 0}\Gamma_k\to \bigoplus_{k\geq 0}((-)^{\otimes k})^{\Sigma_k}$ is a natural ring homomorphism (see \cite{Lak}). We recall that on $x \in (X^{\otimes k_1})^{\Sigma_{k_1}}$ and $y \in (X^{\otimes k_2})^{\Sigma_{k_1}}$ the shuffle product is  defined as
\[x \star y=\sum_{ \sigma \in S_{k_1,k_2}} \sigma (x \otimes y),\]
where $S_{k_1,k_2}$ is the set of $(k_1, k_2)$-shuffles (which are representatives of the left cosets $(\Sigma_{k_1} \times \Sigma_{k_2})\backslash\Sigma_{k_1+k_2}$),
where the left action of $\sigma \in \Sigma_k$ on $X^{\otimes k}$ is defined by the formula
\[\sigma (x_1 \otimes x_2 \otimes \dots \otimes x_k)=x_{\sigma^{-1}(1)} \otimes  x_{\sigma^{-1}(2)} \otimes \dots \otimes  x_{\sigma^{-1}(k)}\]
(there is some confusion between left and right actions in \cite[Section III.5]{Roby1}, see also  \cite{Lak}).
All the rectangles in the last diagram except the ones involving the dashed map commute by naturality, or because $\xi$ is a map of graded algebras.
We define the dashed map by the formula
\[
s(z\otimes w):=\sum_{\sigma\in S_{k_1,k_2}}\sigma c(z\otimes w),
\]
where $z \in F((A^{\times 2})^{\otimes k_1})^{\Sigma_{k_1}}$ and $w \in F((A^{\times 2})^{\otimes k_2})^{\Sigma_{k_2}}$ and $c(z\otimes w)$ is the image of $z\otimes w$ under the composite
\[\xymatrix{F((A^{\times 2})^{\otimes k_1})^{\Sigma_{k_1}}\otimes F((A^{\times 2})^{\otimes k_2})^{\Sigma_{k_2}} \ar[r] & F((A^{\times 2})^{\otimes k_1})\otimes F((A^{\times 2})^{\otimes k_2}) \ar[r]^-c & F((A^{\times 2})^{\otimes k_1+k_2}).}\]
It is immediate to verify that $s$ lands in the $\Sigma_{k_1+k_2}$-invariants.
The rectangle above $s$ commutes since
\begin{align*}
c^{\Sigma_{k_1+k_2}}(x\star y)&=c^{\Sigma_{k_1+k_2}}(\sum_{\sigma\in S_{k_1,k_2}}\sigma(x\otimes y))=\sum_{\sigma\in S_{k_1,k_2}}\sigma c(x\otimes y)
\\
&=\sum_{\sigma\in S_{k_1,k_2}}\sigma c(c^{\Sigma_{k_1}}(x)\otimes c^{\Sigma_{k_2}}(y))
\\
&=s(c^{\Sigma_{k_1}}(x)\otimes c^{\Sigma_{k_2}}(y)).
\end{align*}

Let us denote by $\iota_{n}\colon (X^{\otimes n})^{\Sigma_{n}}\to X^{\otimes n}$ the fixed points inclusion. 
It is sufficient to show that the rectangle below $s$ commutes after composing with $\iota_{k_1+k_2}$.
After postcomposing with this inclusion,  the lower composite of this rectangle sends an element $u\otimes v \in F((A^{\otimes k_1})^{\Sigma_{k_1}})\otimes F((A^{\otimes k_2})^{\Sigma_{k_2}})$ to
\begin{align*}
F(\iota_{k_1+k_2})F((\id,0)_\ast\star (0,\id)_\ast)(c(u\otimes v))&=F(\iota_{k_1+k_2} ((\id,0)_\ast\star (0,\id)_\ast))(c(u\otimes v))
\\&=F(((\id,0)_*\star (0,\id)_*) (\iota_{k_1} \otimes \iota_{k_2}) )(c(u\otimes v))
\\
&=F(\sum_{\sigma\in S_{k_1,k_2}}\sigma ((\id,0)_*\otimes (0,\id)_*) (\iota_{k_1} \otimes \iota_{k_2}))(c(u\otimes v)).
\end{align*}
The last sum is a sum of orthogonal ring homomorphisms since the permutations of $S_{k_1,k_2}$ are shuffles, and by Property \ref{lemma:ext-prop3} of Lemma \ref{lemma:ext-prop} we can write this as
\begin{align*}
&F(\sum_{\sigma\in S_{k_1,k_2}}\sigma ((\id,0)_*\otimes (0,\id)_*) (\iota_{k_1} \otimes \iota_{k_2}))(c(u\otimes v))\\
&=\sum_{\sigma\in S_{k_1,k_2}}F(\sigma ((\id,0)_*\otimes (0,\id)_*) (\iota_{k_1} \otimes \iota_{k_2}))(c(u\otimes v))
\\&=\sum_{\sigma\in S_{k_1,k_2}}\sigma F(((\id,0)_*\otimes (0,\id)_*) (\iota_{k_1} \otimes \iota_{k_2}))(c(u\otimes v))
\\
&=\sum_{\sigma\in S_{k_1,k_2}}\sigma c(F((\id,0)_*)\otimes F((0,\id)_*)) (F(\iota_{k_1})(u) \otimes F(\iota_{k_2})(v))\end{align*}
which is the value of the upper composite.
 \end{proof}

\subsection{The ghost components of a polynomial law}\label{ghostlaw}

A consequence of Theorem \ref{polylawfunct} and Example \ref{Witt vectors} is that any multiplicative polynomial law $f\colon A\to B$ of finite degree induces a multiplicative polynomial law on Witt vectors
\[
W_S(f)\colon W_S(A)\longrightarrow W_S(B),
\]
for every truncation set $S\subset\mathbb{N}$. We will describe the ghost components of $W_S(f)$, and explain how these determine the functoriality of $W_S$ in multiplicative polynomial laws.

\begin{prop} \label{wittnaturalpoly} Let $F,G\colon \Com\to \Com$ be PD-functors and $\alpha \colon F \to G$ a natural transformation. Then $\alpha$ extends to a natural transformation on $\Poly$. That is, for any multiplicative polynomial law $f \colon A \to B$ of finite degree, the diagram of polynomial laws
\[\xymatrix@C=70pt@R=15pt{ F(A) \ar[r]^{F(f)} \ar[d]_{\alpha_A} & F(B) \ar[d]^{\alpha_B} \\ G(A) \ar[r]^{G(f)} & G(B) }\]
commutes. 
\end{prop}

\begin{proof} Suppose that degree of $f$ is at most $n$. Recall that we have a commutative diagram of polynomial laws
\[\xymatrix@C=70pt@R=15pt{A \ar[r]^{f}  \ar[d]_-{\prod_i \gamma_i} & B. \\ \prod_{i=0}^{n} \Gamma_i A \ar[ur]_{\varphi=\oplus_i \varphi_i}    & } \]
We need to show that the outer rectangle in the diagram 
\[\xymatrix@C=40pt{F(A) \ar[d]_{\alpha_A} \ar[r]^-{\prod_i \gamma_i}  & \prod_{i=0}^{n} \Gamma_i F(A) \ar[d]^{\prod_i \Gamma_i (\alpha_A)} \ar[r]^{\prod_i c_i} &  \prod_{i=0}^{n} F(\Gamma_i A) \ar[d]^{\prod_i \alpha_{ \Gamma_i(A)}} \ar[r]^{\cong} & F(\prod_{i=0}^{n} \Gamma_i A) \ar[r]^-{F(\varphi)} \ar[d]^{\alpha_{ \prod_i \Gamma_i(A)}} & F(B) \ar[d]^{\alpha_B} \\ G(A) \ar[r]^-{\prod_i \gamma_i}  & \prod_{i=0}^{n} \Gamma_i G(A) \ar[r]^{\prod_i c_i} &  \prod_{i=0}^{n} G(\Gamma_i A) \ar[r]^{\cong} & G(\prod_{i=0}^{n} \Gamma_i A) \ar[r]^-{G(\varphi)} & G(B) }  \]
commutes. The first square commutes by naturality of the universal polynomial laws $\gamma_i$. The third and fourth squares commute by the naturality of $\alpha$ on $\Com$. It remains to check that the second square commutes. For this it suffices to see that for any $i$, the square
\[\xymatrix@C=70pt@R=15pt{ \Gamma_i F(A) \ar[d]_{\Gamma_i \alpha_A} \ar[r]^{c_i}  & F(\Gamma_i A) \ar[d]^{\alpha_{\Gamma_iA}} \\ \Gamma_i G(A) \ar[r]^{c_i} &  G(\Gamma_i A)} \]
commutes. 
Because $c_i$ is extended to all rings from torsion free rings it suffices to check this when $A$ is torsion-free. By the naturality of $\alpha$ this reduces to showing that for any commutative rings $C_1$ and $C_2$, the diagram
\[\xymatrix@C=70pt@R=15pt{F(C_1) \otimes F(C_2) \ar[d]_{\alpha \otimes \alpha} \ar[r]^c & F(C_1 \otimes C_2) \ar[d]^{\alpha} \\ G(C_1) \otimes G(C_2) \ar[r]^c &  G(C_1 \otimes C_2)}\]
commutes. This also follows from naturality of $\alpha$. 
\end{proof}

Let us recall that for any set $S$, the $S$-fold product functor $(-)^{\times S}\colon \Com\to \Com$ is a PD-functor, and therefore a multiplicative polynomial law $f\colon A\to B$ of finite degree induces a multiplicative polynomial law $f^{\times S}\colon A^{\times S}\to B^{\times S}$. When the set $S$ is finite this is the natural transformation
\[
(f^{\times S})_R\colon (A^{\times S})\otimes_{\Z} R\cong (A\otimes_{\Z}  R)^{\times S}\xrightarrow{(f_R)^{\times S}} (B\otimes_{\Z}  R)^{\times S}\cong (B^{\times S})\otimes_{\Z} R.
\]
If $S$ is infinite there is no obvious direct description of this law as a natural transformation, and one needs to involve divided powers.
%

\begin{prop} \label{levelghost} Let $f \colon A \to B$ be a multiplicative polynomial law of finite degree. Then for any truncation set $S$, the diagram of polynomial laws
\[\xymatrix@C=70pt@R=15pt{W_S(A) \ar[d]_w \ar[r]^{W_S(f)} & W_S(B) \ar[d]^w \\ \prod_S A \ar[r]^{\prod_S f} & \prod_S B  }\]
commutes, where the vertical maps are the ghost coordinates. Moreover, $W_S\colon \Poly\to\Poly$ is the unique extension of $W_S$ with this property.
\end{prop}

\begin{proof} The ghost coordinates $w\colon W_S\to (-)^{\times S}$ form a natural transformation between PD-functors. Thus by Proposition \ref{wittnaturalpoly} the diagram commutes. Let $W_{S}'$ be another extension of $W_S$ to $\Poly$ such that the diagram above commutes. Then we must have $wW_{S}(f)=w'W_{S}(f)$, or equivalently that the composites
\[
\xymatrix@C=50pt{
\prod_{i=0}^n\Gamma_iW_S(A)\ar@<.5ex>[r]^{\varphi}\ar@<-.5ex>[r]_{\varphi'}&W_S(B)\ar[r]^{w}&B^{\times S}
}
\]
agree, where $\varphi$ and $\varphi'$ are the unique ring homomorphisms extending $W_S(f)$ and $W'_S(f)$, respectively. When $B$ is torsion free $w$ is injective, and this shows that $\varphi=\varphi'$, and consequently $W_S(f)=W'_S(f)$. In general, let us consider the commutative square
\[
\xymatrix@C=70pt@R=15pt{
\Z[A]\ar@{->>}[d]_{\varepsilon}\ar[r]^{\prod_i\gamma_i}&\prod_{i=0}^n\Gamma_i\Z[A]\ar[d]^{\varphi}
\\
A\ar[r]^f&B
}
\]
where $\varepsilon$ is the counit from the polynomial ring, the top horizontal map is the universal polynomial law, and $\varphi$ is the unique ring homomorphism extending the polynomial law $f\varepsilon$. By applying $W_S$ and $W'_S$, respectively, to this diagram we obtain analogous commutative squares of polynomial laws, and corresponding commutative squares of ring homomorphisms
\[
\xymatrix@C=40pt{
\prod_{i=0}^n\Gamma_iW_S(\Z[A])\ar@{->>}[d]^{\prod_i\Gamma_iW_S(\varepsilon)}\ar[r]^{\varphi_{W_S(\prod_i\gamma_i)}}&W_S(\prod_{i=0}^n\Gamma_i\Z[A])\ar[d]^{W_S(\varphi)}
\\
\prod_{i=0}^n\Gamma_iW_S(A)\ar[r]_{\varphi_{W_S(f)}}& W_S(B).
}
\ \ \ \ \ \ \  \ \ \ 
\xymatrix@C=40pt{
\prod_{i=0}^n\Gamma_iW_S(\Z[A])\ar@{->>}[d]^{\prod_i\Gamma_iW_S(\varepsilon)}\ar[r]^{\varphi_{W'_S(\prod_i\gamma_i)}}&W_S(\prod_{i=0}^n\Gamma_i\Z[A])\ar[d]^{W_S(\varphi)}
\\
\prod_{i=0}^n\Gamma_iW_S(A)\ar[r]_{\varphi_{W'_S(f)}}&W_S(B).
}
\]
The vertical maps of these squares agree since $W_S$ and $W_S'$ agree on ring homomorphisms by assumption. The top horizontal maps also agree by the previous argument, since $\prod_{i=0}^n\Gamma_i\Z[A]$ is torsion-free. Since the left vertical map is surjective, the bottom horizontal maps are also equal.
\end{proof}

\subsection{The product of polynomial laws}

We show that our construction respects the product of polynomial laws. Given two multiplicative polynomial laws $f \colon A \to B$ of degree at most $n$ and $g \colon A \to B$ of degree at most $m$, there is a multiplicative polynomial law 
\[(f\cdot g)_R:=f_R\cdot g_R\colon  A\otimes_{\Z}R \longrightarrow B\otimes_{\Z} R\]
defined by the pointwise product in the ring $B\otimes_{\Z} R$, of degree at most $n+m$. 

\begin{prop} \label{productroby} Let $F$ be a PD-functor. For any pair of  finite degree multiplicative polynomial laws $f, g \colon A \to B$, we have
$F(f \cdot g)=F(f) \cdot F(g)$. \end{prop}

\begin{proof} This is very similar to the proof of Theorem \ref{functorialityinroby}. Suppose that $f$ is of degree at most $n$ and $g $ of degree at most $m$. Then we have commutative diagrams
\[\xymatrix@C=70pt@R=18pt{A \ar[r]^{f}  \ar[d]_-{\prod_i \gamma_i} & B, \\ \prod_{i=0}^{n} \Gamma_i A \ar[ur]_{\varphi}    & } 
\ \ \ \ \ \ \ \ \ \ \ \ \ \ \ \ \ \ 
\xymatrix@C=70pt@R=18pt{A \ar[r]^{g}  \ar[d]_-{\prod_j \gamma_j} & B, \\ \prod_{j=0}^{m} \Gamma_j A \ar[ur]_{\psi},    & }\]
where $\varphi$ and $\psi$ are not necessarily unital ring homomorphisms. In order to prove the proposition it suffices to show that the diagram
\[\xymatrix@C=28pt{& & & F(\prod_{i=0}^{n} \Gamma_i A) \otimes  F(\prod_{j=0}^{m} \Gamma_i A) \ar[d]^c \ar[rr]^-{F(\varphi) \otimes F(\psi)} & & F(B) \otimes F(B) \ar@/^{4pc}/[dd]^{\mu} \ar[d]^c \\ F(A) \ar[urrr]^{\hspace{-1.5cm}F(\prod_i \gamma_i) \hat{\otimes} F(\prod_j \gamma_j)} \ar[drrrrr]_-{F(f \cdot g)} \ar[rrr]^-{\hspace{1.2cm} F((\prod_i \gamma_i) \hat{\otimes} (\prod_j \gamma_j)) } & & & F(\prod_{i=0}^{n} \Gamma_i A \otimes \prod_{j=0}^{m} \Gamma_j A) \ar[rr]^-{F(\varphi \otimes \psi)} & & F(B \otimes B) \ar[d]^{F(\mu)} \\ & & & & & F(B)  }\]
commutes, where $\mu$ denotes the ring multiplications and $\hat{\otimes}$ is the tensor product of polynomial laws of the proof of Theorem  \ref{functorialityinroby}. Indeed, the outer composite through the upper right hand corner is $F(f) \cdot F(g)$ by functoriality of $F$. The lower triangle commutes by functoriality of $F$, and the square commutes by Lemma \ref{lemma:ext-prop} \ref{lemma:ext-prop1}. The upper left hand triangle commutes by the final step in the proof of Theorem  \ref{functorialityinroby}.
\end{proof}

\section{On the functoriality of the Witt vectors in polynomial maps}\label{secpolymap}

In this section we show that in certain circumstances the functoriality of the Witt vectors functors in polynomial laws extends to polynomial maps. We start by reviewing some material on polynomial maps.

\subsection{Review of polynomial maps}

This is a recollection of results on polynomial maps and their relationship with polynomial laws that we will use throughout the paper. The content is classical and we do not claim originality for these results. Some of these results can be found in \cite{Passi1, Passi2, Leib, Xan}.

Let $A$ and $B$ be abelian groups and $f\colon A\to B$ a function which is not necessarily a group homomorphism, and $n\geq 0$ an integer. We recall that the \emph{$n$-th cross-effect}, or \emph{$n$-th deviation}, of $f$ is the function $\cross_n\colon A^{\times n}\to B$ defined by
\[
\cross_{n} f(a_1,\dots, a_{n}):=\sum_{\substack{U\subset \{1,\cdots, n \} }} (-1)^{n-|U|}f(\sum_{l\in U}a_l).
\]

\begin{defn}
A function $f\colon A\to B$ of abelian groups is called \emph{polynomial of degree} $\leq n$, or $n$-\emph{polynomial}, if $\cross_{n+1}f=0$. It is called $n$-homogeneous if it is $n$-polynomial and $f(ka)=k^nf(a)$ for every $k\in\Z$ and $a\in A$.

A function of rings  $f\colon A\to B$ is called a \emph{multiplicative $n$-polynomial map} if it is $n$-polynomial as a map of abelian groups, and $f(aa')=f(a)f(a')$ for any $a, a' \in A$. Similarly, it is multiplicative $n$-homogeneous if it is multiplicative $n$-polynomial and $n$-homogeneous.
\end{defn}

The composition of an $n$-polynomial map and an $m$-polynomial map is $nm$-polynomial, by \cite{Leib}, and similarly for $n$-homogeneous maps. Similarly the product of an $n$-polynomial map and an $m$-polynomial map is $(n+m)$-polynomial.

\begin{example}\
\begin{enumerate}
\item The only $0$-polynomial maps are the constant maps.
\item A multiplicative polynomial map $f\colon A\to B$ of degree $1$ is precisely a map of the form $f(a)=c + \varphi(a)$, $a \in A$, where $c$ is a constant idempotent and $\varphi$ is a  not necessarily unital ring homomorphism from $A$ to $B$ which is orthogonal to $c$.
\item The exponentiation map $(-)^n\colon A\to A$ is a multiplicative $n$-homogeneous map for every commutative ring $A$. This is the case since it is the $n$-fold product of the identity map with itself, which is $1$-homogeneous.
\item The map $N\colon \Z\to \Z[x]/(x^2-2x)$ defined by $N(a)=a+\frac{a(a-1)}{2}x$ is multiplicative $2$-polynomial, but not homogeneous. This is the multiplicative norm of the Burnside Tambara functor of the group $\Z/2$. It is an instance of the following more general example.
\item Let $T$ be a Tambara functor for a finite group $G$. The multiplicative transfer 
\[N(f)\colon T(G/H)\to T(G/K)\] 
induced by a $G$-equivariant map $f\colon G/H\to G/K$ is $[K:H^g]$-polynomial, where $H^g$ is a subgroup of $K$ conjugate to $H$ with $f(eH)=gK$. This is proved in \cite{Tambara} \cite[13.22]{Strickland}. For example for the groups $0=H\leq K=G=\Z/2$, the Tambara reciprocity formula for $N=N(\Z/2\to\ast)$ gives
 \[N(a+b)=N(a)+N(b)+\tran(a\overline{b})\]
 where $\tran$ is the additive transfer and the bar denotes the involution of $T(\Z/2)$. In this case one can explicitly calculate that
\begin{align*}
cr_3 N(a,b,c)&=N(a+b+c)-N(a+b)-N(b+c)-N(a+c)+N(a)+N(b)+N(c)
\\ &=N(a)+N(b)+N(c)+\tran(a\overline{b})+\tran(a\overline{c})+\tran(b\overline{c})
\\&-N(a+b)-N(b+c)-N(a+c)+N(a)+N(b)+N(c)\\
&=0.
\end{align*}
\end{enumerate}
\end{example}

Not all polynomial maps can be decomposed into a sum of homogeneous maps. A well-known counterexample is the degree $2$ map
\[n \mapsto  {n\choose 2} \colon \mathbb{Z} \to  \mathbb{Z},\]
 (see e.g., \cite{Gaud-Hart}). However, this is possible when sufficiently many integers are invertible. 

\begin{prop}\label{homdecomp}
Let $f\colon A\to B$ be $n$-polynomial, and suppose that every integer $1 \leq k \leq n$ is invertible in $B$ (e.g., if $B$ is $p$-local for some prime $p>n$). Then $f$ decomposes uniquely as
\[
f=\sum_{k=0}^n\varphi_k
\]
where each $\varphi_k$ is $k$-polynomial and homogeneous. If moreover $f$ is multiplicative, so are the $\varphi_k$ and
\[\varphi_i(x) \varphi_j(y)=0\]
for any $x, y \in A$ and $i \neq j$. 
\end{prop}

\begin{proof}
By definition the $n$-th cross-effect of an  $n$-homogeneous function $h$ satisfies
\[(\cross_n  h)(x,\dots,x)=\sum_{i=0}^{n}(-1)^{n-i}{n\choose i}i^nh(x)=n!h(x).\]
First let us show that the decomposition into homogeneous summands is unique. For simplicity we introduce the notation
\[(\crr_n \alpha)(x):=(\cross_n \alpha) (x, \cdots, x).\]
If  $f=\sum_{k=0}^n\varphi_k$, then $f-\varphi_n$ is $(n-1)$-polynomial, and
\[
0=\crr_{n}(f-\varphi_n)=\crr_n f-\crr_{n}\varphi_n=\crr_n f-n!\varphi_n.
\]
Since by assumption $n!$ is invertible, $\varphi_n$ is uniquely determined by $f$. Similarly, $f-\varphi_n-\varphi_{n-1}-\dots-\varphi_{k}$ is $(k-1)$-polynomial, and
\begin{align*}
0&=\crr_{k}(f-\varphi_n-\varphi_{n-1}-\dots-\varphi_{k})=\crr_k (f-\varphi_n-\varphi_{n-1}-\dots-\varphi_{k+1})-\crr_{k}\varphi_k
\\
&=\crr_k (f-\varphi_n-\varphi_{n-1}-\dots-\varphi_{k+1})-k!\varphi_k
\end{align*}
shows that $\varphi_k$ is inductively determined by the $\varphi_j$ for $k<j\leq n$.

The proof of uniqueness gives us an inductive procedure to define the maps $\varphi_k$.
We recall that the $n$-th cross-effect $\cross_n f$ of an $n$-polynomial map is additive in each variable, and therefore the diagonal $\crr_n f$ is $n$-homogeneous.
We set
\[
\varphi_n:=\frac{1}{n!}\crr_n f,
\]
and we define inductively 
\[
\varphi_k:=\frac{1}{k!}\crr_k(f-\varphi_n-\dots-\varphi_{k+1})
\]
for every $0\leq k<n$. A simple inductive argument shows that each $\varphi_k$ is $k$-homogeneous.

Now let us show that if $f$ is multiplicative, then so is each $\varphi_k$ and the different summands in the decomposition are orthogonal to each other. Indeed, from the equation $f(xy)=f(x)f(y)$ we see that
\[
\sum_{k=0}^n\varphi_k(xy)=\sum_{k,j=0}^n\varphi_k(x)\varphi_j(y)=\sum_{k=0}^n\varphi_k(x)(\sum_{j=0}^n\varphi_j(y)).
\]
The function $f(xy)$ is $n$-polynomial in $x$ and $\varphi_k(xy)$ is $k$-homogeneous in $x$ for any fixed $y$. By the uniqueness of the homogeneous decomposition it follows that
\[
\varphi_k(xy)=\varphi_k(x)(\sum_{j=0}^n\varphi_j(y))
\]
for every $k$ and every fixed $y\in A$. This is now a $k$-homogeneous polynomial map in $y$ for any fixed $x$, and again by the uniqueness of the homogeneous decomposition
\[\varphi_k(xy)=\varphi_k(x)\varphi_k(y) \ \ \ \ \ \ \ \ \ \ 
\mbox{and} 
\ \ \ \ \ \ \ \ \ \ 0=\varphi_i(x)\varphi_j(y)
\]
when $i \neq j$.
\end{proof}

\begin{example}
Let $N^{G}_H\colon T(G/H)\to T(G/G)$ be the norm-map of a $G$-Tambara functor $T$, for some finite group $G$. Then $N^{G}_H$ is a multiplicative polynomial map whose degree is equal to the index $[G:H]$. This map does not in general extend to a polynomial law (See Example \ref{cex} below), nor does it decompose into a sum of homogeneous pieces. After inverting the group order there is an isomorphism of rings
\[
T(G/G)[|G|^{-1}]\stackrel{\cong}{\longrightarrow} \prod_{(K\leq G)}(T(G/K)[|G|^{-1}]/J_K)^{W_G(K)},
\]
where the product runs though the conjugacy classes of subgroups of $G$, $J_K$ is the sum of the images of the additive transfers $T(G/L)\to T(G/K)$ where $L$ is a proper subgroup of $K$, and $W_G(K)$ is the Weyl group of $K$ in $G$ \cite[Proposition 3.4.18]{Schwedeglobal}. The isomorphism is induced by the restrictions $\res^{G}_K$. The $K$-component of the composite map
\[
T(G/H)\stackrel{N}{\longrightarrow}T(G/G)\longrightarrow T(G/G)[|G|^{-1}]\cong \prod_{(K\leq G)}(T(G/K)[|G|^{-1}]/J_K)^{W_G(K)}
\]
is then homogeneous of degree $|K\backslash G/H|$. 
It suffices to check this after postcomposing with the inclusion 
\[\prod_{(K\leq G)}(T(G/K)[|G|^{-1}]/J_K)^{W_G(K)} \subset \prod_{(K\leq G)}T(G/K)[|G|^{-1}]/J_K. \]
By the double coset formula modulo $J_K$ we get
\[
\res_{K}^GN_{H}^G=\prod_{[g]\in K\backslash G/H}N^{K}_{H^g\cap K}c_g\res_{H\cap K^g}^H,
\]
where the conjugation map $c_g$ and $\res_{H\cap K^g}^H$ are ring-homomorphisms. By reciprocity
\[
N^{K}_{H^g\cap K}(mx)=mN^{K}_{H^g\cap K}(x)+t(x)
\]
for $m \in \Z$, where $t(x)$ is a sum of proper transfers which vanishes in the quotient.
\end{example}

We remark that if $f\colon A\to B$ is an $n$-homogeneous polynomial law, the underlying map of abelian groups 
\[f_{\Z} \colon A\longrightarrow B\] 
obtained by evaluating the law at the ring $R=\Z$, is an $n$-homogeneous polynomial map. This can be verified for the universal $n$-homogeneous law $\gamma_n\colon A\to \Gamma_n(A)$, where
\begin{align*}
\cross_{n+1}\gamma_n(a_1,\dots,a_n,a_{n+1})&=\cross_{n}(\gamma_n(a_{n+1}+(-)))(a_1,\dots ,a_n)-\cross_{n}\gamma_n(a_1,\dots, a_n)\\
&=\cross_{n}(\sum_{i+j=n}\gamma_i(a_{n+1})\gamma_{j}(-))(a_1,\dots,a_n)-\cross_{n}\gamma_n(a_1,\dots,a_n)\\
&=\sum_{i+j=n}\gamma_i(a_{n+1})(\cross_{n}\gamma_{j})(a_1,\dots,a_n)-\cross_{n}\gamma_n(a_1,\dots,a_n)\\
&=\gamma_0(a_{n+1})(\cross_{n}\gamma_{n})(a_1,\dots,a_n)-\cross_{n}\gamma_n(a_1,\dots,a_n)=0
\end{align*}
since the terms of the sum with $j<n$ vanish by induction on $n$, and $\gamma_0$ is the constant function with value $1$. 

\begin{prop} \label{polymapvspolylaw}
Let $f\colon A\to B$ be an $n$-polynomial map, and $p>n$ a prime. Then the composite $f_{(p)}\colon A\to B\to B_{(p)}$ extends to a unique polynomial law of degree $n$. If $f$ is multiplicative, then so is the extension. 
\end{prop}

\begin{proof} For convenience we will sometimes denote an element $\gamma_i(a) \in \Gamma_i(A)$ just by $a^{(i)}$. 

We start by showing that any $n$-homogeneous polynomial map $\varphi \colon A\to B_{(p)}$ extends to an $n$-homogeneous law, which is multiplicative if the original map was. We will now define a unique additive extension $\overline{\varphi}\colon \Gamma_n(A) \to B_{(p)}$ of $\varphi$ as follows. The generators of the form $a^{(1)}_1\dots a^{(1)}_n$ can always be expressed as a sum of generators of the form $a^{(n)}$, by the formula
\[
a^{(1)}_1\dots a^{(1)}_n=(\cross_{n}\gamma_n)(a_1,\dots, a_n)=\sum_{U\subset n}(-1)^{n-|U|}(\sum_{i\in U}a_i)^{(n)}.
\]
This formula can easily be proven by induction on $n$. Indeed it clearly holds for $n=0$ and $n=1$, and in general
\begin{align*}
\cross_{n}\gamma_n(a_1,\dots,a_{n-1},a_{n})&=\cross_{n-1}(\gamma_n(a_{n}+(-)))(a_1,\dots, a_{n-1})-\cross_{n-1}\gamma_n(a_1,\dots,a_{n-1})\\
&=\cross_{n-1}(\sum_{i+j=n}\gamma_i(a_{n})\gamma_{j}(-))(a_1,\dots,a_{n-1})-\cross_{n-1}\gamma_n(a_1,\dots,a_{n-1})\\
&=\sum_{i+j=n}\gamma_i(a_{n})(\cross_{n-1}\gamma_{j}(-))(a_1,\dots,a_{n-1})-\cross_{n-1}\gamma_n(a_1,\dots,a_{n-1})\\
&=\sum_{\substack{i+j=n\\ j<n}}\gamma_i(a_{n})(\cross_{n-1}\gamma_{j})(a_1,\dots,a_{n-1})\\
&=\gamma_1(a_{n})(\cross_{n-1}\gamma_{n-1})(a_1,\dots,a_{n-1})=\gamma_1(a_{n})a^{(1)}_1\dots a^{(1)}_{n-1}
\\
&=a^{(1)}_1\dots a^{(1)}_n,
\end{align*}
where the second to last equality holds if we inductively assume that the identity holds for $n-1$. Here we used that $\cross_{n-1}\gamma_j=0$ for $j<n-1$. 
Therefore we define
\[
\overline{\varphi}(a^{(1)}_1\dots a^{(1)}_n):=(\cross_{n}\varphi)(a_1,\dots, a_n).
\]
Then we observe that if an additive extension $\overline{\varphi}\colon \Gamma_n(A) \to B_{(p)}$ exists, then it will factor over $\Gamma_n(A)_{(p)}$.  A generic generator $a_{1}^{(n_1)}\dots a_{l}^{(n_l)}$ of $\Gamma_n (A)$ with $\sum n_i=n$ will decompose in $\Gamma_n(A)_{(p)}$ as
\[
a_{1}^{(n_1)}\dots a_{l}^{(n_l)}=\frac{1}{n_1!\dots n_l!}\underbrace{a_{1}^{(1)}\dots a^{(1)}_1}_{n_1}\dots \underbrace{a_{l}^{(1)}\dots a^{(1)}_l}_{n_l},
\]
where the positive integers smaller or equal to $n$ are invertible since $p>n$. We therefore define
\[
\overline{\varphi}(a_{1}^{(n_1)}\dots a_{l}^{(n_l)}):=\frac{1}{n_1!\dots n_l!}(\cross_{n}\varphi)(\underbrace{a_{1},\dots, a_1}_{n_1},\dots, \underbrace{a_{l},\dots, a_l}_{n_l}).
\]
This map extends $\varphi$, since 
\[
\overline{\varphi}(a^{(n)})=\frac{1}{n!}(\cross_{n}\varphi)(a,\dots, a)=\varphi(a),
\]
where the second equality holds because $\varphi$ is $n$-homogeneous.
It remains to verify that this map respects the relations of the divided power algebra. We recall that the $n$-th cross-effect of an $n$-polynomial map is additive in each variable. It follows that
\begin{align*}
\overline{\varphi}((ka)^{(n)})&=\frac{1}{n!}(\cross_{n}\varphi)(ka,\dots, ka)=\frac{k^n}{n!}(\cross_{n}\varphi)(a,\dots, a)=k^n\overline{\varphi}(a^{(n)}),
\\
\overline{\varphi}((a+b)^{(n)})&=\frac{1}{n!}(\cross_{n}\varphi)(a+b,\dots, a+b)
=\frac{1}{n!}\sum_{i+j=n}{n\choose i} (\cross_{n}\varphi)(\underbrace{a,\dots, a}_i,\underbrace{b,\dots ,b}_j)
\\
&=\sum_{i+j=n}\frac{1}{i!j!}(\cross_{n}\varphi)(\underbrace{a,\dots, a}_{i}, \underbrace{b,\dots, b}_{j})
 =\sum_{i+j=n} \overline{\varphi}(a^{(i)}b^{(j)}),
\\
\overline{\varphi}(a^{(i)}a^{(j)})&=\frac{1}{i!j!}(\cross_{n}\varphi)(a,\dots,a)
={i+j\choose i}
\frac{1}{(i+j)!}(\cross_{n}\varphi)(a,\dots,a)
=
{i+j\choose i}\overline{\varphi}(a^{(i+j)}).
\end{align*}
The additive map $\overline{\varphi}$ is unique since it factors through an additive map
\[\widehat{\overline{\varphi}} \colon  \Gamma_n(A)_{(p)} \longrightarrow B_{(p)} \]
and by the above observations the $\mathbb{Z}_{(p)}$-module $\Gamma_n(A)_{(p)}$ is generated by the image of the canonical map
\[\xymatrix{A \ar[r]^-{\gamma_n} &  \Gamma_n(A) \ar[r] & \Gamma_n(A)_{(p)}. }\]

Next we check that if $\varphi$ is multiplicative, then so is $\overline{\varphi}$. For this it suffices to show that $\widehat{\overline{\varphi}}$ is multiplicative. We have
\[\widehat{\overline{\varphi}}(a^{(n)})=\varphi(a)\]
and the elements of the form $a^{(n)}$ generate the $\mathbb{Z}_{(p)}$-module $\Gamma_n(A)_{(p)}$. Hence $\widehat{\overline{\varphi}}$ is multiplicative on additive generators and hence in general. 

Now we complete the proof for a general an $n$-polynomial map $f\colon A\to B$. By Proposition \ref{homdecomp} the map $f_{(p)}$ decomposes into an orthogonal sum of homogeneous polynomial maps:
\[f_{(p)}=\varphi_0+\varphi_1+\cdots+\varphi_n.\]
The previous paragraph provides an additive extension
\[ \oplus_i \overline{\varphi_i} \colon \prod_{i=0}^n \Gamma_i(A) \to B_{(p)}\]
of $f_{(p)}$.
By Proposition \ref{homdecomp} and the previous paragraph the latter extension is unique. Further if $f$ is multiplicative, then by Proposition \ref{homdecomp} so are $\varphi_i$-s and the images of $\varphi_i$ and $\varphi_j$ are orthogonal if $i \neq j$. We already saw that $\overline{\varphi_i}$ is multiplicative for every $i$. Similarly, using the fact that $\Gamma_i(A)_{(p)}$ is generated by the elements of the form $a^{(i)}$, we can see the desired orthogonality property. Altogether we get that $\oplus_i \overline{\varphi_i}$ is multiplicative.
\end{proof}

\begin{rem} Any abelian group $A$ admits a universal polynomial map $A\to P_n(A)$ which classifies $n$-polynomial maps out of $A$. The construction can be found in \cite{Passi1, Passi2}. If $A$ is a ring, then the ring multiplication on $A$, makes $P_n(A)$ also into a ring. The latter proposition implies that if $p >n$, then the $p$-localisation of $P_n(A)$ is isomorphic as a ring to $\prod_{i=0}^n \Gamma_i(A)_{(p)}$. In other words, the composite
\[\xymatrix{ A \ar[r]^-{\prod_i \gamma_i}  & \prod_{i=0}^n \Gamma_i(A) \ar[r] & \prod_{i=0}^n \Gamma_i(A)_{(p)}}\]
is the universal multiplicative $n$-polynomial map with a $p$-local target. 
\end{rem}

\begin{rem}\label{invertnfact}
We remark that Proposition \ref{polymapvspolylaw} holds also if we replace the localisation $B_{(p)}$ with $B[\tfrac{1}{n!}]$. Thus by Theorem \ref{polylawfunct}, for any $n\geq 1$, any PD-functor $F$ extends to the category with partially defined composition of $\Z[\tfrac{1}{n!}]$-algebras and multiplicative $n$-polynomial maps. In the next section we will see how for the $p$-typical Witt vectors we can further extend this result integrally.
\end{rem}

\subsection{Functoriality in polynomial maps}\label{secWpoly}

Let us fix once and for all a prime number $p$. Let $1\leq m\leq\infty$ be an integer or infinity. We denote by $W_m(A)$ the ring of $p$-typical $m$-truncated Witt vectors of $A$. The case $m=\infty$ gives the full ring of $p$-typical Witt vectors for which we will use the usual notation $W(A):=W_{\infty}(A)$.

\begin{theorem}\label{Wittpoly}
The functor $W_m \colon \Com \to \Com$ extends to the partial category of multiplicative $(p-1)$-polynomial maps. That is, a multiplicative $n$-polynomial map $f\colon A\to B$, with $n<p$, induces a multiplicative $n$-polynomial map $W_m(f)\colon W_m(A)\to W_m(B)$, with the property that if $f\colon A\to B$ and $g\colon B\to C$ are multiplicative $n$ and $k$-polynomial, respectively, and $nk<p$, then
\[
W_m(g)\circ W_m(f)=W_m(g\circ f).
\]
This extension is unique with the property that the ghost coordinates of the map $W_m(f)$ is the product map $\prod_m f$, i.e., the square
\[
\xymatrix@C=70pt@R=15pt{
W_m(A)\ar[r]^-{W_m(f)}\ar[d]_{w}&W_m(B)\ar[d]^w
\\
\prod_{m} A\ar[r]_-{\prod_{m} f}&\prod_{m} B
}
\]
commutes. If moreover $f,g\colon A\to B$ are multiplicative $n$ and $k$ polynomial, respectively, and $n+k<p$, then $W_m(f\cdot g)=W_m(f)\cdot W_m(g)$.
\end{theorem}

\begin{example}\label{formula} Using the uniqueness of the functoriality of Theorem \ref{Wittpoly}, one can go ahead and try to calculate the components of a polynomial map $W_m(f)$ by inductively solving the equations provided by the description in ghost components.

Let us consider $p$ odd so that $f$ can have degree greater than $1$.
The first component of the image of a Witt vector $(a_0,a_1,\dots)$ in $W(A)$ by a polynomial map $W(f)\colon W(A)\to W(B)$ is $b_0=f(a_0)$. The next component $b_1$ must be the unique natural solution to the equation
\[
f(a_0)^{p}+pb_1=w_1(b_0,b_1)=f(w_1(a_0,a_1))=f(a^{p}_0+pa_1).
\]
Since $p>n$, the map $f$ is also $p$-polynomial
and from the equation $\cross_{p}f(a_0^p+(-))=0$ one can calculate that
\[
f(a^{p}_0+pa_1)=f(a_{0}^{p})+\sum_{i=1}^{p-1}(-1)^i{p\choose i}f(a_{0}^p+ia_1).
\]
The binomial coefficients of this sum are all divisible by $p$, and the unique natural solution to the equation above is
\[
b_1=\sum_{i=1}^{p-1}(-1)^i({p\choose i}/p)f(a_{0}^p+ia_1).
\]
We remark that when $f$ is a ring homomorphism, this sum is in fact equal to $f(a_1)$ so that we indeed recover the usual functoriality in ring homomorphisms.
\end{example}

The proof of Theorem \ref{Wittpoly} will use the following well-known lemma:

\begin{lemma} \label{wittpullback} For any commutative ring $A$ and integer or infinity $1\leq m\leq \infty$, the commutative diagram
\[\xymatrix@C=70pt@R=15pt{W_m(A) \ar[d]_w \ar[r] & W_m(A_{(p)}) \ar[d]^w \\ \prod_{m} A \ar[r] & \prod_{m} A_{(p)}  }\]
is a pullback of rings, where the horizontal maps are induced by the canonical localisation homomorphism $A \to A_{(p)}$.
\end{lemma}

\begin{proof} The case $m=\infty$ follows from the case $m < \infty$ by passing to inverse limits, since pullbacks commute with inverse limits. For finite $m$, we prove the statement by induction. The case $m=1$ is obvious. Next we observe that for any commutative ring (in fact, for any abelian group), the diagram of abelian groups
\[\xymatrix@C=70pt@R=15pt{A \ar[r]^-\pi \ar[d]_{p^m} & A_{(p)} \ar[d]^{p^m} \\ A \ar[r]_-\pi & A_{(p)} } \]
is a pullback. This follows from the fact that the induced maps on kernels and cokernels of the vertical maps are isomorphisms. Using the latter commutative diagram we see that the sequence
\[\xymatrix@C=30pt{0 \ar[r] & A \ar[r]& W_{m+1}(A_{(p)}) \times_{\prod_{m+1} A_{(p)}} \prod_{m+1} A \ar[r]^-{(R,\proj)} & W_{m}(A_{(p)}) \times_{\prod_{m} A_{(p)}} \prod_{m} A \ar[r] & 0}\]
is exact, where the first map sends $a$ to $(V^m(\pi(a)),(0,\dots,0,a))$, $V^m$ is the Verschiebung operator, $R$ is the restriction operator and $\proj$ projects off the last factor. The proof is then completed by induction.
\end{proof}

\begin{proof}[Proof of \ref{Wittpoly}]
%
We denote by $\lambda_R \colon R \to R_{(p)}$ the natural localisation homomorphism.
We first define $W_m(f)\colon W_m(A)\to W_m(B)$ using the pullback of Lemma \ref{wittpullback}. Since pull-backs of rings are pullbacks of sets, it suffices to construct a multiplicative $n$-polynomial map $W_m(A)\to W_m(B_{(p)})$ which will make the diagram of polynomial maps
\[\xymatrix@C=70pt@R=15pt{
W_m(A) \ar[d]_{(\prod_{m} f) \circ w} \ar[r] & W_m(B_{(p)}) \ar[d]^w \\ \prod_{m} B \ar[r]^-{\prod_{m} \lambda_B} & \prod_{m} B_{(p)}  
}\] 
commute. By Proposition \ref{polymapvspolylaw} we know that the composite $A \stackrel{f}{\longrightarrow}B\stackrel{\lambda_B}{\longrightarrow} B_{(p)}$
extends uniquely to a multiplicative polynomial law. We can therefore take the map underlying the multiplicative polynomial law $W_m(\lambda_B \circ f) \colon W_m(A)\to W_m(B_{(p)})$ provided by Theorem \ref{polylawfunct}. Corollary \ref{levelghost} implies that the underlying multiplicative $n$-polynomial map of this polynomial law has the desired description in ghost components. We therefore obtain a map 
\[W_m(f)\colon W_m(A)\longrightarrow W_m(B)\] such that $w \circ W_m(f)= (\prod_{m} f) \circ w$ and $ W_m(\lambda_B) \circ W_m(f)=W_m(\lambda_B \circ f)$. If we evaluate the polynomial law induced by the map $\lambda_B\circ f$
on the ring $\mathbb{Z}_{(p)}$, we get a multiplicative $n$-polynomial map $f_{(p)} \colon A_{(p)} \to B_{(p)}$ making the diagram
\[\xymatrix@C=70pt@R=15pt{A \ar[d]_{\lambda_A} \ar[r]^f & B \ar[d]^{\lambda_B} \\ A_{(p)} \ar[r]^{f_{(p)}} & B_{(p)} }\]
commute. By Proposition \ref{polymapvspolylaw}, the maps $f_{(p)}$ and $\lambda_B \circ f$ are underlying maps of unique polynomial laws. Hence Theorem \ref{functorialityinroby} implies that the diagram
\[\xymatrix@C=70pt@R=15pt{W_m(A) \ar[d]_{W_m(\lambda_A)} \ar[r]^{W_m(f)} & W_m(B) \ar[d]^{W_m(\lambda_B)} \\ W_m(A_{(p)}) \ar[r]^{W_m(f_{(p)})} & W_m(B_{(p)}) }\]
commutes.

Next, we check the identity $W_m(g)\circ W_m(f)=W_m(g\circ f)$. Using the pullback of Lemma \ref{wittpullback}, it suffices to check that $W_m(g)\circ W_m(f)=W_m(g\circ f)$ holds after postcomposing with the canonical map 
\[W_m(\lambda_C) \colon  W_m(C) \to W_m(C_{(p)}).\]
Consider the commutative diagram
\[\xymatrix@C=70pt@R=15pt{A \ar[d]_{\lambda_A} \ar[r]^f & B \ar[d]^{\lambda_B} \ar[r]^g & C \ar[d]^{\lambda_C} \\ A_{(p)} \ar[r]^{f_{(p)}} & B_{(p)} \ar[r]^{g_{(p)}} & C_{(p)}.}\]
The maps $g_{(p)}$, $f_{(p)}$ and $\lambda_{(-)}$ uniquely extend to multiplicative polynomial laws by Proposition \ref{polymapvspolylaw}. So do their composites, and the polynomial laws corresponding to the compositions are the compositions of the polynomial laws associated to the individual maps. Hence Theorem \ref{functorialityinroby} implies that
\[W_m(\lambda_C) \circ W_m(g\circ f)=W_m(\lambda_C \circ g\circ f)= W_m(g_{(p)} \circ f_{(p)} \circ \lambda_A)=W_m(g_{(p)}) \circ W_m(f_{(p)}) \circ W_m(\lambda_A).\]
Finally, using the commutative diagram
\[\xymatrix@C=70pt@R=15pt{W_m(A) \ar[d]_{W_m(\lambda_A)}  \ar[r]^{W_m(f)} & W_m(B) \ar[d]^{W_m(\lambda_B)} \ar[r]^{W_m(g)} & W_m(C) \ar[d]^{W_m(\lambda_C)} \\ W_m(A_{(p)}) \ar[r]^{W_m(f_{(p)})} & W_m(B_{(p)}) \ar[r]^{W_m(g_{(p)})} &  W_m(C_{(p)}),}\]
we conclude  that $W_m(\lambda_C) \circ W_m(g\circ f)=W_m(\lambda_C) \circ W_m(g) \circ W_m(f)$, which shows that $W_m$ is a partial functor.


The uniqueness is immediate  in the torsion-free case since the ghost maps are injective. In the general case one  reduces to the torsion-free case by choosing a resolution similar to the one of Proposition \ref{levelghost}.

Finally, arguing as above and using Proposition \ref{productroby} we see that under our conditions the functor $W_m$ respects multiplications of multiplicative polynomial maps.
\end{proof}

Now we provide examples which show that the conditions of Theorem \ref{Wittpoly} are optimal. 

\begin{example}\label{cex}
The following counterexample shows that the $p-1$ bound on the degree is necessary. Let us consider the multiplicative $p$-polynomial map
\[
N\colon \Z\longrightarrow \Z[x]/(x^2-px)
\]
defined by $N(a)=a+\frac{a^p-a}{p}x$. This is the norm of the Burnside Tambara functor of the cyclic group $C_p$ of order $p$. We show that the norm of the first two ghost components $(Nw_0,Nw_1)$ of $\Z$ is not in the image of the first two ghost coordinates $\langle w_0,w_1 \rangle$ of the ring $\Z[x]/(x^2-px)$.  Indeed, suppose that there are elements $b_0,b_1\in  \Z[x]/(x^2-px)$ such that $\langle w_0,w_1 \rangle(b_0,b_1)=( Nw_0,Nw_1) (0,1)$. Then
\begin{align*}
(b_0,b^{p}_0+pb_1)=( Nw_0,Nw_1 ) (0,1)=(N(0),N(p))=(0,p+(p^{p-1}-1)x),
\end{align*}
and we must have that $pb_1=p+(p^{p-1}-1)x$. This equation has no solution in $\Z[x]/(x^2-px)$ since $p^{p-1}-1$ is not divisible by $p$.

This shows that there cannot be a map on $p$-typical Witt vectors which in ghost components is the map $N$ in each coordinate. That is, for any $m \geq 2$ there cannot be any map (of sets) $f \colon W_m(\Z) \to W_m(\Z[x]/(x^2-px))$ making the diagram
\[\xymatrix@C=60pt@R=17pt{
W_m(\Z)\ar@{-->}[rr]^-f \ar[d]_{w}&&W_m(\Z[x]/(x^2-px))\ar[d]^w
\\
\prod_{m} \Z\ar[rr]_-{\prod_{m} N}&&\prod_{m} \Z[x]/(x^2-px)
}
\]
commute.
\end{example}

\begin{rem}
Another piece of evidence on the optimality of the theorem is provided by the exponentiation maps. The map $(-)^n\colon A\to A$ induces $W((-)^n)=(-)^n\colon W(A)\to W(A)$ for $n<p$. If we try to go beyond the bound $n<p$ we see that the map $(-)^p=\id\colon \F_p\to \F_p$ should simultaneously induce the $p$-th power map and the identity on the $p$-adic integers $W(\F_p)=\Z_p$, contradicting the functoriality of $W$.
\end{rem}

\begin{example}
The functor $W$ does not extend to a functor on the subcategory of commutative rings and set maps generated by the multiplicative $(p-1)$-polynomial maps. The reason is that in general, a map can have several factorisations as compositions of $(p-1)$-polynomial maps, and the extension will depend on this choice. For example, let us take $p=5$  and the multiplicative polynomial maps $(-)^2,(-)^{3}\colon \F_5\to \F_5$, so that $2\cdot 3\geq 5$. These maps compose to the map $(-)^6=(-)^2\colon \F_5\to \F_5$, but the composite of  $(-)^2$ and $(-)^{3}$ is not $(-)^2$ on $W(\F_5) \cong \Z_5$. This shows that the condition $nk<p$ for the composition formula in Theorem \ref{Wittpoly} is necessary.
\end{example}

\begin{rem}
Exponentiation illustrates well the different roles played on Witt vectors by polynomial maps and polynomial laws. The key fact used in the previous examples is that $(-)^p=\id$ as a self polynomial map of $\mathbb{F}_p$. This equality does not however hold as polynomial laws. Indeed, the polynomial law on $\mathbb{F}_p$ defined by exponentiation by $n$ is the natural transformation
\[
(-)^n\colon \mathbb{F}_p\otimes R\longrightarrow \mathbb{F}_p\otimes R,
\]
defined by the exponentiation of the ring $\mathbb{F}_p\otimes R$, where $R$ ranges through all commutative rings. When $R=\Z[t]$ we clearly have that $(-)^p$ on $\mathbb{F}_p\otimes \Z[t]=\mathbb{F}_p[t]$ is not the identity.
\end{rem}


\begin{example}
Similar to the case of composition, for the product of maps one cannot remove the hypothesis that $n+k<p$ and only require $n,k<p$. For example, the map $(-)^6=(-)^2\colon \F_5\to \F_5$ decomposes both as the product of $(-)^3$ and $(-)^3$, and of the identity with itself. However $(-)^6$ and $(-)^2$ are different on $W(\F_5) \cong \Z_5$.
\end{example}

\begin{example}
One cannot extend the functor $W$ additively on sums of multiplicative $(p-1)$-polynomial maps. For example, the identity map $\id\colon \F_p\to \F_p$ decomposes as the sum of $(p+1)$-identities, but $(p+1)W(\id)=(p+1) \id$ is different than $\id$ on $\Z_p$, i.e., such an extension will not be well-defined. 
\end{example}

\begin{rem}
Theorem \ref{Wittpoly} gives a potential obstruction for decomposing an $n$-polynomial map $f\colon A\to B$ into a composition of polynomial maps of lower degree. Indeed, let $n=lk$ with $l,k\neq 1$ and choose a prime $p$ such that $l,k< p\leq  lk$ (say that $l\leq k$, then a prime $p$ with $k<p<2k$ can be used, since $2k\leq lk$). 
By Theorem \ref{Wittpoly}, if $f$ is the composition of two polynomial maps of respective degree $l$ and $k$, it will induce a map on ($m$-truncated) $p$-typical Witt vectors whose ghost coordinates are $\prod f$. In particular  $\prod f\circ w$ is in the image of the ghost coordinates of $B$ and one can attempt to contradict this fact as in Example \ref{cex}.
An example of this will be provided in \S\ref{factor}.
\end{rem}

\section{Applications}

\subsection{The factorisation problem for polynomial maps}\label{factor}

As we pointed out above, Theorem \ref{Wittpoly} can provide an obstruction for detecting if a polynomial map of degree $nk$ is the composite of some polynomial maps of degree $k$ and $n$. In this subsection we give a non-trivial example, where Theorem \ref{Wittpoly} tells us that such a factorisation is not possible.

Let $A_4$ be the fourth alternating group, and $A_3\leq A_4$ any copy of the third alternating group. We denote by $\mathbb{A}(A_n)$ the corresponding Burnside rings.

\begin{prop}\label{exobstruction}
The $4$-polynomial norm map $N_{A_3}^{A_4}\colon \mathbb{A}(A_3)\to \mathbb{A}(A_4)$  does not decompose as the composition of two multiplicative polynomial maps of degree $2$.
\end{prop}

\begin{rem}
The norm $N^{G}_H$ of a Tambara functor, corresponding to a subgroup inclusion $H\leq G$, is always $[G:H]$-polynomial. Suppose that $[G:H]=nk$. If there is a subgroup $H\leq K\leq G$ of index $[G:K]=k$ we would obtain that
\[N^{G}_H=N^{G}_KN^{K}_H\]
decomposes as a composite of an $n$-polynomial map and a $k$-polynomial map. Hence, in order to get an interesting example out of such norms, we need to know that $H$ is maximal in $G$ and has a non-prime index. The group $A_4$ is the smallest group which admits a maximal subgroup $A_3$ of a non-prime index. Hence the norm $N_{A_3}^{A_4}$ does not factor as the composite of two norm maps of subgroups of index $2$. Our theorem shows that it does not even decompose abstractly as the composition of two  $2$-polynomial maps.
\end{rem}

\begin{proof}[Proof of Proposition \ref{exobstruction}] Let $W_2(-;3)$ denote the ring of $2$-truncated $3$-typical Witt vectors.
Suppose that $N:=N_{A_3}^{A_4}$ is the composite 
\[\xymatrix{ \mathbb{A}(A_3) \ar[r]^-f & C \ar[r]^-g & \mathbb{A}(A_4) }\]
of multiplicative polynomial maps $f$ and $g$ of degree $2$, for some commutative ring $C$. Then (since $2<3$) by Theorem \ref{Wittpoly}, we get a commutative diagram
\[\xymatrix@C=70pt@R=15pt{W_2(\mathbb{A}(A_3)) \ar[d]^{w} \ar[r]^{W_2(g) \circ W_2(f)}  &  W_2(\mathbb{A}(A_4)) \ar[d]^{w} \\ \mathbb{A}(A_3) \times \mathbb{A}(A_3) \ar[r]^{N \times N} & \mathbb{A}(A_4) \times \mathbb{A}(A_4). }\]
We claim that the top horizontal map making this diagram commute cannot exist. This will follow immediately if we show that
\[((N \times N) \circ w) (1,1) \]
is not in the image of the ghost map of $\mathbb{A}(A_4)$. The argument is similar to the one in Example \ref{cex}. If
\[((N \times N) \circ w) (1,1)=(N \times N)(1,4)=(1, N(4))\]
was in the image of $w$, then $N(4)$ would be congruent to $1$ mod $3$. 

Let us then show that $N(4)$ is not congruent to $1$ mod $3$. Recall that the norm $N_{H}^G$ in the Burnside Tambara functor is defined as follows. Choose orbit representatives $g_1,\dots g_n$ for the quotient $G/H$ of cardinality $n$. Then each group element $g\in G$ determines a permutation $\sigma\in \Sigma_n$ and elements $h_1,\dots,h_n\in H$, defined by the relations $gg_i=g_{\sigma(i)}h_i$ for $1\leq i\leq n$. Then $N_{H}^G$ sends an $H$-set $X$ to the set $X^{\times n}$ with the $G$ action obtained by restricting the natural $\Sigma_n\wr H$-action along the group homomorphism
\[
G\longrightarrow \Sigma_n\wr H \ \ \ \ \ \ \ \ \ \ \ \ g\longmapsto (\sigma,h_1,\dots,h_n).
\]
In particular an integer $m\in \mathbb{A}(H)$, represented by the trivial $H$-set with $m$-elements $X$, is sent to the set $X^{\times n}$ where $G$ acts by the group homomorphism $\rho\colon G\to \Sigma_n$ that sends $g$ to $\sigma$. In the case of the groups $H=A_3\leq A_4=G$ this is the standard inclusion $A_4\to \Sigma_4$ up to an automorphism of $A_4$ (automorphisms of $A_4$ act trivially on $ \mathbb{A}(A_4)$). Indeed, the kernel of $\rho$ consists of those elements $g\in A_4$ such that $gg_i=g_{i}h_i$ for all $1\leq i\leq 4$, that is
\[
g\in \bigcap_{i=1}^4(g_{i}A_3g_{i}^{-1})=1.
\]
Thus a trivial $A_3$-set $X$ with $m$-elements is sent by $N$ to the set $X^{\times 4}$, where $A_4$ acts by permuting the components via the standard inclusion into $\Sigma_4$.
Up to conjugacy the stabiliser of an element $(x,y,z,w)\in X^{\times 4}$ only depends on how many of the components are equal. This leads to the orbits decomposition of the $A_4$-set $X^{\times 4}$, as
\[
N(m)=X^{\times 4}=m+(m^2-m)A_4/A_{3}+\left(\begin{array}{c}m\\ 2\end{array}\right)A_4/\Z/2+m\left(\begin{array}{c}m-1\\ 2\end{array}\right)A_4/e+2\left(\begin{array}{c}m\\ 4\end{array}\right)A_4/e.
\]
In particular
\[
N(4)=4+12A_4/A_{3}+6A_4/\Z/2+14A_4/e\equiv 1+14A_4/e \ \mbox{mod}\ 3,
\]
with $14$ not divisible by $3$.
\end{proof}

%
%
%
\begin{rem} The above example shows that if a product map lifts to the Witt vectors, then it must satisfy certain congruences. Indeed if $f \colon A \to B$ is a map that fits in a commutative diagram
\[\xymatrix@C=60pt@R=17pt{W_2(A;p) \ar[d]^w \ar[r]^{W_2(f)} & W_2(B;p) \ar[d]^w \\ A \times A \ar[r]^{f \times f} & B \times B}\]
for some map $W_2(f)$, then for any $a_0, a_1 \in A$ we must have
\[f(a_0^p+pa_1) \equiv f(a_0)^p \ \mbox{mod}\ p. \]
Similarly, an analogous argument for $W_m$ with $m > 2$ forces $f$ to satisfy higher versions of these congruences. One can in fact directly show by induction on $k$ that an $n$-polynomial map $f$ satisfies
\[
f(a+p^kc)=f(a)+\sum_{i_1,\dots,i_k=1}^{p-1}(-1)^{i_1+\dots+i_k}{p\choose i_1}\dots{p\choose i_k}f(a+i_1\cdots i_kc)
\]
for every odd prime $p>n$ and $k\geq 0$, where the sum is divisible by $p^k$. This congruence is then preserved by the composition of maps, and our proof of Proposition \ref{exobstruction} shows that $N_{A_3}^{A_4}$ does not satisfy this congruence for $p=3$ and $k=1$.

The formula above is in fact sufficient to lift the product map of $f\colon A\to B$ to the Witt vectors when $A$ and $B$ are torsion-free with compatible Frobenius lifts, as this guarantees that the congruences of the Dwork Lemma are satisfied (this is for example the case for the universal $n$-polynomial map of  \cite{Passi1, Passi2}). Thus these congruences are closely related to the lifts of the product map on the Witt-vectors.
%
%
%
\end{rem}

\subsection{Witt vectors of  \texorpdfstring{$\Z/2$}{Z/2}-Tambara functors
}\label{secTambara}

We use Theorem \ref{Wittpoly} to extend the $p$-typical Witt vectors functor to the category of $\Z/2$-Tambara functors, for odd primes $p$.
We will schematically display a $\Z/2$-Tambara functors $T$ as
\[T=\xymatrix@C=70pt{\big(A \ar@<1ex>[r]^-{\tran}\ar@<-1ex>[r]_-{N}& B\ar[l]|-{\res}\big),}\]
where we keep in mind, but suppress from the notation, that $A$ has an involution which is part of the structure. 
Given a $\Z/2$-Tambara functor $T$ as above, we let
\[\prod_{n}T=\xymatrix@C=70pt{\big(\prod_{n}A \ar@<1ex>[r]^-{\prod\tran}\ar@<-1ex>[r]_-{\prod N}& \prod_{n}B\ar[l]|-{\prod\res}\big)}\]
be the $n$-fold product of $T$ in the category of $\Z/2$-Tambara functors. The ring structures, involution, restriction, transfer and norm are all defined componentwise. The classical Witt polynomials for a prime $p$ define ghost coordinates $w\colon \prod_n A\to\prod_{n}A$ and $w\colon \prod_{n}B\to \prod_{n}B$ which are ring homomorphisms precisely when the sources are endowed with the ring structure of the $p$-typical Witt vectors.

\begin{theorem}\label{WittTambara} Let $p$ be an odd prime and $1\leq n \leq\infty$ an integer or infinity.
There is a unique structure of $\Z/2$-Tambara functor $\xymatrix@C=40pt{W_{n}(A) \ar@<1ex>[r]^-{\tran}\ar@<-1ex>[r]_-{N}&W_{n}(B)\ar[l]|-{\res}}$ such that the Witt polynomials define a natural transformation of Tambara functors
\[
\xymatrix@C=70pt@R=17pt{W_{n}(A)\ar[d]_w \ar@<1ex>[r]^-{\tran}\ar@<-1ex>[r]_-{N}&W_{n}(B)\ar[d]^w\ar[l]|-{\res}
\\
\prod_{n}A \ar@<1ex>[r]^-{\prod\tran}\ar@<-1ex>[r]_-{\prod N}& \prod_{n}B\ar[l]|-{\prod\res}\rlap{\ .}
}
\]
We denote the resulting Tambara functor by $W_{n}(T)$.
\end{theorem}

\begin{rem}
Theorem \ref{WittTambara} cannot be extended to the prime $p=2$. Indeed, Example \ref{cex} shows that the norm of the $\Z/2$-Burnside Tambara functor does not induce a map on $2$-typical Witt vectors with the above description in ghost components.
\end{rem}

We start by defining the maps of the Tambara functor $W_n(T)$. The restriction map $\res\colon B\to A$ is a ring homomorphism, and therefore it induces a ring homomorphism $W_n(\res)\colon W_n(B)\to W_n(A)$, and we define this to be the restriction of $W_n(T)$. Similarly, the involution on $A$ induces an involution on  $W_n(A)$. The multiplicative transfer $N\colon A\to B$ of $T$ is multiplicative $2$-polynomial, and Theorem \ref{Wittpoly} provides an induced multiplicative $2$-polynomial map for odd primes
\[
W_n(N)\colon W_{n}(A)\longrightarrow W_{n}(B)
\]
which in ghost components is the product map. We declare this to be the norm of $W_{n}(T)$. The additive transfer $\tran\colon A\to B$ of a $\Z/2$-Tambara functor is always determined by the norm $N$ by the Tambara reciprocity formula
\[
\tran(a)=N(a+1)-N(a)-1.
\]
Therefore we define $W_n(\tran):=W_n(N)\circ (1+\id)+W_n(N)-1$. By Theorem \ref{Wittpoly} the ghost components of these maps are the product maps, and therefore it remains to show that this structure indeed defines a Tambara functor, and that it is unique.

\begin{proof}[Proof of \ref{WittTambara}]
The main tool for proving this theorem is the existence of a resolution of Tambara functors
\[
\xymatrix@C=70pt@R=17pt{
\Z[A]\ar@{->>}[d]_\epsilon \ar@<1ex>[r]^-{\tran}\ar@<-1ex>[r]_-{N}& \ar[l]|-{\res} \mathbb{A}[A;B]\ar@{->>}[d]^\delta
\\
A \ar@<1ex>[r]^-{\tran}\ar@<-1ex>[r]_-{N}& B\ar[l]|-{\res}
}
\]
where the vertical arrows are surjective and where $\mathbb{A}[A;B]$ is torsion free. Such a resolution is constructed in the Appendix by explicitly calculating the left adjoint of the forgetful functor from $\Z/2$-Tambara functors to presheaves of sets on the orbit category of $\Z/2$. The construction of $W_n(T)$ described above gives a diagram
\[
\xymatrix@C=70pt@R=19pt{
W_n(\Z[A])\ar@{->>}[d]_{W_n(\epsilon)} \ar@<1ex>[r]^-{W_n(\tran)}\ar@<-1ex>[r]_-{W_n(N)}& \ar[l]|-{W_n(\res)} W_n(\mathbb{A}[A;B])\ar@{->>}[d]^{W_n(\delta)}
\\
W_n(A) \ar@<1ex>[r]^-{W_n(\tran)}\ar@<-1ex>[r]_-{W_n(N)}& W_n(B)\ar[l]|-{W_n(\res)}
}
\]
where the vertical arrows are surjective ring homomorphisms, which commute with $W_n(N)$, $W_n(\tran)$ and $W_n(\res)$. Since the vertical maps are surjective it is sufficient to show that the relations needed for $W_n(T)$ are satisfied by the upper Tambara functor. Since $\Z[A]$ and $\mathbb{A}[A;B]$ are torsion free their ghost coordinates are injective, where the relations hold since in ghost coordinates the maps become the product maps. The commutativity of this diagrams also shows uniqueness. 
\end{proof}


\subsection{The components of the dihedral fixed-points of  \texorpdfstring{$\THR$}{THR} for odd primes}\label{secTHR}

In this section we apply the main theorem of the previous section to topology. We describe the $\pi_0$ ring of the dihedral fixed-points of real topological Hochschild homology and their multiplicative transfers in terms of the algebraic constructions of the previous section, for odd primes.

We assume that the reader is familiar with basic notions of equivariant stable homotopy theory of \cite{ManMay, Schwede}. We will be mostly interested in the group $O(2)$ and its dihedral subgroups $D_{p^n}$ of order $2p^n$. We fix once and for all a group isomorphism $O(2) \cong S^{1}\rtimes \Z/2$ by sending the generator of $\Z/2$ to the reflection with respect to the $x$-axis.

Let $E$ be an \emph{orthogonal ring spectrum with anti-involution} (see \cite[Section 2.1]{THRmodels} for a definition), which is \emph{flat} (flat here refers to being underlying cofibrant in the flat model structure of \cite{Sto, BrDuSt} on $\Z/2$-spectra). 
We recall from \cite{THRmodels} and \cite{Amalie} that the \emph{real topological Hochschild homology} of $E$ is the $O(2)$-spectrum defined as the geometric realisation of the dihedral nerve
\[
\THR(E):=B^{di}E=|[k]\longmapsto E^{\wedge k+1}|
\]
with the usual cyclic structure of the cyclic nerve, and the involution of $E^{\wedge k+1}$ defined as the indexed smash product over the $\Z/2$-set $\{0,1,\dots,k\}$ with the involution which reverses the order of $\{1,\dots,k\}$ and keeps $0$ fixed. 
We will regard $\THR(E)$ as a genuine $D_{p^n}$-spectrum for all primes $p$ and $n\geq 0$. For every integer $n\geq 0$, we define a $\Z/2$-equivariant spectrum
\[
\TRR^{n+1}(E;p):=\THR(E)^{C_{p^n}},
\]
where $(-)^{C_{p^n}}$ stands for the derived fixed points and $\Z/2$ acts via the Weyl action of $D_{p^n}/C_{p^n}\cong \Z/2$.

\begin{rem}
The $\Z/2$-spectra $\TRR^n(E;p)$ were constructed in \cite{Amalie} using a version of B\"{o}kstedt's model for $\THR$. This model is compared to the dihedral nerve in \cite{THRmodels}, where an explicit zig-zag of equivalences is constructed. The same zig-zag is used in \cite{DMPSW} to see that the cyclotomic structures, reviewed below, agree. By combining these results we see that the two models 
of $\TRR^n(E;p)$ defined here and in \cite{Amalie} are equivalent.
\end{rem}

When $E$ is a $\Z/2$-equivariant flat commutative orthogonal ring spectrum, then $\THR(E)$ is a $D_{p^n}$-equivariant commutative orthogonal ring spectrum. In this case $\TRR^{n+1}(E;p)$ is a commutative $\Z/2$-equivariant ring spectrum, and we are interested in calculating its $\Z/2$-Tambara functor of components. In \cite[Cor. 5.2]{THRmodels} it is shown that for $n=0$
\[
\underline{\pi}_0\THR(E)=\underline{\pi}_0\TRR^{1}(E;p)=\xymatrix@C=40pt{\big(\pi_0E \ar@<1ex>[r]^-{\tran}\ar@<-1ex>[r]_-{N}& \pi_0(E^{\Z/2})\otimes_\phi\pi_0(E^{\Z/2}) \ar[l]|-{\res}\big),}
\]
where $\otimes_\phi$ indicates the quotient of the tensor product by the subgroup (which is in fact an ideal) generated by the elements of the form
\[x \otimes N(e)y-xN(e) \otimes y,\] 
where $N\colon\pi_0E \to \pi_0^{\Z/2}(E)$ is the multiplicative norm. We recall that a $\Z/2$-Tambara functor $T$ is called \emph{cohomological} if $N\res=(-)^2$.

\begin{theorem}\label{pi0TRR} Let $E$ be a connective $\Z/2$-equivariant flat commutative orthogonal ring spectrum, such that $\underline{\pi}_0E$ is cohomological.
For an odd prime $p$ and $n \geq 0$, there is a natural isomorphism of $\Z/2$-Tambara functors
\[\underline{\pi_0}\TRR^{n+1}(E;p)\cong W_{n+1}(\underline{\pi}_0\THR(E)),\]
with the Tambara functor $W_{n+1}(-)$ of $p$-typical Witt vectors of \S\ref{secTambara}. In particular, the ring $\pi_0\THR(E)^{D_{p^n}}$ is isomorphic to $W_{n+1}( \pi_0(E^{\Z/2})\otimes_\phi\pi_0(E^{\Z/2}))$.
\end{theorem}

\begin{example}\
\begin{enumerate}
\item Since the transfer of a $\Z/2$-Tambara functor is determined by the Tambara reciprocity formula $\tran(x)=N(x+1)-N(x)-1$, for cohomological Tambara functors we also have $\tran(1)=2$. 
\item
The Tambara functor associated to a commutative ring with involution $A$ defined by the fixed-points ring $A^{\Z/2}$ is always cohomological, since $N\res (a)=a\overline{a}=a^2$ for all $a\in A^{\Z/2}$.
\item
The Tambara functor
\[
T=\xymatrix@C=40pt{\big(\Z/2 \ar@<1ex>[r]^-{\tran}\ar@<-1ex>[r]_-{N}& \Z/4\ar@{->>}[l]|-{\res}\big)}
\]
is also cohomological, and it is not associated to a ring with involution. The restriction is the canonical projection, $\tran$ sends $0$ to $0$ and $1$ to $2$, and $N$ preserves $0$ and $1$.
\item The Burnside Tambara functor is not cohomological, since $\tran(1)=[\Z/2]$ is the free transitive $\Z/2$-set, which does not represent $2$ in the $\Z/2$-Burnside ring.
\end{enumerate}
\end{example}

An extension of Theorem  \ref{pi0TRR} to non-cohomological $\underline{\pi}_0E$ is provided at the end of the section.
The proof of Theorem \ref{pi0TRR} will use an inductive argument based on maps $R \colon \TRR^{n+1}(E;p) \to \TRR^{n}(E;p)$. These maps are constructed using the \emph{real cyclotomic structure} on $\THR(E)$ which we now recall. We denote by $T$ the $O(2)$-spectrum $T:=\THR(E)$. As $E$ is flat, there is an isomorphism of $O(2)$-spectra
\[\delta\colon \Phi^{C_{p}}T\stackrel{\cong}{\longrightarrow} T,\]
where $\Phi^{C_{p}}$ is the relative monoidal geometric fixed-points functor of \cite[Section 5.4]{ManMay} for a complete $O(2)$-universe (for convenience, we choose the universe $\mathcal{U}$ of \cite[Section 2]{Amalie}). Here $\Phi^{C_{p}}T$ has the residual $O(2)$-action given by the isomorphism
\[O(2)/C_p \cong (S^{1}\rtimes \Z/2)/C_p \cong (S^1/C_p) \rtimes \Z/2 \cong S^1 \rtimes \Z/2 \cong O(2).\]
The map $\delta$ is an $S^1$-equivariant isomorphism by work of \cite{sixauthors} when $E$ is cofibrant as an associative or commutative algebra, and by \cite{DMPSW} when $E$ is underlying flat, based on results of \cite{Sto, BrDuSt}. The map $\delta$ is moreover $\Z/2$-equivariant, and therefore an $O(2)$-equivariant isomorphism when $E$ is an orthogonal ring spectrum with anti-involution which is flat as a $\Z/2$-spectrum.
When $E$ is commutative, the map $\delta$ is an isomorphism of $\Z/2$-equivariant commutative orthogonal ring spectra. This map $\delta$ is the real cyclotomic structure of $\THR(E)$. We also note that $\Phi^{C_{p}}T$ is already derived as a $\Z/2$-spectrum when $E$ is flat. This is a real analog of \cite[Theorem 4.7]{sixauthors}. One first checks that $\Phi^{C_{p}}$ is derived on the levels of the dihedral bar construction using \cite[\S 3.4.3]{Sto} and the cofibrant replacement functor of \cite[Appendix A.1]{THRmodels}. After this one passes to the geometric realisations as in the proof of \cite[Theorem 4.7]{sixauthors}. By iterating this argument we see that in fact $\Phi^{C_{p}}T$ is derived as a $D_{p^{n+1}}/C_p \cong D_{p^n}$-spectrum for any $n \geq 0$.


\begin{defn} Let $E$ be a flat orthogonal ring spectrum with anti-involution and $p$ a prime.
The restriction map of $T=\THR(E)$ is the zig-zag of $O(2)$-equivariant maps
\[
R\colon T^{C_{p^n}}\stackrel{\simeq}{\longrightarrow}(T^{C_p})^{C_{p^{n-1}}}\longrightarrow (\Phi^{C_{p}}(T_f))^{C_{p^{n-1}}}\stackrel{\simeq}{\longleftarrow}  (\Phi^{C_{p}}(T))^{C_{p^{n-1}}}\xrightarrow{\delta^{C_{p^{n-1}}}} T^{C_{p^{n-1}}},
\]
where $(-)_f$ is a functorial fibrant replacement in the model category of $O(2)$-equivariant orthogonal spectra. We define $\TRR(A;p)$ to be the homotopy limit of the diagram
\[
\dots
\xrightarrow{\delta^{C_{p^2}}} T^{C_{p^2}}\to (\Phi^{C_{p}}(T_f))^{C_{p}}\stackrel{\simeq}{\leftarrow}  (\Phi^{C_{p}}(T))^{C_{p}}\xrightarrow{\delta^{C_{p}}}
 T^{C_{p}}\to \Phi^{C_{p}}(T_f)\stackrel{\simeq}{\leftarrow}  \Phi^{C_{p}}T\xrightarrow{\delta} T
\]
in the category of $\Z/2$-spectra.
\end{defn}
The third map in the zig-zag of the definition of $R$ is an equivalence since  $X \to X_f$ is an acyclic cofibration, which are preserved by the relative geometric fixed points functor.
It will be crucial for the proof of Theorem \ref{pi0TRR} to understand the interaction between the map $\pi_0R$ and the norm of the Tambara functor $\underline{\pi_0}\TRR^{n+1}(E;p)$, when $E$ is commutative.

\begin{rem}
When $E$ is commutative the zig-zag defining $R$ can be arranged to take place in the category of $O(2)$-equivariant orthogonal ring spectra. This is achieved by taking $(-)_f$ to be a functorial fibrant replacement functor for the model category of $O(2)$-equivariant orthogonal ring spectra, and use that the acyclic cofibrations of ring spectra are underlying acyclic cofibrations, which are preserved by geometric fixed points. This in particular gives $\TRR(E;p)$ the structure of an $O(2)$-equivariant orthogonal ring spectrum, and the map
\[\underline{\pi}_0 (R) \colon \underline{\pi}_0 \TRR^{n+1}(E;p) \to  \underline{\pi}_0\TRR^{n}(E;p) \]
is a map of commutative $\Z/2$-Green functors. This is however not obviously a map of $\Z/2$-Tambara functors, since it is not clear if one can represent the $R$ map in the category of $\Z/2$-equivariant \emph{commutative} orthogonal ring spectra. The problem is that the cofibrations of $\Z/2$-equivariant commutative orthogonal ring spectra are not cofibrations of underlying $\Z/2$-spectra, 
and hence the third map defining $R$ will not be a weak equivalence in general. One should be able to solve this problem by working with algebras over an equivariant $E_{\infty}$-operad instead of strictly commutative $O(2)$-ring spectra. 
This is however outside of the main scope of this paper and we will instead explicitly show that $\underline{\pi}_0R$ is compatible with the norms, and hence is a map of $\Z/2$-Tambara functors.
\end{rem}

For ease of notation we write $R\colon=\pi_0(R) \colon \pi_0T^{C_{p^n}} \to \pi_0T^{C_{p^{n-1}}}$ and $R_2\colon=\pi_0^{\Z/2}(R) \colon \pi_0 T^{D_{p^n}} \to \pi_0T^{D_{p^{n-1}}}$.

\begin{lemma} \label{Rtambaramap} Let $E$ be a connective $\Z/2$-equivariant flat commutative orthogonal ring spectrum. Then:
\begin{enumerate}
\item For any $n \geq 1$ and $p$ odd, there is a commutative triangle
\[\xymatrix@C=60pt@R=17pt{\pi_0 \THR(E)^{\Z/2} \ar[r]^{N_{\Z/2}^{D_{p^n}} } \ar[dr]_{N_{\Z/2}^{D_{p^{n-1}}}}  &  \pi_0 \THR(E)^{D_{p^n}} \ar[d]^{R_2} \\   &  \pi_0 \THR(E)^{D_{p^{n-1}}}  \rlap{\ .}}\]
\item For any $n \geq 1$ there is a commutative square
\[\xymatrix@C=60pt@R=17pt{ \pi_0  \THR(E)^{C_{p^n}} \ar[d]^{R} \ar[r]^{N_{C_{p^n}}^{D_{p^n}} }  &  \pi_0 \THR(E)^{D_{p^n}} \ar[d]^{R_2}  \\ \pi_0  \THR(E)^{C_{p^{n-1}}} \ar[r]^{N_{C_{p^{n-1}}}^{D_{p^{n-1}}} }  &  \pi_0 \THR(E)^{D_{p^{n-1}}} \rlap{\ .}} \]
Hence $\underline{\pi}_0 (R) \colon \underline{\pi}_0 \TRR^{n+1}(E;p) \to  \underline{\pi}_0 \TRR^{n}(E;p)$ is a map of $\Z/2$-Tambara functors.
\end{enumerate}
\end{lemma}

\begin{proof} For simplicity we will set $T\colon=\THR(E)$. 
We need the following construction for both parts of the lemma. For any orthogonal $G$-spectrum $X$ and a normal subgroup $N$ of $G$, there is natural transformation
\[ \xymatrix{ \phi\colon \pi_0^GX \cong  \pi_0^GX^c \ar[r]^-{\Phi^N} &  \pi_0^{G/N}\Phi^N X^c \ar[r] &  \pi_0^{G/N} \Phi^N X,}\]
where $(-)^c$ is a functorial cofibrant replacement (of equivariant spectra or associative algebras). In other words this construction first takes derived geometric fixed points and then composes with the canonical map of the cofibrant replacement. The maps $R$ and $R_2$ can be described by the composites
\[\xymatrix@C=20pt{R\colon \pi_0^{C_{p^n}}T \ar[r]^-\phi & \pi_0^{C_{p^{n-1}}} \Phi^{C_p}T \ar[r]^-{\delta_\ast} & \pi_0^{C_{p^{n-1}}} T }, \ \ \ \ \ \ \ \ \ \ \ \ \xymatrix@C=20pt{R_2\colon \pi_0^{D_{p^n}}T \ar[r]^-\phi & \pi_0^{D_{p^{n-1}}} \Phi^{C_p}T \ar[r]^-{\delta_\ast} & \pi_0^{D_{p^{n-1}}} T,} \]
respectively. 
We begin with Statement i). As a $D_{p^n}$-spectrum $\THR(E)$ is isomorphic to the geometric realisation of a simplicial $D_{p^n}$-spectrum, defined as the Segal subdivision \cite[Appendix 1]{Segalsub} of the $p^n$-fold edgewise subdivision of the dihedral bar construction. For odd $p$, its zero simplices are isomorphic to the norm $N_{\Z/2}^{D_{p^n}} (E \wedge E)$, and we let
\[v_{p^n} \colon N_{\Z/2}^{D_{p^n}} (E \wedge E) \to \THR(E)\]
be the canonical map from the zero simplices to the geometric realisation. Consider the commutative diagram
\[\hspace{-2.5cm} \xymatrix{\pi_0^{\Z/2} (E \wedge E) \ar[d]^{(v_1)_\ast} \ar[r]^-{N^{ex}} &  \pi_0^{D_{p^n}} N_{\Z/2}^{D_{p^n}} (E \wedge E) \ar[d]^{N_{\Z/2}^{D_{p^n}} (v_1)_\ast} \ar@{=}[r] & \pi_0^{D_{p^n}} N_{\Z/2}^{D_{p^n}} (E \wedge E) \ar[r]^-{\phi} \ar[d]^{(v_{p^n})_\ast} & \pi_0^{D_{p^{n-1}}} \Phi^{C_p} N_{\Z/2}^{D_{p^n}} (E \wedge E) \ar[d]^{\Phi^{C_p} (v_{p^n})_\ast} & \pi_0^{D_{p^{n-1}}} N_{\Z/2}^{D_{p^{n-1}}} (E \wedge E) \ar[l]^-{\Delta_\ast}_-{\cong} \ar[d]^{(v_{p^{n-1}})_\ast}   
\\
 \pi^{\Z/2}_0 T \ar[r]^-{N^{ex}} &  \pi_0^{D_{p^n}} N_{\Z/2}^{D_{p^n}} T \ar[r]^{\epsilon_\ast} &  \pi_0^{D_{p^n}} T \ar[r]^-{\phi} &   \pi_0^{D_{p^{n-1}}} \Phi^{C_p} T \ar[r]^{\delta_\ast} & \pi_0^{D_{p^{n-1}}} T\rlap{\ .}  }\]
Here $N^{ex}$ is the external norm and $\epsilon$ is the counit. The first and third squares commute by naturality. The second square commutes by definition of the subdivisions. The right hand square commutes by the construction of $\delta$. We also note that the diagonal $\Delta$ is an isomorphism by results of \cite{Sto, BrDuSt, sixauthors}. 
By definition the composite $\epsilon_\ast N^{ex}$ is the norm $N_{\Z/2}^{D_{p^n}}$ of $\underline{\pi_0}T$. Moreover the external norm satisfies $\phi N^{ex}=\Delta_\ast N^{ex}$, and therefore
\[
R_2 \circ N_{\Z/2}^{D_{p^n}}\circ(v_1)_\ast=(v_{p^{n-1}})_\ast\Delta^{-1}_\ast\phi N^{ex}=(v_{p^{n-1}})_\ast N^{ex}=N_{\Z/2}^{D_{p^{n-1}}}\circ(v_1)_\ast,
\]
where the last equality uses the first two squares for $n-1$ instead of $n$. Since $(v_1)_\ast$ is surjective this proves Part i).

For Part ii), we consider diagram
\[\xymatrix@C=35pt{\pi_0^{C_{p^n}}T \ar[rr]^-{N^{ex}} \ar[d]^-\phi & &  \pi_0^{D_{p^n}} N_{C_{p^n}}^{D_{p^n}} T \ar[d]^-\phi \ar[r]^{\epsilon_\ast} & \pi_0^{D_{p^n}} T \ar[d]^-\phi 
\\
 \pi_0^{C_{p^{n-1}}} \Phi^{C_p} T \ar[r]^-{N^{ex}} &  \pi_0^{D_{p^{n-1}}}  N_{C_{p^{n-1}}}^{D_{p^{n-1}}}  \Phi^{C_p} T \ar[r]^-{\Delta_\ast}  & \pi_0^{D_{p^{n-1}}} \Phi^{C_p}   N_{C_{p^n}}^{D_{p^n}} T \ar[r]^-{(\Phi^{C_p} \epsilon)_\ast} & \pi_0^{D_{p^{n-1}}}  \Phi^{C_p} T,} \]
where $\Delta \colon N_{C_{p^{n-1}}}^{D_{p^{n-1}}}  \Phi^{C_p} T \to  \Phi^{C_p}  N_{C_{p^n}}^{D_{p^n}} T$ is the relative version of the Hill-Hopkins-Ravenel diagonal constructed in \cite{sixauthors}. We do not claim that it is an equivalence since $T$ is only flat rather than cofibrant. The first square commutes since it commutes after replacing $T$ cofibrantly. The second square commutes by naturality. Moreover the composite $\Phi^{C_p} (\epsilon) \circ \Delta$ is equal to the counit $\epsilon \colon N_{C_{p^{n-1}}}^{D_{p^{n-1}}}  \Phi^{C_p} T \to \Phi^{C_p} T$. This can be seen by explicitly computing the adjoint of the latter composite and identifying it with the identity. Thus the lower horizontal composite is equal to the norm
\[\xymatrix{N_{C_{p^{n-1}}}^{D_{p^{n-1}}} \colon \pi_0^{C_{p^{n-1}}} \Phi^{C_p} T  \ar[r] &  \pi_0^{D_{p^{n-1}}}  \Phi^{C_p} T.}\]
Finally, the claim follows from the fact that the map $\delta \colon \Phi^{C_p} T \to T$ is a map of commutative $O(2)$-ring spectra, and therefore $\delta_\ast$ is compatible with the norms, and from the observation that $R_2$ and $R$ are the composites of $\delta_\ast$ and $\phi$. 
\end{proof}

\begin{proof} [Proof of Theorem \ref{pi0TRR}]  We start by calculating the components of the fixed-points
\[
\pi_0(\TRR^{n+1}(E;p)^{\Z/2})=\pi_0((\THR(E)^{C_{p^n}})^{\Z/2})=\pi_0(\THR(E)^{D_{p^n}}),
\]
using an argument analogous to \cite[\S 3.3]{Wittvect}.
Let us denote the components of the underlying ring spectrum by $A:=\pi_0(E)$ and of the derived fixed-points by $B:=\pi_0(E^{\Z/2})$.
Let $\mathcal{R}$ be the family of subgroups of $D_{p^n}$ generated by the reflections, together with the trivial group. Let $E\mathcal{R}$ be a universal space for this family, for concreteness one could take the unit sphere $E\mathcal{R}=S(\mathbb{C}^{\infty})$ in the countably infinite direct sum of copies of $\mathbb{C}$, where $O(2)$ acts on $\mathbb{C}\cong \mathbb{R}^2$ by the standard action. Since $\mathcal{R}$  is the family of subgroups of $O(2)$ which do not contain $C_p$, using the Adams isomorphism (see e.g., \cite{RV16}) and isotropy separation, we get a cofiber sequence of $D_{p^{n-1}}$-spectra 
\[
E\mathcal{R}_+\wedge_{C_p} \THR(E)\longrightarrow \THR(E)^{C_p}\longrightarrow \Phi^{C_p}\THR(E).
\]
(see \cite{Amalie} where this sequence for B\"okstedt's model is used). By postcomposing the second map with the equivalence $\delta \colon \Phi^{C_p}\THR(E)\stackrel{\simeq}{\to}\THR(E)$ and taking derived $C_{p^{n-1}}$-fixed points (and again using the Adams isomorphism), we obtain a fibre sequence of $\Z/2$-spectra
\[
\xymatrix{E\mathcal{R}_+\wedge_{C_{p^{n}}} \THR(E)\ar[r] & \THR(E)^{C_{p^{n}}}\ar[r]^-{R} & \THR(E)^{C_{p^{n-1}}}.}
\]
By the homotopy orbits spectral sequence induced by the standard filtration of $S(\mathbb{C}^{\infty})$, it is clear that 
\[\pi_0(E\mathcal{R}_+\wedge_{C_{p^{n}}} \THR(E))=\pi_0(\THR(E))_{C_{p^n}} \cong \pi_0(E)=A,\]
where the $C_{p^{n}}$ acts trivially on $A\cong \pi_0\THH(E)$ since the cyclic actions are restricted from the circle. A similar analysis on the spectral sequence converging to the homotopy groups of $(E\mathcal{R}_+\wedge_{C_{p^{n}}} \THR(E))^{\Z/2}$ shows that there is an isomorphism
\begin{align*}
\pi_0(E\mathcal{R}_+\wedge_{C_{p^{n}}} \THR(E))^{\Z/2}&\cong H^{\Z/2}_0(E\mathcal{R}/C_{p^n};\underline{\pi}_0\THR(E))\cong H^{D_{p^n}}_0(E\mathcal{R};\underline{\pi}_0\THR(E))
\\
&\cong \colim_{\mathcal{O}_R} \underline{\pi}^{(-)}_0\THR(E)\cong \pi_0\THR(E)^{\Z/2}
\end{align*}
where the second isomorphism holds because $C_{p^n}$ acts freely on $E\mathcal{R}$. The Bredon homology group $H^{D_{p^n}}_0$ is computed as the colimit over the full subcategory $\mathcal{O}_R$ of the orbit category of $D_{p^{n}}$ generated by $\mathcal{R}$. The final isomorphism holds because when $p$ is odd there is only one conjugacy class of reflections in $D_{p^n}$. We recall that $\pi_0\THR(E)^{\Z/2}$ is a quotient of $B\otimes B$, which we denote by $B\otimes_{\phi} B$.

On homotopy groups the above fibre sequence induces a long exact sequence
\[
\xymatrix@C=15pt{
\dots \ar[r]&
\pi_1\THR(E)^{D_{p^{n-1}}}\ar[r]^-{\partial}&B\otimes_{\phi} B\ar[r]^-{V_{2}^n}&\pi_0\THR(E)^{D_{p^n}}\ar[r]^-{R_2}&\pi_0\THR(E)^{D_{p^{n-1}}}\ar[r]&0
}
\]
where the map $V^{n}_2=\tran_{\Z/2}^{D_{p^{n}}}$ is the transfer map from $\Z/2$ to $D_{p^{n}}$. 
We claim that $V^{n}_2$ is injective, and therefore that the $\pi_0$ terms form a short exact sequence for every $n$. The restriction for the subgroup inclusion $\Z/2\subset D_{p^n}$ defines a map $F_{2}^n\colon \pi_0\THR(E)^{D_{p^n}}\to B\otimes_{\phi} B$. By the double coset formula we see that $F_{2}^nV_{2}^n$ is the map
\[
F_{2}^nV_{2}^n=\sum_{[g]\in {_{\Z/2}}\backslash D_{p^n}/_{\Z/2}}\tran _{{}^g\Z/2\cap \Z/2}^{\Z/2}c_g\res^{\Z/2}_{\Z/2\cap \Z/2^g}=\id+\frac{(p^n-1)}{2}\tran _{1}^{\Z/2}\res^{\Z/2}_{1}
\]
where the conjugations are trivial since $C_{p^n}$ acts trivially on $\pi_\ast\THH(E)$. It follows that any element $x$ in the kernel of $V_{2}^n$ must satisfy
\[
x+\frac{(p^n-1)}{2}\tran _{1}^{\Z/2}\res^{\Z/2}_{1}(x)=0.
\]
Now let us consider the commutative square
\[
\xymatrix@C=60pt{
B\otimes_{\phi} B\ar[d]_{\res^{\Z/2}_{1}}\ar[r]^-{V_{2}^n}&\pi_0\THR(E)^{D_{p^n}}\ar@<-1ex>[d]^{\res^{D_{p^n}}_{C_{p^n}}}
\\
A\ar@{>->}[r]^-{V^n}&\pi_0\THR(E)^{C_{p^n}}
}
\]
where the bottom horizontal map $V^n=\tran^{C_{p^n}}_1$ is injective by \cite[Proposition 3.3]{Wittvect}. By the commutativity of this diagram if $x$ is in the kernel of  $V_{2}^n$ we have that $\res^{\Z/2}_{1}(x)=0$, and therefore by the formula above
\[
x=x+\frac{(p^n-1)}{2}\tran _{1}^{\Z/2}\res^{\Z/2}_{1}(x)=0.
\]
This shows that $V_{2}^n$ is injective.

We define maps $I^{n}_2\colon W_{n+1}(B\otimes_{\phi}B)\to \pi_0\THR(E)^{D_{p^n}}$ for every $n\geq 0$, by the formula
\[
I^{n}_2(x_0,\dots, x_n)=\sum_{i=0}^nV^{i}_2N^{p^{n-i}}_2(x_i),
\]
where $N^{p^{n-i}}_2\colon =N^{D_{p^{n-i}}}_{\Z/2}\colon B\otimes_{\phi} B\cong \pi_0\THR(E)^{\Z/2}\to \pi_0\THR(E)^{D_{p^{n-i}}}$ is a short notation for the norm. We claim that the following diagram commutes and that its rows are exact:
\[
\xymatrix@C=40pt{
0\ar[r]&W_{n}(B\otimes_{\phi} B)\ar[d]_{I^{n-1}_2}\ar[r]^-{V}&W_{n+1}(B\otimes_{\phi} B)\ar[r]^-{R^n}\ar[d]^{I_{2}^n}&B\otimes_{\phi} B\ar[d]^{I_{2}^0}_{\cong}\ar[r]&0
\\
0\ar[r]&\pi_0\THR(E)^{D_{p^{n-1}}}\ar[r]_{V_2}&\pi_0\THR(E)^{D_{p^{n}}}\ar[r]_-{R_{2}^n}&\pi_0\THR(E)^{\Z/2}\ar[r]&0\rlap{\ .}
}
\]
The exactness of the top row follows from the definition of the Witt vectors. An inductive diagram chase using the latter exact sequence for $R_2$ shows that the lower row is also exact. The commutativity of the first square is clear. The second square commutes by the definition of $R_2$ and Lemma \ref{Rtambaramap} (i). Further, again using the same lemma, we know that $N^{p^i}_2$ splits $R_2^i$. This implies inductively that the maps $I^{n}_2$ are all bijections, and it remains to show that $I^{n}_2$ is a ring homomorphism.

We follow the strategy of \cite{Wittvect} and define topological version of the ghost coordinates. That is, we consider the commutative diagram
\[
\xymatrix@C=70pt{
W_{n+1}(B\otimes_{\phi} B)\ar[r]^-{w}\ar[d]_-{I^{n}_2}^{\cong}&\prod_{i=0}^nB\otimes_{\phi} B
\\
\pi_0\THR(E)^{D_{p^{n}}}\ar[ur]_{\overline{w}}
}
\]
where $\overline{w}$ is the ring homomorphism with components $\overline{w}_j=R_{2}^{n-j}\res^{D_{p^n}}_{D_{p^{n-j}}}=R_{2}^{n-j}F^{j}_2$. We prove that this diagram commutes in Lemma \ref{topghost} below.
If $B\otimes_{\phi} B$ is $p$-torsion free the ghost map $w$ is injective. Since $I^{n}_2$ is bijective $\overline{w}$ is also injective, and it is sufficient to show that $\overline{w}\circ I^{n}_2=w$ is a ring homomorphism. This is clear by the definition of the Witt vectors.

Now suppose that  $B\otimes_{\phi} B$ possibly has $p$-torsion. In Lemma \ref{cohTambres} we construct a $\Z/2$-set $X$ and a map of cohomological Tambara functors  $S:=\Z[X]\twoheadrightarrow\underline{\pi}_0E$ which is pointwise surjective, where $\Z[X]$ is regarded as a Tambara functor by the involution induced by the functoriality in $X$. Using the Eilenberg-MacLane functor $H$ from \cite{Ullman}, we can form a homotopy pullback of commutative $\Z/2$-ring spectra
\[\xymatrix{\overline{E} \ar[d] \ar[r]^\epsilon & E \ar[d] \\ H \Z[X] \ar@{->>}[r] & H\underline{\pi_0}E.}\]
Then  $\underline{\pi}_0\overline{E}$ is isomorphic to the Tambara functor associated to $\Z[X]$, and $\epsilon$ induces surjection of Tambara functors on $\underline{\pi}_0$. Let $\overline{B}$ denote $\pi_0 \overline{E}^{\Z/2}$ and $\overline{A}$ denote $\pi_0 \overline{E}$ . By Lemma \ref{ptorfree} below $\pi_0\THR( \overline{E})^{\Z/2}=\overline{B}\otimes_\phi\overline{B}$ is torsion-free. 
The induced map $\overline{B}\otimes_\phi\overline{B}\to B\otimes_{\phi} B$ is also surjective, and we have a commutative diagram
\[
\xymatrix@C=70pt{
W_{n+1}(\overline{B}\otimes_\phi\overline{B})\ar[r]^-{I^{n}_2}\ar@{->>}[d]&\pi_0\THR(\overline{E})^{D_{p^n}}\ar[d]
\\
W_{n+1}(B\otimes_{\phi} B)\ar[r]_-{I^{n}_2}&\pi_0\THR(E)^{D_{p^n}}
}
\]
where the vertical maps and the top horizontal map are ring homomorphisms, and where the left vertical map is surjective. It follows that the bottom horizontal map is also a ring homomorphism.

Let us now identify the $\Z/2$-Tambara structure.
Since the restriction map $\res$ and the involution $w\colon A\to A$ are ring homomorphisms, their induced maps on Witt vectors $W(\res)$ and $W(w)$ are defined coordinatewise. We can therefore verify by direct calculation that the squares
\[
\xymatrix@R=20pt{
W_{n+1}(B\otimes_\phi B)\ar[d]_-{\cong}^-{I_{2}^n}\ar[r]^-{W(\res)}&W_{n+1}(A)\ar[d]^{I_n}_{\cong}
\\
\pi_0\THR(E)^{D_{p^n}}\ar[r]_-{\res_{C_{p^n}}^{D_{p^n}}}&\pi_0\THR(E)^{C_{p^n}}
}
\ \ \ \ \ \ \ \ \
\xymatrix@R=20pt{
W_{n+1}(A)\ar[d]_-{\cong}^-{I^n}\ar[r]^-{W(w)}&W_{n+1}(A)\ar[d]^{I_n}_{\cong}
\\
\pi_0\THR(E)^{C_{p^n}}\ar[r]_-{c_{r}}&\pi_0\THR(E)^{C_{p^n}}
}
\]
commute, where $c_{r}$ is conjugation by the preferred reflection $r=(0,\tau)\in D_{p^n}$, where $\tau$ is the generator of $\Z/2$. The commutativity of the first square is obtained by the double coset formula for the additive and multiplicative transfers
\begin{align*}
\res_{C_{p^n}}^{D_{p^n}} I^{n}_2(a_0,\dots, a_n)&=\res_{C_{p^n}}^{D_{p^n}}\sum_{i=0}^n\tran_{D_{p^{n-i}}}^{D_{p^n}}N^{D_{p^{n-i}}}_{\Z/2}(a_i)=\sum_{i=0}^n\tran_{C_{p^{n-i}}}^{C_{p^n}}\res^{D_{p^{n-i}}}_{C_{p^{n-i}}}N^{D_{p^{n-i}}}_{\Z/2}(a_i)
\\&=\sum_{i=0}^n\tran_{C_{p^{n-i}}}^{C_{p^n}}N^{C_{p^{n-i}}}_{e}\res^{\Z/2}_{e}(a_i)=I^n(\res a_0,\dots, \res a_n),
\end{align*}
where we used that the double cosets $C_{p^n}\backslash D_{p^n}/D_{p^{n-i}}$ and $C_{p^{n-i}}\backslash D_{p^{n-i}}/\Z/2$ are trivial. The second square commutes because  conjugations commute with transfers and norms
\begin{align*}
c_r I^{n}(a_0,\dots, a_n)&=c_r\sum_{i=0}^n\tran_{C_{p^{n-i}}}^{C_{p^n}}N^{C_{p^{n-i}}}_{e}(a_i)=\sum_{i=0}^n\tran_{C_{p^{n-i}}}^{C_{p^n}}c_rN^{C_{p^{n-i}}}_{e}(a_i)
\\&=\sum_{i=0}^n\tran_{C_{p^{n-i}}}^{C_{p^n}}N^{C_{p^{n-i}}}_{e}(c_r a_i)=I^n(c_r a_0,\dots, c_ra_n)=I^n(w(a_0),\dots, w(a_n)),
\end{align*}
where we used that $C_{p^i}$ is normal in $D_{p^n}$.

Since the norm map $W(N)$ is not defined componentwise we are not able to directly show that $W(N)$ and $N_{C_{p^n}}^{D_{p^n}}$ coincide. Instead, 
we show that these agree in ghost components and conclude by reducing to the universal case. We observe that since the transfer is determined by the norm, this will conclude the proof.
We show that the outer part of the diagram
\[
\xymatrix@R=5pt@C=60pt{
W_{n+1}(A)\ar[dd]_-{\cong}^-{I^n}\ar[r]^-{W(N)}&W_{n+1}(B\otimes_\phi B)\ar[dd]_{\cong}^{I_{2}^n}\ar[dr]^{w}
\\
&&\prod_{i=0}^nB\otimes_\phi B
\\
\pi_0\THR(E)^{C_{p^n}}\ar[r]_-{N_{C_{p^n}}^{D_{p^n}}}&\pi_0\THR(E)^{D_{p^n}}\ar[ur]_-{\overline{w}}
}
\]
commutes. By construction (Theorem \ref{WittTambara}), $w\circ W(N)=(\prod N)\circ w$, and the lower composite has components
\begin{align*}
\overline{w}_jN_{C_{p^n}}^{D_{p^n}}I^{n}&=R_{2}^{n-j}\res^{D_{p^n}}_{D_{p^{n-j}}}N_{C_{p^n}}^{D_{p^n}}I^{n}
=\res^{D_{p^j}}_{\Z/2}R_{2}^{n-j}N_{C_{p^n}}^{D_{p^n}}I^{n}=\res^{D_{p^j}}_{\Z/2}N_{C_{p^j}}^{D_{p^j}}R^{n-j}I^n,
\end{align*}
where we use Lemma \ref{Rtambaramap} and that $R_2$ is induced by the fixed-points of the map $R$ of equivariant spectra. By applying the double coset formula for the norm we obtain
\begin{align*}
\res^{D_{p^j}}_{\Z/2}N_{C_{p^j}}^{D_{p^j}}R^{n-j}I^n&=\prod_{g\in \Z/2\backslash D_{p^j}/C_{p^j}} N_{\Z/2 \cap C_{p^j}^g}^{\Z/2}c_g\res^{C_{p^j}}_{\Z/2^g \cap C_{p^j}}R^{n-j}I^n=N_{e}^{\Z/2}\res^{C_{p^j}}_{e}R^{n-j}I^n
\\
&=N_{e}^{\Z/2}\overline{w}_jI^n=N_{e}^{\Z/2}w_j,
\end{align*}
where the last equality is from \cite[Theorem 3.3]{Wittvect}.
Since the triangle in the diagram above commutes by Lemma \ref{topghost} below, this proves the claim when $B\otimes_\phi B$ is $p$-torsion free. In general, the resolution $\epsilon\colon \overline{E}\to E$ above induces a diagram
\[
\xymatrix@C=25pt@R=10pt{
W_{n+1}(\overline{A})\ar@{->>}[dr]^{\epsilon}\ar[dd]_-{\cong}^-{I^n}\ar[rr]^-{W(N)}&&W_{n+1}(\overline{B}\otimes_\phi\overline{B})\ar[dd]_(.3){\cong}^(.3){I_{2}^n }\ar[dr]^{\epsilon}
\\
&W_{n+1}(A)\ar[rr]^(.3){W(N)}\ar[dd]_(.3){I^n}&&W_{n+1}(B\otimes_\phi B)\ar[dd]^{I_{2}^n}
\\
\pi_0\THR(\overline{E})^{C_{p^n}}\ar[dr]^{\epsilon}\ar[rr]_(.65){N_{C_{p^n}}^{D_{p^n}}}&&\pi_0\THR(\overline{E})^{D_{p^n}}\ar[dr]^{\epsilon}
\\
&\pi_0\THR(E)^{C_{p^n}}\ar[rr]_-{N_{C_{p^n}}^{D_{p^n}}}&&\pi_0\THR(E)^{D_{p^n}}\rlap{\ .}
}
\]
The top and bottom faces commute since $\epsilon$ induces a morphism of Tambara functors. The side faces commute by naturality of $I^{n}$ and $I^{n}_2$. The back face commutes by the argument above, since $\overline{B}\otimes_\phi\overline{B}$ is torsion free. Since the maps $\epsilon$ are surjective the front face commutes as well.
\end{proof}

\begin{lemma}\label{ptorfree}
Let $X$ be a $\Z/2$-set. The abelian group $\pi_0\THR(\mathbb{Z}[X])^{\Z/2}$ is $p$-torsion free for every odd prime $p$.
\end{lemma}

\begin{proof}
We observe that $\mathbb{Z}[X]$ is the monoid-ring on the free commutative monoid $M(X)$ generated by the set $X$, with the involution induced functorially by the involution on $X$. It follows from \cite[Proposition 5.12]{THRmodels} that the real topological Hochschild homology spectrum of $\Z[X]$ decomposes as
\[
\THR(\mathbb{Z}[X])\simeq \THR(\Z)\wedge \Sigma^{\infty}N^{di}M(X)_+
\]
where $N^{di}$ is the dihedral nerve with respect to the product of spaces. In particular
\[
\underline{\pi}_0\THR(\mathbb{Z}[X])=\underline{\pi}_0\THR(\mathbb{Z})\Box \underline{\pi}_0(\Sigma^{\infty}N^{di}M(X)_+)\cong \underline{\Z}\Box \underline{\pi}_0(\Sigma^{\infty}N^{di}M(X)_+),
\]
where the identification of $\underline{\pi}_0\THR(\mathbb{Z})$ with the constant Mackey functor $\underline{\Z}$  is in \cite[Corollary 5.2]{THRmodels}. We observe that the underlying group $\pi_0\Sigma^{\infty}N^{di}M(X)_+=\Z[X]$ is torsion-free. Thus the result follows from the following general claim: if $L=(\xymatrix{A\ar@<.5ex>[r]^{\tran}&B\ar@<.5ex>[l]^{\res}})$ is a $\Z/2$-Mackey functor such that $A$ is $p$-torsion-free, then the box product Mackey functor $\underline{\Z}\Box L$ is $p$-torsion free.

Clearly the value at the trivial group $(\underline{\Z}\Box L)(e)=\mathbb{Z}\otimes A\cong A$ is $p$-torsion-free by assumption. The value at $\Z/2$ is the abelian group
\[
(\underline{\Z}\Box L)(\Z/2)=(A\oplus B)/I,
\]
where $I$ is the ideal generated by the elements of the form $2b-\res(b)$, $a-\tran(a)$, $\tau a-a$, for every $a\in A$ and $b\in B$, where $\tau$ is the involution of $A$ (this follows from the Frobenius reciprocity relations for the box product, see e.g., \cite[Section 1.5]{Bouc}). We notice that the second relation collapses the $A$-summand, and that the box product has value isomorphic to
\[
(\underline{\Z}\Box L)(\Z/2)\cong B/J
\]
where $J$ is generated by the elements $2b-\tran\res(b)$. Suppose that $b\in B/J$ is $p$-torsion for some odd prime $p$. The restriction map $\res\colon B/J\to A$ is additive, and since $A$ is $p$-torsion-free we must have $\res(b)=0$. It follows that in $B/J$
\[
0=2b-\tran \res(b)=2b,
\]
that is that $b$ is also $2$-torsion. Since $p$ is odd $b=0$.
\end{proof}

\begin{lemma}\label{topghost}
For every odd prime $p$ and connective $\Z/2$-equivariant flat commutative orthogonal ring spectrum $E$ with $\underline{\pi}_0E$ cohomological, we have $\overline{w}I^{n}_2=w$. 
\end{lemma}

\begin{proof}
First we observe that since $R^{n-j}_2$ is induced by a map of $O(2)$-spectra it commutes with transfers, in the sense that
\[
R^{n-j}_2\tran_{D_{p^{n-i}}}^{D_{p^n}}=\tran_{D_{p^{j-i}}}^{D_{p^{j}}}R^{n-j}_2
\]
if $i\leq j$. Since $R^{n-j}$ is induced by the canonical map to the geometric fixed-points, which kills the proper transfers, one can directly verify that $R^{n-j}_2\tran_{D_{p^{n-i}}}^{D_{p^n}}=0$ if $i>j$.
It also commutes with restrictions, and therefore
\begin{align*}
\overline{w}_jI^{n}_2(a_0,\dots,a_n)&=R_{2}^{n-j}\res^{D_{p^n}}_{D_{p^{n-j}}}\sum_{i=0}^n\tran_{D_{p^{n-i}}}^{D_{p^n}}N^{D_{p^{n-i}}}_{\Z/2}(a_i)
\\
&=\res^{D_{p^{j}}}_{\Z/2}R_{2}^{n-j}\sum_{i=0}^n\tran_{D_{p^{n-i}}}^{D_{p^n}}N^{D_{p^{n-i}}}_{\Z/2}(a_i)
\\
&=\res^{D_{p^{j}}}_{\Z/2}\sum_{i=0}^j\tran_{D_{p^{j-i}}}^{D_{p^{j}}}R_{2}^{n-j}N^{D_{p^{n-i}}}_{\Z/2}(a_i).
\end{align*}
Moreover by Lemma \ref{Rtambaramap} we have that $R_{2}^{n-j}N^{D_{p^{n-i}}}_{\Z/2}=N^{D_{p^{j-i}}}_{\Z/2}$ for $i\leq j$.
It follows that
\begin{align*}
\overline{w}_jI^{n}_2(a_0,\dots,a_n)=\res^{D_{p^{j}}}_{\Z/2}\sum_{i=0}^j\tran_{D_{p^{j-i}}}^{D_{p^{j}}}N^{D_{p^{j-i}}}_{\Z/2}(a_i).
\end{align*}
Now we apply the double coset formula for restrictions and transfers:
\begin{align*}
\overline{w}_jI^{n}_2(a_0,\dots,a_n)&=\sum_{i=0}^j\sum_{[g]\in (\Z/2)\backslash D_{p^j}/D_{p^{j-i}}}\tran^{\Z/2}_{\Z/2\cap (D_{p^{j-i}})^g}c_g\res^{D_{p^{j-i}}}_{(\Z/2)^g\cap D_{p^{j-i}}}N^{D_{p^{j-i}}}_{\Z/2}(a_i).
\end{align*}
The set of double cosets is isomorphic to the quotient $(\Z/2)\backslash C_{p^i}$ of the inversion action, which has representatives $\{1,\theta,\theta^2,\dots, \theta^{(p^i-1)/2}\}$, where $\theta=\sigma^{p^{j-i}}$ for $\sigma \in C_{p^j}$ a generator. The intersection $(\Z/2)^g\cap D_{p^{j-i}}$ is equal to $\Z/2$ if $g=1$, and to the trivial group otherwise (since $p$ is odd). Moreover the conjugations $c_g$ are trivial for the elements of the cyclic group. Therefore we have
\begin{align*}
\overline{w}_jI^{n}_2(a_0,\dots,a_n)&=\sum_{i=0}^j (\res^{D_{p^{j-i}}}_{\Z/2} N^{D_{p^{j-i}}}_{\Z/2}(a_i)+ (p^i-1)/2\tran^{\Z/2}_{e}\res^{D_{p^{j-i}}}_{e}N^{D_{p^{j-i}}}_{\Z/2}(a_i))
\\
&=\sum_{i=0}^j (\res^{D_{p^{j-i}}}_{\Z/2} N^{D_{p^{j-i}}}_{\Z/2}(a_i)+ (p^i-1)/2\tran^{\Z/2}_{e}\res^{\Z/2}_{e}\res^{D_{p^{j-i}}}_{\Z/2}N^{D_{p^{j-i}}}_{\Z/2}(a_i))
\\&=\sum_{i=0}^j p^i\res^{D_{p^{j-i}}}_{\Z/2}N^{D_{p^{j-i}}}_{\Z/2}(a_i).
\end{align*}
where the last equality holds since $\tran^{\Z/2}_{e}\res^{\Z/2}_{e}(a)=\tran^{\Z/2}_{e}(1)\cdot a$, and $\tran^{\Z/2}_{e}(1)=2$ (since $\underline{\pi}_0(E)$ and hence $\underline{\pi}_0\THR(E)$ is cohomological). Similarly, by applying the double coset formula for the norm we have that
\begin{align*}
\res^{D_{p^{j-i}}}_{\Z/2}N^{D_{p^{j-i}}}_{\Z/2}(a_i)&=\prod_{[g]\in (\Z/2)\backslash D_{p^{j-i}}/(\Z/2)}N^{\Z/2}_{\Z/2 \cap (\Z/2)^g}\res^{\Z/2}_{(\Z/2)^g \cap \Z/2}(a_i)
\\&=a_i(N^{\Z/2}_{e}\res^{\Z/2}_{e}(a_i))^{(p^{j-i}-1)/2}.
\end{align*}
Similarly since $\underline{\pi}_0\THR(E)$ is cohomological we have that $N^{\Z/2}_{e}\res^{\Z/2}_{e}(a_i)=a_{i}^2$, and thus
\begin{align*}
\overline{w}_jI^{n}_2(a_0,\dots,a_n)&=\sum_{i=0}^j p^i\res^{D_{p^{j-i}}}_{\Z/2}N^{D_{p^{j-i}}}_{\Z/2}(a_i)=\sum_{i=0}^j p^ia_{i}^{p^{j-i}}=w_j(a_0,\dots,a_n).\qedhere
\end{align*}
\end{proof}

\begin{cor}Let $E$ be a connective $\Z/2$-equivariant flat commutative orthogonal ring spectrum, such that $\underline{\pi}_0E$ is cohomological. Then the Green functor $\underline{\pi_0}\TRR(E;p)$ admits a structure of Tambara functor, and the isomorphisms of Theorem \ref{pi0TRR} induce an isomorphism of Tambara functors
\[\underline{\pi_0}\TRR(E;p)\cong W(\underline{\pi}_0\THR(E))\]
for every odd prime $p$.
\end{cor}

\begin{proof} By the proof of Theorem \ref{pi0TRR} the connecting homomorphism $\partial\colon \pi_1\THR(E)^{D_{p^{n-1}}}\to \pi_0\THR(E)^{\Z/2}$ is zero, and the $R$ maps induce surjective group homomorphisms in $\pi_1$. Thus the Mittag-Leffler condition is satisfied and there are induced ring isomorphisms
\[
\pi_0\TRR(E)^{\Z/2}\cong \lim_n \pi_0\THR(E)^{D_{p^n}}\cong W(B\otimes_\phi B).
\]
Combining this with \cite[Proof of Proposition 3.3]{Wittvect}, yields an isomorphism of Green functors
\[\underline{\pi}_0 \TRR(E) \cong \lim_n \underline{\pi}_0 \TRR^n(E). \]
The right hand side of this isomorphism is canonically a Tambara functor since $\TRR^n(E)$ are $\Z/2$-equivariant commutative ring spectra and the $R$ maps are compatible with the norms by Lemma \ref{Rtambaramap}. Since limits of Tambara and Green functors are computed pointwise, $\underline{\pi}_0 \TRR(E)$ inherits a norm which defines a Tambara functor. The rest follows from Theorem \ref{pi0TRR}.
\end{proof}

Let us now address the case where the flat commutative $\Z/2$-equivariant orthogonal ring spectrum $E$ has a Tambara functor of components $\underline{\pi}_0E$ which is not necessarily cohomological. For any $\Z/2$-Tambara functor $T=(\xymatrix@C=15pt{A \ar@<.5ex>[r]\ar@<-.5ex>[r]&B\ar[l]})$ and odd prime $p$, we define twisted ghost coordinates $\tilde{w}_j\colon \prod_{i=0}^{n}B\to B$ by the formula
\[
\tilde{w}_j(x_0,\dots,x_n):=\sum_{i=0}^j(1+\frac{(p^i-1)}{2}\tran(1))x_i(N\res(x_i))^{\frac{p^{j-i}-1}{2}},
\]
for all $0\leq j<n+1$.
When $T$ is cohomological this is the usual ghost map $w_j$ of the Witt vectors of the commutative ring $B$. If $E$ is a connective commutative  $\Z/2$-equivariant orthogonal ring spectrum, we denote by $A:=\pi_0 E$ and $B:=\pi_0 E^{\Z/2}$. 

\begin{theorem}\label{THRTamb}
Let $E$ be a connective $\Z/2$-equivariant flat commutative orthogonal ring spectrum, and $p$ an odd prime.  There is a unique ring structure $\tilde{W}_{n+1}(B\otimes_\phi B)$ on the set $\prod_{i=0}^{n}B\otimes_\phi B$ such that the maps $\tilde{w}_j$ are natural ring homomorphisms, and a natural ring isomorphism
\[\pi_0\THR(E)^{D_{p^n}}\cong  \tilde{W}_{n+1}(B\otimes_\phi B)\]
for every $1\leq n\leq \infty$.
\end{theorem}

\begin{proof}
The proof of the bijectivity of $I_{2}^n\colon \prod_{i=0}^{n}B\otimes_\phi B\to \pi_0\THR(E)^{D_{p^n}}$ from the proof of Theorem \ref{pi0TRR} does not use that $\underline{\pi}_0E$ is cohomological. Moreover the calculation of Lemma \ref{topghost} shows in fact that the topological ghost coordinates correspond to the twisted algebraic ghost maps, that is $\overline{w}_jI_{2}^n=\tilde{w}_j$. The maps  $\overline{w}_j=R_{2}^{n-j}F^{j}_2$ are natural ring homomorphisms. It is therefore sufficient to show that a ring structure on $\pi_0\THR(E)^{D_{p^n}}$ such that the maps $\overline{w}_j$ are ring homomorphisms is unique.


The product $\overline{w}$ of the maps $\overline{w}_j$ for $T:=\THR(E)$ fits into a commutative diagram
\[
\xymatrix@C=60pt@R=13pt{
\pi_0T^{D_{p^n}}\ar[r]^-{\Psi} \ar[d]_{\overline{w}}&\prod_{i=0}^{n}(\pi_0\Phi^{C_{p^i}}T\times \pi_0\Phi^{D_{p^i}}T)\ar[d]_{\cong}^{\delta}
\\
\prod_{i=0}^{n}\pi_0T^{\Z/2}\ar[r]&\prod_{i=0}^{n}(\pi_0T\times \pi_0\Phi^{\Z/2}T)
}
\]
where $\delta$ is the product of the cyclotomic structure maps. The top horizontal map has components the composites of the restrictions from $D_{p^n}$ to $C_{p^i}$, and from $D_{p^n}$ to $D_{p^i}$, with the canonical projections to the geometric fixed-points. The bottom horizontal map is the product of the map $\pi^{\Z/2}_0T\to \pi_0T\times \pi_0\Phi^{\Z/2}T$ which is the restriction on the first factor and the canonical projection on the second factor. If the top horizontal map $\Psi$ is injective, then the map $\overline{w}$ is injective, and the ring structure on $\pi^{D_{p^n}}_0T$ is unique. We show that $E$ can be resolved by a commutative $\Z/2$-ring spectrum whose map $\Psi$ is injective. 

Given $E$, we resolve the Eilenberg-MacLane ring spectrum of the Tambara functor $\underline{\pi}_0E$ by the free (genuine) $\Z/2$-$E_\infty$-algebra $\mathbb{P}(\Sigma^{\infty} Z_{+})$ of some $\Z/2$-CW complex $Z$, and form the homotopy pullback of commutative $\Z/2$-ring spectra
\[\xymatrix@C=60pt@R=13pt{\overline{E} \ar[r] \ar[d]  & E \ar[d] \\  \mathbb{P}(\Sigma^{\infty} Z_{+}) \ar[r] & H \underline{\pi}_0E. }\]
Here the vertical maps induce isomorphisms on $\underline{\pi}_0$ and the horizontal maps induce surjections. We also note that the free $\Z/2$-$E_{\infty}$-algebra $\mathbb{P}(\Sigma^{\infty} Z_{+})$ can be modelled via a strictly commutative flat $\Z/2$-equivariant ring spectrum and as an associative algebra it is stably equivalent to the spherical monoid-ring $\mathbb{S}\wedge M_+$ of a monoid with anti-involution $M$. The real topological Hochschild homology of the spherical monoid-ring $\THR(\mathbb{S}\wedge M_+)$ is a suspension spectrum by \cite{Amalie} and \cite[5.12]{THRmodels}, and it follows from \cite[3.3.15]{Schwedeglobal} that the map $\Psi$ of $\mathbb{P}(\Sigma^{\infty} Z_{+})$ is injective. By the naturality of $\Psi$ the diagram
\[\xymatrix@C=50pt@R=13pt{\pi_0 \THR(\overline{E})^{D_{p^n}} \ar[d]_-{\cong} \ar[r]^-{\Psi} &  \prod_{i=0}^{n}(\pi_0\Phi^{C_{p^i}}\THR(\overline{E})\times \pi_0\Phi^{D_{p^i}}\THR(\overline{E})) \ar[d] \\ \pi_0 \THR(\mathbb{P}(\Sigma^{\infty} Z_{+}))^{D_{p^n}} \ar@{>->}[r]^-{\Psi} &  \prod_{i=0}^{n}(\pi_0\Phi^{C_{p^i}}\THR(\mathbb{P}(\Sigma^{\infty} Z_{+}))\times \pi_0\Phi^{D_{p^i}}\THR(\mathbb{P}(\Sigma^{\infty} Z_{+}))) }\]
commutes.
The left vertical map is an isomorphism since we already know that $I_2^n$ is a bijection. Hence we conclude that $\Psi$ is injective for $\overline{E}$ as well.
Finally, since $\overline{E} \to E$ induces a surjection on $\underline{\pi}_0$, we get a surjective ring homomorphism $\pi_0\THR(\overline{E})^{\Z/2}\cong \overline{B}\otimes_\phi \overline{B} \twoheadrightarrow B\otimes_\phi B\cong\pi_0\THR(E)^{\Z/2}$, and hence a surjective ring homomorphism
\[
\pi_0\THR(\overline{E})^{D_{p^n}}\cong \prod_{i=0}^{n}\overline{B}\otimes_\phi \overline{B} \longrightarrow \prod_{i=0}^{n}B\otimes_\phi B\cong\pi_0\THR(E)^{D_{p^n}}.
\]
Thus $\pi_0\THR(E)^{D_{p^n}}$ is a quotient of $\pi_0\THR(\overline{E})^{D_{p^n}}$ and the ring structure of the former is determined by the latter.
\end{proof}

\begin{example} \label{different-witt} We demonstrate on an explicit example that the rings  $\tilde{W}_{m+1}(B\otimes_\phi B)$ and $W_{m+1}(B\otimes_\phi B)$ are different in general. Consider the sphere spectrum $\mathbb{S}$ with the trivial $\Z/2$-action. We know that $\THR(\mathbb{S})=\mathbb{S}$ as real cyclotomic spectra and therefore
\[
\pi_0\TRR^{m+1}(\mathbb{S};p)^{\Z/2}=\pi_0\mathbb{S}^{D_{p^m}}\cong \mathbb{A}(D_{p^m})
\]
is the Burnside ring of the dihedral group. Hence $\tilde{W}_{m+1}(B\otimes_\phi B)$ in this case is just $\mathbb{A}(D_{p^m})$. We claim that this ring is not isomorphic to $W_{m+1}(\mathbb{A}(\Z/2))$ for $m>0$. To see this we check that the groups of units $\mathbb{A}(D_{p^m})^{*}$ and $W_{m+1}(\mathbb{A}(\Z/2))^{*}$ are different.

The unit group of the ring $\mathbb{A}(\Z/2)=\Z[x]/(x^2-2x)$ is isomorphic to $\Z/2 \times \Z/2$, containing the elements $\pm 1$ and $\pm (x-1)$, where $x=\tran(1)$. The ring $\mathbb{A}(\Z/2)$ is torsion-free. Hence the ghost map
\[w \colon W_{m+1}(\mathbb{A}(\Z/2)) \to \prod_{m+1} \mathbb{A}(\Z/2)  \]
is injective. Thus $W_{m+1}(\mathbb{A}(\Z/2))^{*}$ injects into the units $(\prod_{m+1} \mathbb{A}(\Z/2))^{*}$. Now using Dwork's Lemma one can see that the image of the latter injection is isomorphic to $\Z/2 \times \Z/2$. Indeed, since $p$ is odd, the identity map of $\mathbb{A}(\Z/2)$ is a Frobenius lift. This implies that a tuple $(x_0, \dots, x_m)$ is in the image of $w$ if and only if $x_{i-1} \equiv x_i \mod p^i$ for all $1 \leq i \leq m$. The latter condition implies that a tuple consisting of the elements of $\{\pm 1, \; \pm (x-1) \}$ is in the image of $w$ if and only if all the elements are equal. We conclude that $W_{m+1}(\mathbb{A}(\Z/2))^{*} \cong \Z/2 \times \Z/2$. We can in fact explicitly check that all the units of $W_{m+1}(\mathbb{A}(\Z/2))$ are given by the Teichm\"uler lifts of the four units of $\mathbb{A}(\Z/2)$:
\[\{\pm 1, \; \pm [x-1] \}.\] 

On the other hand the unit group $\mathbb{A}(D_{p^m})^{*}$ is isomorphic to $(\Z/2)^{m+2}$ (see e.g., \cite{BolPf}). Each subgroup of $C_{p^m}$ gives two units and additionally we also have $\pm 1$. The group $C_{p^m}$ has $(m+1)$ subgroups and in total we get $(\Z/2)^{m+2}$ as the unit group. One can make these elements more explicit when $m=1$. The eight units of $\mathbb{A}(D_{p})$ are given by:
\[\{\pm 1, \; \pm (1+[D_p]-[D_p/C_p]-2[D_p/\Z/2]), \; \pm (1+[D_p]-2[D_p/\Z/2]), \; \pm (1-[D_p/C_p]) \}.\]
\end{example}

\appendix

\section{The free  \texorpdfstring{$\Z/2$}{Z/2}-Tambara functor on a presheaf of sets}

Let $\OO_{\Z/2}$ denote the orbit category of $\Z/2$. We compute the left adjoint of the forgetful functor that sends a Tambara functor to the underlying $\OO_{\Z/2}$-diagram of sets. Let $X\leftarrow Y\colon \res$ be an $\OO_{\Z/2}$-diagram of sets ($X$ is in particular a $\Z/2$-set but we suppress this from the notation for simplicity). 
We let $M(X)$ denote the free multiplicative  abelian monoid generated by a set $X$ and $\Z(M(X))$ its monoid ring, so that the polynomial ring $\Z[X]=\Z(M(X))$. More generally, we denote by $\Z(-)$ the free abelian group functor. We define a Tambara functor
\[
\xymatrix@C=40pt{\mathbb{A}[X;Y]:=\big(
\Z[X] \ar@<1ex>[r]^-{\tran}\ar@<-1ex>[r]_-{N}& \Z(M(Y\amalg X/\Z/2))\oplus \Z(M(X)^{\Z/2})\oplus \Z(M(X)^{free}/\Z/2) \ar[l]|-{\res}
\big),}
\]
where $M(X)^{free}=M(X)\setminus M(X)^{\Z/2}$. Let us first discuss how these different summands arise, and postpone the definition of the structure maps. We represent an element in the first summand as a linear combination of monomials of the form $m(x)g(y)m(\overline{x})$ where $g$ is a monomial in $Y$ and $m$ is a monomial in a set of representatives of the orbits of $X$. An element in the second summand is represented by a linear combination of monomials of the form $m(x)k(x')m(\overline{x})$ where $k$ is a monomial in the fixed-points set $X^{\Z/2}$. An element in the third summand is represented by a linear combination of formal sums $h(x)+h(\overline{x})$ where $h$ is a monomial in $X$ which is not fixed by the involution ($h(\overline{x})+h(x)$ and $h(x)+h(\overline{x})$ are identified).

\begin{rem} \label{remA}
The motivation behind this definition is the following. First one can form the free semi-Mackey-functor on  $X\leftarrow Y\colon \res$ by taking the free commutative monoids on $X$ and $Y$ and freely add a norm. This is the diagram
\[
\xymatrix{
M(X;Y)=\big(M(X) \ar@<-1ex>[r]_-{N}&M(Y\amalg X/{\Z/2}) \ar[l]_-{\res}\big),
}
\]
where the restriction is induced by the restriction on $Y$ and the map $X/{\Z/2}\to X$ that sends $[x]$ to $x\overline{x}$. The norm is induced by the projection $X\to X/{\Z/2}$. The free Tambara functor on this diagram is heuristically $\pi_0$ of the ``spherical monoid ring $\mathbb{S}\wedge M(X;Y)$''. The underlying abelian group of components is then $\Z[X]$. The fixed-points can be additively calculated using the tom Dieck splitting as
\[
\pi_0(\mathbb{S}\wedge (M(X;Y)^{\Z/2}))\oplus \pi_0(\mathbb{S}\wedge M(X;Y))_{h\Z/2}=\Z(M(Y\amalg X/{\Z/2}) )\oplus (\Z[X]/{\Z/2}),
\]
where the quotient in the second summand is a quotient of abelian groups. This is therefore isomorphic to
\[
\Z[X]/{\Z/2}=\Z(M(X)/\Z/2)=\Z(M(X)^{\Z/2})\oplus \Z(M(X)^{free}/\Z/2),
\]
hence the formula. In order to calculate the multiplicative structure from this formula rigorously, one needs to find an actual (or homotopy coherent) topological commutative monoid with involution $M$, such that the associated semi-Mackey functor of components is isomorphic to $M(X;Y)$. We have not investigated if this is possible. Instead of doing this, we directly define the Tambara structure on $\mathbb{A}[X;Y]$ keeping the latter heuristics in mind.
\end{rem}

The commutative multiplication on $\Z(M(Y\amalg X/\Z/2))\oplus \Z(M(X)^{\Z/2})\oplus \Z(M(X)^{free}/\Z/2)$ is defined on additive generators by the following multiplication table:
\begin{center}
\begin{tabular}{c|c|c|c}
$\bullet$&$mg\overline{m}$ &$mk\overline{m}$&$h+\overline{h}$
\\
\hline
$m'g'\overline{m}'$ & $mm'gg'\overline{mm'}$ &$mm'\res(g')k\overline{mm'}$ & $m'\res(g')\overline{m}'h+\overline{m'\res(g')\overline{m}'h}$
\\
\hline
$m'k'\overline{m}'$& & $2mm'kk'\overline{mm'}$&$2(m'k'h\overline{m}'+\overline{m'k'h\overline{m}'})$
\\
\hline
$h'+\overline{h}'$&&&
$\begin{array}{ll}
(hh'+\overline{hh'})+(h\overline{h'}+\overline{h\overline{h'}})& \mbox{if }h\overline{h'}\neq \overline{h\overline{h'}}, \ hh' \neq \overline{hh'}
\\
hh'+(h\overline{h'}+\overline{h\overline{h'}})&\mbox{if }h\overline{h'} \neq \overline{h\overline{h'}}, \ hh'= \overline{hh'}
\\
(hh'+\overline{hh'})+h\overline{h'}& \mbox{if } h\overline{h'}= \overline{h\overline{h'}}, \ hh' \neq \overline{hh'}
\end{array}$
\end{tabular}
\end{center}
where the elements $hh'$ and $h\overline{h'}$ in the last two cases of the bottom right case belong to the second summand. Thus, we see that $\Z(M(Y\amalg X/\Z/2))$ is a subring and that $ \Z(M(X)^{\Z/2})\oplus \Z(M(X)^{free}/\Z/2)$ is an ideal. The restriction, transfer and norm maps are defined on additive generators by the formulas
\begin{align*}
\res(mg\overline{m})=m\res(g)\overline{m} &&\res(mk\overline{m})=2mk\overline{m}&&\res(h+\overline{h})=h+\overline{h}
\\
&&\tran(k)=k&&\tran(h)=h+\overline{h}
\\
N(m)=m1\overline{m}
\end{align*}
where $k$ is a fixed monomial in $X$, and $h$ is a non-fixed monomial. Here $m$ is any monomial in $X$ and its norm belongs to the first summand. The transfer and restriction are extended additively, and the norm by Tambara reciprocity:
\[N(a+b)=N(a)+N(b)+\tran(a\overline{b}).\]
It is easy to verify that the restriction and the norm are multiplicative, and that this indeed defines a Tambara functor.

\begin{lemma} \label{Tamresfree}
The functor $\mathbb{A}[-;-]$ is left adjoint to the forgetful functor that sends a Tambara functor $T=\xymatrix@C=40pt{\big(A \ar@<1ex>[r]^-{\tran}\ar@<-1ex>[r]_-{N}& B\ar[l]|-{\res}\big)}$ to the $\OO_{\Z/2}$-diagram of sets $A\leftarrow B\colon\res$ with the involution of $A$. 
\end{lemma}

\begin{proof}
Clearly a map of Tambara functors $\mathbb{A}[X;Y]\to T$ induces maps of sets $X\to A$ and $Y\to B$, by restricting the maps of underlying rings respectively to the set of polynomial generators $X$ and along the inclusion $Y\to M(Y\amalg X/\Z/2)$. It is easy to verify that these commute with the restriction and the involution.

 Conversely, let us show that maps of sets $\alpha\colon X\to A$ and $\beta\colon Y\to B$ commuting with the restriction and the involution induce a unique map of Tambara functors
\[
\xymatrix@C=40pt{\mathbb{Z}[ X]\ar[d]_{\alpha_\ast} \ar@<1ex>[r]^-{\tran}\ar@<-1ex>[r]_-{N}& \Z(M(Y\amalg X/\Z/2))\oplus \Z(M(X)^{\Z/2})\oplus \Z(M(X)^{free}/\Z/2) \ar@{-->}[d]^{\beta_\ast}\ar[l]|-{\res}
\\
A \ar@<1ex>[r]^-{\tran}\ar@<-1ex>[r]_{N}& B,\ar[l]|-{\res}
}
\]
whose map of underlying presheaves is again given by $\alpha$ and $\beta$. Clearly $\alpha_\ast$ is the unique map of rings which restricts to $\alpha$, and it is equivariant. Since $\alpha_\ast$ and $\beta_\ast$ must commute with transfers and norms, the map $\beta_\ast$ must satisfy the conditions
\begin{align*}
\beta_\ast(mg\overline{m})=\beta_\ast(m1\overline{m})\beta_\ast(g)=\beta_\ast(N(m))\beta_\ast(g)=N(\alpha_\ast(m))\beta_\ast(g)
\\
\beta_\ast(mk\overline{m})=\beta_\ast(\tran(mk\overline{m}))=\tran (\alpha_\ast(mk\overline{m}))
\\
\beta_\ast(h+\overline{h})=\beta_\ast(\tran(h))=\tran(\alpha_\ast(h))
\end{align*}
(here we abuse the notation and denote the monoid maps induced by $\alpha$ and $\beta$ also by $\alpha_\ast$ and $\beta_*$, respectively). Thus if $\beta_\ast$ exists it must be unique. In order to show existence we need to verify that these formulas indeed define a morphism of Tambara functors. It is immediate that $\alpha_\ast$ and $\beta_\ast$ commute with the transfers and the norms. For the restrictions, we verify that
\begin{align*}
\res\beta_\ast(mg\overline{m})&=\res(N(\alpha_\ast(m))\beta(g))=\res(N(\alpha_\ast(m)))\res(\beta_\ast(g))
\\
&=\alpha_\ast(m)\overline{\alpha_\ast(m)}\alpha_\ast \res(g)=\alpha_\ast \res(mg\overline{m})
\\
\res \beta_\ast(mk\overline{m})&=\res\tran (\alpha_\ast(mk\overline{m}))=2\alpha_\ast(mk\overline{m})=\alpha_\ast\res(mk\overline{m})
\\
\res \beta_\ast(h+\overline{h})&=\res \tran(\alpha_\ast(h))=\alpha_\ast(h)+\overline{\alpha_\ast(h)}=\alpha_\ast (\res (h+\overline{h})).
\end{align*}
It remains to verify that $\beta_\ast$ is multiplicative. We do this directly on generators:
\begin{align*}
\beta_\ast(mg\overline{m}\cdot m'g'\overline{m}')&=\beta_\ast(mm'gg'\overline{mm'})=N(\alpha_\ast(mm'))\beta_\ast(gg')
\\
&=N(\alpha_\ast(m))N(\alpha_\ast(m'))\beta_\ast(g)\beta_\ast(g')=\beta_\ast(mg\overline{m})\beta_\ast(m'g'\overline{m}')
\end{align*}
\begin{align*}
\beta_\ast(mg\overline{m}\cdot m'k\overline{m}')&=\beta_\ast(mm'\res(g)k\overline{mm'})=\tran (\alpha_\ast(mm'\res(g)k\overline{mm'}))
\\
&=
\tran \big(\alpha_\ast(m'k\overline{m'})\alpha_\ast(m)\alpha_\ast(\overline{m})\alpha_\ast \res(g)\big)
\\&
=\tran \big(\alpha_\ast(m'k\overline{m'})\res(N(\alpha_\ast(m))\beta_\ast(g))\big)
\\
&
=N(\alpha_\ast(m))\beta(g)\tran (\alpha_\ast(m'k\overline{m'}))
=\beta_\ast(mg\overline{m})\beta_\ast(m'k\overline{m}')
\end{align*}
\begin{align*}
\beta_\ast(mg\overline{m}\cdot (h+\overline{h}))&=\beta_\ast(m\res(g)\overline{m}h+\overline{m\res(g)\overline{m}h})=\tran (\alpha_\ast(m\res(g)\overline{m}h))
\\
&
=\tran(\alpha_\ast(h)\alpha_\ast(m)\overline{\alpha_\ast(m)}\alpha_\ast\res(g))
=\tran(\alpha_\ast(h)\res (N(\alpha_\ast(m))\beta_\ast(g)))
\\
&=
N(\alpha_\ast(m))\beta_\ast(g)\tran(\alpha_\ast(h))
=\beta_\ast(mg\overline{m})\beta_\ast(h+\overline{h})
\end{align*}
\begin{align*}
\beta_\ast(mk\overline{m}\cdot m'k'\overline{m}')&=\beta_\ast(2mm'kk'\overline{mm'})=2\tran (\alpha_\ast(mm'kk'\overline{mm'}))
\\
&=2\tran (\alpha_\ast(mk\overline{m})\alpha_\ast(m'k'\overline{m'}))
=\tran (\alpha_\ast(mk\overline{m})\res(\tran (\alpha_\ast(m'k'\overline{m'}))))
\\
&=\tran (\alpha_\ast(mk\overline{m}))\tran (\alpha_\ast(m'k'\overline{m'}))=\beta_\ast(mk\overline{m})\beta_\ast(m'k'\overline{m}')
\end{align*}
\begin{align*}
\beta_\ast(mk\overline{m}\cdot (h+\overline{h}))&=2\beta_\ast(mk\overline{m}h+\overline{mk\overline{m}h})=2\tran(\alpha_\ast(mk\overline{m}h))
\\
&=\tran (\alpha_\ast(mk\overline{m})(\alpha_\ast(h)+\overline{\alpha_\ast(h)}))=\tran (\alpha_\ast(mk\overline{m})\res\tran(\alpha_\ast(h)))
\\
&=
\tran (\alpha_\ast(mk\overline{m}))\tran(\alpha_\ast(h))=\beta_\ast(mk\overline{m})\beta_\ast(h+\overline{h})
\end{align*}
\begin{align*}
\beta_\ast((h+\overline{h})\cdot(h'+\overline{h}'))&=\beta_\ast\left\{
\begin{array}{ll}
(hh'+\overline{hh'})+(h\overline{h'}+\overline{h\overline{h'}})& \mbox{if }h\overline{h'}\neq \overline{h\overline{h'}}, \ hh' \neq \overline{hh'}
\\
hh'+(h\overline{h'}+\overline{h\overline{h'}})&\mbox{if }h\overline{h'} \neq \overline{h\overline{h'}}, \ hh'= \overline{hh'}
\\
(hh'+\overline{hh'})+h\overline{h'}& \mbox{if } h\overline{h'}= \overline{h\overline{h'}}, \ hh' \neq \overline{hh'}
\end{array}
\right.
\\
&=
\tran(\alpha_\ast (hh'))+\tran(\alpha_\ast (h\overline{h'}))
=\tran(\alpha_\ast (h)(\alpha_\ast (h')+\alpha_\ast (\overline{h'})))
\\
&=
\tran(\alpha_\ast (h))\tran(\alpha_\ast (h'))=\beta_\ast(h+\overline{h})\beta_\ast(h'+\overline{h}')
\end{align*}
\end{proof}

\begin{rem} \label{toprestamb} There is a way of constructing the left adjoint $\mathbb{A}[-;-]$ using $\Z/2$-equivariant stable homotopy theory. Given an $\OO_{\Z/2}$-diagram of sets  $X\leftarrow Y\colon \res$ (with an involution on $X$), one can functorially construct a $\Z/2$-CW complex $Z$ whose fixed point $\OO_{\Z/2}$-diagram is weakly equivalent to the latter. This uses Elmendorf's Theorem \cite{Elmendorf}. Let $\mathbb{P}$ denote the free (genuine) $\Z/2$-$E_{\infty}$-algebra functor. Then using the adjunction of \cite[Theorem 5.2]{Ullman}, we can see that the components of $\mathbb{P}(\Sigma^{\infty} Z_{+})$ are the free Tambara functor, and therefore
\[
\underline{\pi}_0 (\mathbb{P}(\Sigma^{\infty} Z_{+}))\cong\mathbb{A}[X;Y].
\]
We do not go into the details of this construction, but note that we can deduce from this abstract construction that the free Tambara functor is torsion-free. Indeed, the spectrum $\mathbb{P}(\Sigma^{\infty} Z_{+})$ is equivalent to a wedge of suspension spectra of $\Z/2$-CW complexes. By the tom Dieck splitting, the groups of components of such spectra are in fact free-abelian. 

This strategy is different than the heuristics presented in Remark \ref{remA}. There we first pass to the free semi-Mackey functor and then to the free Tambara functor. Here we first go in the additive direction and create the free (honest) Mackey functor and then the associated free Tambara functor. The functors $\Sigma^{\infty}(-)_{+}$ and $\mathbb{P}$ are just topological analogs of the latter two constructions.
\end{rem}

We can use the construction $\mathbb{A}[X;Y]$ to build resolutions for cohomological Tambara functors. Let
\[T=\xymatrix@C=70pt{\big(A \ar@<1ex>[r]^-{\tran}\ar@<-1ex>[r]_-{N}& B\ar[l]|-{\res}\big),}\]
be a cohomological $\Z/2$-Tambara functor, that is one for which $N\res=(-)^2$ (and in particular from Tambara reciprocity $\tran(1)=2$). Any cohomological Tambara functor can be resolved with a torsion-free commutative ring with involution. 

\begin{lemma} \label{cohTambres} For any cohomological Tambara functor $T$ as above, there exists a polynomial ring with involution $S$ and a surjection of Tambara functors
\[
\xymatrix@C=70pt@R=13pt{S\ar@{->>}[d] \ar@<1ex>[r]^-{\tran}\ar@<-1ex>[r]_-{N}& S^{\Z/2} \ar@{->>}[d] \ar[l]|-{\res}
\\
A \ar@<1ex>[r]^-{\tran}\ar@<-1ex>[r]_{N}& B. \ar[l]|-{\res}
}
\]
\end{lemma}

\begin{proof}
Consider the free $\Z/2$-set $X=\Z/2 \times A$ which consists of two disjoint copies of the underlying set of $A$. The $\Z/2$-action exchanges the two copies of $A$. The obvious $\Z/2$-equivariant counit map $ \alpha \colon X \to A$ gives a map of $\OO_{\Z/2}$-diagram of sets
\[
\xymatrix@C=60pt@R=17pt{X\ar[d]_{\alpha} & \emptyset \ar[d]^{\beta} \ar[l]^-{\res}
\\
A & B,\ar[l]^-{\res}
}
\]
which in turn by the adjunction of Lemma \ref{Tamresfree}  gives a map of Tambara functors $\mathbb{A}[X;\emptyset]\to T$:
\[
\xymatrix@C=40pt@R=17pt{\mathbb{Z}[ X]\ar[d]_{\alpha_\ast} \ar@<1ex>[r]^-{\tran}\ar@<-1ex>[r]_-{N}& \Z(M(X/\Z/2))\oplus \Z(M(X)^{\Z/2})\oplus \Z(M(X)^{free}/\Z/2) \ar[d]^{\beta_\ast}\ar[l]|-{\res}
\\
A \ar@<1ex>[r]^-{\tran}\ar@<-1ex>[r]_{N}& B. \ar[l]|-{\res}
}
\]
The top row is not a cohomological Tambara functor. We make some identifications which will transform the top row into a cohomological Tambara functor, and which are preserved by $\beta_\ast$. For any fixed monomial $m\overline{m} \in M(X)^{\Z/2}$, consider the difference
\[2N(m)-m\overline{m}.\]
These elements generate an ideal $I$ of the upper right corner. It follows from the definitions that $\res(I)=0$. Since the lower row is cohomological, we know that $\tran(1)=2$ which implies that $\beta_*(I)=0$. Hence we get a morphism of cohomological Tambara functors
\[
\xymatrix@C=40pt@R=17pt{\mathbb{Z}[ X]\ar[d]_{\alpha_\ast} \ar@<1ex>[r]^-{\tran}\ar@<-1ex>[r]_-{N}& (\Z(M(X/\Z/2))\oplus \Z(M(X)^{\Z/2})\oplus \Z(M(X)^{free}/\Z/2))/I \ar[d]^{\beta_\ast}\ar[l]|-{\res}
\\
A \ar@<1ex>[r]^-{\tran}\ar@<-1ex>[r]_{N}& B, \ar[l]|-{\res}
}
\]
where we keep the notation $\beta_*$ to denote the the induced map on the quotient. An elementary calculation now shows that there is a ring isomorphism
\[(\Z(M(X/\Z/2))\oplus \Z(M(X)^{\Z/2})\oplus \Z(M(X)^{free}/\Z/2))/I \cong \mathbb{Z}[ X]^{\Z/2}.\] 
In fact the Tambara functor associated to the commutative ring with involution $\mathbb{Z}[ X]$ is isomorphic to the top row of the latter diagram. Hence we get a map of Tambara functors
\[
\xymatrix@C=60pt@R=17pt{\mathbb{Z}[ X]\ar[d]_{\alpha_\ast} \ar@<1ex>[r]^-{\tran}\ar@<-1ex>[r]_-{N}& \mathbb{Z}[ X]^{\Z/2}  \ar[d]^{\beta'_\ast}\ar[l]|-{\res}
\\
A \ar@<1ex>[r]^-{\tran}\ar@<-1ex>[r]_{N}& B, \ar[l]|-{\res}
}
\]
where $\beta'_\ast(m\overline{m})=N(\alpha_*(m))$ and $\beta'_\ast(h+\overline{h})=\tran(\alpha_*(h))$, for any monomial $m$ and a non-fixed monomial $h$. 

This almost proves the desired result, except the map $\beta'_\ast$ is not necessarily surjective. By tensoring with $\Z[B]$, where $B$ has the trivial $\Z/2$ action, we obtain a morphism of Tambara functors
\[
\xymatrix@C=40pt{\mathbb{Z}[ X] \otimes \Z[B] \ar@{->>}[d] \ar@<1ex>[r]^-{\tran}\ar@<-1ex>[r]_-{N}& \mathbb{Z}[ X]^{\Z/2} \otimes \Z[B] \ar@{->>}[d] \ar[l]|-{\res}
\\
A \ar@<1ex>[r]^-{\tran}\ar@<-1ex>[r]_{N}& B. \ar[l]|-{\res}
}
\]
Here we define $N(mb)=N(m)b^2$, where $b$ is a monomial in the elements of $B$. For linear combinations, the norm is extended using Tambara reciprocity. The left vertical map is induced by $\alpha_*$ and the restriction $\res \colon B \to A$. The right vertical map is induced by $\beta'_\ast$ and the identity on $B$. Clearly, both vertical maps are surjective. That these maps are compatible with the norm uses that $N \res=(-)^2$, i.e., that our Tambara functor $T$ is cohomological.

It is now easy to verify that the Tambara functor associated to the commutative ring with involution $S=\Z[X] \otimes \Z[B]\cong \Z[X\amalg B]$, where $B$ has a trivial $\Z/2$ action, is isomorphic to the top row of the latter diagram. \end{proof}

\phantomsection\addcontentsline{toc}{section}{References} 
\bibliographystyle{amsalpha}
\bibliography{bib}

\noindent
\begin{tabular}{l}Emanuele Dotto\\
Mathematical Institute, University of Warwick\\
\textit{e-mail address:} \href{mailto:emanuele.dotto@warwick.ac.uk}{emanuele.dotto@warwick.ac.uk}
\end{tabular}
\vspace{.5cm}
\\
\begin{tabular}{l}
Kristian Jonsson Moi\\
Department of Mathematics, KTH\\
\textit{e-mail address:}  \href{mailto:krjm@kth.se}{krjm@kth.se}
\end{tabular}
\vspace{.5cm}
\\
\begin{tabular}{l}
Irakli Patchkoria\\
Department of Mathematics, University of Aberdeen \\
\textit{e-mail address:} \href{mailto:irakli.patchkoria@abdn.ac.uk}{irakli.patchkoria@abdn.ac.uk}
\end{tabular}

\end{document}